\documentclass[12pt]{article}
\usepackage{amssymb}
\usepackage{amsmath}
\usepackage{amsthm}
\usepackage{comment}
\usepackage{color}
\usepackage{graphicx}
\begin{document}

\def\bV{{\bf V}}
\def\bW{{\bf W}}
\def\bx{{\bf x}}
\def\bX{{\bf X}}
\def\by{{\bf y}}
\def\bY{{\bf Y}}
\def\dd{d}
\def\ee{\varepsilon}
\def\rD{{\rm D}}

\newcommand{\removableFootnote}[1]{}

\newtheorem{theorem}{Theorem}
\newtheorem{conjecture}[theorem]{Conjecture}
\newtheorem{lemma}[theorem]{Lemma}

\title{Stochastic Perturbations of Periodic Orbits with Sliding.}
\author{
D.J.W.~Simpson$^{\dagger}$ and R.~Kuske$^{\ddagger}$\\\\
$^{\dagger}$Institute of Fundamental Sciences\\
Massey University\\
Palmerston North\\
New Zealand\\\\
$^{\ddagger}$Department of Mathematics\\
University of British Columbia\\
Vancouver, BC\\
Canada
}
\maketitle



\begin{abstract}

Vector fields that are discontinuous on codimension-one surfaces are known as Filippov systems
and can have attracting periodic orbits involving segments
that are contained on a discontinuity surface of the vector field.
In this paper we consider the addition of small noise to a general Filippov system
and study the resulting stochastic dynamics near such a periodic orbit.
Since a straight-forward asymptotic expansion in terms of the noise amplitude is
not possible due to the presence of discontinuity surfaces,
in order to quantitatively determine the basic statistical properties of the dynamics,
we treat different parts of the periodic orbit separately.
Dynamics distant from discontinuity surfaces is analyzed by the use of a series expansion of the transitional probability density function.
Stochastically perturbed sliding motion is analyzed through stochastic averaging methods.
The influence of noise on points at which the periodic orbit escapes a
discontinuity surface is determined by zooming into the transition point.
We combine the results to quantitatively determine the effect of
noise on the oscillation time for a three-dimensional canonical model of relay control.
For some parameter values of this model,
small noise induces a significantly large reduction in the average oscillation time.
By interpreting our results geometrically, we are able to identify four features
of the relay control system that contribute to this phenomenon.

\end{abstract}

\section{Introduction}
\label{sec:INTRO}
\setcounter{equation}{0}

Filippov systems are vector fields with codimension-one surfaces, termed switching manifolds,
on which the vector field is discontinuous.
Subsets of switching manifolds at which
the vector field on either side of the manifold points towards the manifold
are known as stable sliding regions.
Whenever a trajectory of the system arrives at a stable sliding region,
future evolution is constrained to the switching manifold until it exits the sliding region.
This evolution is known as sliding motion \cite{DiBu08,Fi88}.
In Filippov models of stick-slip oscillators,
sliding motion corresponds to the sticking phase of the dynamics \cite{OeHi96},
and for relay control, sliding motion models extremely rapid switching \cite{DiBu08,Jo03,ZhMo03}.
So-called sliding-mode controllers specifically utilize sliding motion
to achieve superior control objectives \cite{TaLa12,Su06}.
Stable periodic orbits that involve segments of sliding motion arise in models of
stick-slip oscillators \cite{LuGe06,SzOs08},
relay control \cite{DiJo01,JoRa99,JoBa02,ZhFe10},
and population dynamics \cite{DeGr07,AmOl13,TaLi12}.
The purpose of the present paper is to quantitatively determine
the effects of noise on such periodic orbits.

\begin{figure}[b!]
\begin{center}
\setlength{\unitlength}{1cm}
\begin{picture}(7.2,6)
\put(0,0){\includegraphics[height=6cm]{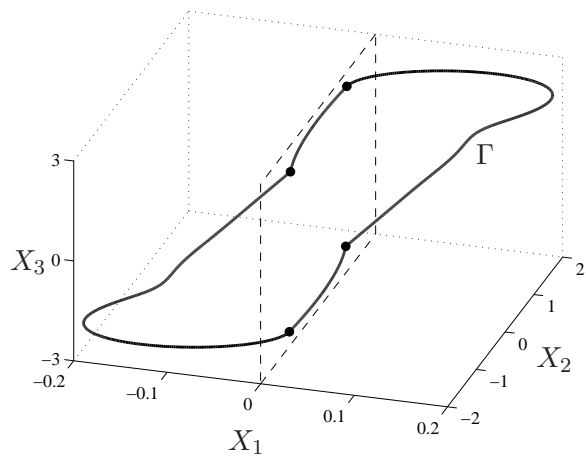}}
\put(6.2,3.9){\small $\Gamma$}
\put(2.9,.1){\small $X_1$}
\put(7,1.2){\small $X_2$}
\put(0,2.5){\small $X_3$}
\end{picture}
\caption{
The attracting periodic orbit $\Gamma$ of (\ref{eq:relayControlSystem})-(\ref{eq:paramValues}).
\label{fig:ppSketch}
}
\end{center}
\end{figure}

Throughout this paper, we use a canonical relay control model, given in \cite{DiBu08,Jo03,ZhMo03}, as an example.
The general model equations are
\begin{equation}
\begin{split}
\dot{\bX} &= A\bX + B\nu \;, \\
\varphi &= C^{\sf T} \bX \;, \\
\nu &= -{\rm sgn}(\varphi) \;,
\end{split}
\label{eq:relayControlSystem}
\end{equation}
where $\bX \in \mathbb{R}^N$ represents the state of the system,
$\varphi$ is the control measurement,
and $\nu$ is the control response.
We consider the following three-dimensional example
($N=3$) of (\ref{eq:relayControlSystem}), given in \cite{DiBu08,DiJo01},
\begin{equation}
A = \left[ \begin{array}{ccc}
-2 \zeta \omega - \lambda & 1 & 0 \\
-2 \zeta \omega \lambda - \omega^2 & 0 & 1 \\
-\lambda \omega^2 & 0 & 0
\end{array} \right] \;, \qquad
B = \left[ \begin{array}{c} 1 \\ -2 \\ 1 \end{array} \right] \;, \qquad
C = \left[ \begin{array}{c} 1 \\ 0 \\ 0 \end{array} \right] \;,
\label{eq:ABCDvalues3d}
\end{equation}
with parameter values,
\begin{equation}
\zeta = 0.5 \;, \qquad \lambda = 0.05 \;, \qquad \omega = 5 \;.
\label{eq:paramValues}
\end{equation}
The system (\ref{eq:relayControlSystem})-(\ref{eq:paramValues})
has an attracting symmetric periodic orbit, call it $\Gamma$, with two sliding segments, Fig.~\ref{fig:ppSketch}\removableFootnote{
When I export a figure that uses fill3 and has axes, to eps, the image is blurry.
}.
The parameter values (\ref{eq:paramValues}) are typical in the sense that 
(\ref{eq:relayControlSystem}) with (\ref{eq:ABCDvalues3d}) exhibits
an attracting periodic orbit with one or more sliding segments over a relatively large range of parameter values \cite{DiJo01}.

From the viewpoint of control,
it is important to understand the robustness of (\ref{eq:relayControlSystem})
to random fluctuations, parameter uncertainty, and unmodelled nonlinear dynamics.
Conditions ensuring the robustness of equilibria of
general hybrid control systems under various assumptions have been established \cite{RaMi10,FeZh06,ChLi04,SkEv99}.
In \cite{DiJo02}, the robustness of attracting periodic orbits of (\ref{eq:relayControlSystem}) is investigated numerically
by altering the switching condition in different ways, such as by incorporating time delay.
The authors conclude that periodic orbits with sliding appear to be less robust
than periodic orbits that only have transversal intersections with the switching manifold.
In \cite{TaOs09}, a model of anti-lock brakes is shown to exhibit attracting periodic orbits with sliding and the
robustness of the periodic orbits is correlated with the size of their basins of attraction.
In addition, unlike attracting periodic orbits of smooth systems,
attracting periodic orbits with sliding segments may be destroyed by stable singular perturbations \cite{SiKo10}\removableFootnote{
Would be nice to include a short paragraph on stochastically perturbed periodic orbits in smooth systems.
However, I actually don't know of any paper that studies this specifically.
One can infer basic results from Theorem 2.3 of Chapter 2 of \cite{FrWe12},
and similar descriptions in \cite{Sc10,GrVa99}.
For coloured noise a stochastic Poincar\'{e} map is derived in \cite{WeKn90}.
Large deviations from a periodic orbit in a smooth system are described in \cite{HiMe13}.
}.

Randomness or uncertainty enters into relay control systems in various ways,
such as via the input and output of the controlling component
or through the action of circuit elements,
and is present in modelling by means of parameter uncertainty
and modelling approximations \cite{Ts84,FrPo02,DoBi01,AsMu08}.
For simplicity, 
we incorporate randomness in (\ref{eq:relayControlSystem})
by adding white Gaussian noise to the control response.
Specifically, the stochastic model is
\begin{equation}
d\bX(t) = \left( A \bX(t) - B \,{\rm sgn} \left( C^{\sf T} \bX(t) \right) \right) \, dt +
\sqrt{\ee} B \, dW(t) \;,
\label{eq:relayControlSystem2}
\end{equation}
where $W(t)$ is standard Brownian motion and $0 < \ee \ll 1$.

Let us consider a sample solution to (\ref{eq:relayControlSystem2}) with (\ref{eq:ABCDvalues3d})-(\ref{eq:paramValues})
from an arbitrary initial point.
Once sufficient time has passed to allow the solution to become close to $\Gamma$,
with high probability the solution follows a random path near $\Gamma$ for a long period of time.
Throughout this paper we ignore transient dynamics.
We define an {\em oscillation time} of a sample solution to (\ref{eq:relayControlSystem2})
as the difference between successive times at which the solution returns to the switching manifold
after a large excursion with $X_1 > 0$. 
The oscillation time represents a stochastic analogue of the period of $\Gamma$.

To investigate the effect of the noise, for a handful of fixed values of $\ee$ we numerically
solved (\ref{eq:relayControlSystem2}) with (\ref{eq:ABCDvalues3d})-(\ref{eq:paramValues})
over a long time frame
and recorded the oscillation times, $t_{\rm osc}$.
For all Monte-Carlo simulations in this paper we used the Euler-Maruyama method with a fixed step size.
We found that different step sizes produced essentially the same results.
We let $t_{{\rm osc},\Gamma}$ denote the period of $\Gamma$, and let
\begin{equation}
{\rm Diff}(t_{\rm osc}) \equiv \mathbb{E}[t_{\rm osc}] - t_{{\rm osc},\Gamma} \;,
\label{eq:Difftosc}
\end{equation}
where $\mathbb{E}$ denotes expectation.
Roughly speaking, we say that the noise alters the oscillation time {\em significantly} if
$|{\rm Diff}(t_{\rm osc})|$ is larger than, or comparable to, ${\rm Std}(t_{\rm osc})$.

\begin{figure}[t!]
\begin{center}
\setlength{\unitlength}{1cm}
\begin{picture}(7.2,6)
\put(0,0){\includegraphics[height=6cm]{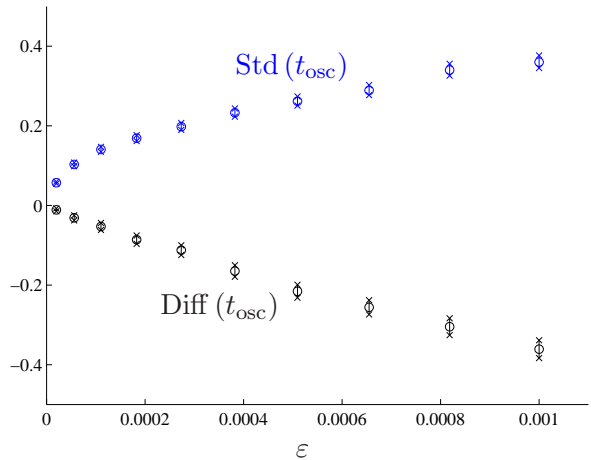}}
\put(3.8,0){\small $\ee$}
\put(2,1.9){\small ${\rm Diff} \left( t_{\rm osc} \right)$}
\put(3,5.1){\small \color{blue} ${\rm Std} \left( t_{\rm osc} \right)$}
\end{picture}
\caption{
Plots of ${\rm Diff}(t_{\rm osc})$
(the difference between $\mathbb{E}[t_{\rm osc}]$ and the period of $\Gamma$ (\ref{eq:Difftosc})),
and ${\rm Std}(t_{\rm osc})$ (the standard deviation of $t_{\rm osc}$)
for the system (\ref{eq:relayControlSystem2}) with (\ref{eq:ABCDvalues3d})-(\ref{eq:paramValues}).
To create this figure, for several different values of $\ee$,
we solved (\ref{eq:relayControlSystem2}) using the
Euler-Maruyama method with a fixed step size, $\Delta t = 0.00001$,
and recorded 1000 consecutive oscillation times, $t_{\rm osc}$.
For each $\ee$, we have plotted the mean value (as a small circle)
and a 95\% confidence interval (as a line segment ending at small crosses)
for both ${\rm Diff}(t_{\rm osc})$ and ${\rm Std}(t_{\rm osc})$ as determined from the 1000 sample values of $t_{\rm osc}$.
\label{fig:manyPeriod}
}
\end{center}
\end{figure}

Fig.~\ref{fig:manyPeriod} shows the variation in ${\rm Diff}(t_{\rm osc})$ and ${\rm Std}(t_{\rm osc})$ with $\ee$,
given by our numerical experiment.
As for an analogous stochastic perturbation of a periodic orbit in a smooth system \cite{FrWe12},
${\rm Diff}(t_{\rm osc}) \sim K_1 \ee$,
and ${\rm Std}(t_{\rm osc}) \sim K_2 \sqrt{\ee}$,
for some constants $K_1, K_2$.
Consequently, as $\ee \to 0$, ${\rm Std}(t_{\rm osc})$ is large relative to ${\rm Diff}(t_{\rm osc})$.
Yet the noise significantly alters the oscillation time for relatively small values of $\ee$ because
$|{\rm Diff}(t_{\rm osc})| \approx {\rm Std}(t_{\rm osc})$ for $\ee = 0.001$.
As we may infer from Fig.~\ref{fig:manyPeriod}, this is because the magnitude of $K_1$ is extremely large.

In our earlier work \cite{SiKu13c},
we gave numerical results similar to Fig.~\ref{fig:manyPeriod} for parameter values different to (\ref{eq:paramValues}).
We analyzed stochastically perturbed sliding motion
and showed that the noise may cause this motion to be significantly faster (or slower) than without noise, on average.
We suggested that this mechanism may be the cause for the reduction in oscillation time.
In this paper we use analytical methods to approximate
${\rm Diff}(t_{\rm osc})$ and ${\rm Std}(t_{\rm osc})$ 
and explain why we may have $|{\rm Diff}(t_{\rm osc})| \approx {\rm Std}(t_{\rm osc})$ for relatively small values of $\ee$.
We find that the mechanism described in \cite{SiKu13c}
is one of four phenomena that induce a significant reduction in oscillation time for (\ref{eq:relayControlSystem2}) with (\ref{eq:ABCDvalues3d})-(\ref{eq:paramValues}).

The remainder of the paper is organized as follows.
As detailed in \S\ref{sec:OUTLINE},
to perform our analysis we split the stochastic dynamics into three phases.
A {\em regular} phase corresponds to dynamics near a section of the periodic orbit that does not involve sliding motion.
This is analyzed in \S\ref{sec:EXCUR} for which the dynamics
are described by a stochastic differential equation with a smooth drift coefficient.
A {\em sliding} phase corresponds to random motion about the switching manifold
near a sliding section of the periodic orbit and is strongly influenced by the discontinuity.
The stochastically perturbed dynamics are analyzed using stochastic averaging principles in \S\ref{sec:SLIDE}.
Our methods for both of these phases are invalid at points where the periodic orbit escapes from the switching manifold,
and the methods are impractical near such points.
Consequently, in \S\ref{sec:ESCAPE} we provide a separate analysis for 
the transition from sliding to regular motion that we refer to as an {\em escaping} phase.
In \S\ref{sec:COMB} we combine the results to determine the statistics of
${\rm Diff}(t_{\rm osc})$ and ${\rm Std}(t_{\rm osc})$
for (\ref{eq:relayControlSystem2}) with (\ref{eq:ABCDvalues3d})-(\ref{eq:paramValues}).
Sections \ref{sec:EXCUR}-\ref{sec:COMB}
involve fundamentally different analytical methods and may be read independently.
Conclusions are given in \S\ref{sec:CONC}.

\section{General equations and three phases of stochastic dynamics}
\label{sec:OUTLINE}
\setcounter{equation}{0}

In this section we begin by introducing general equations and a coordinate system that is most convenient for our analysis, \S\ref{sub:COORDASSUM}.
In \S\ref{sub:DIVISON} we precisely partition the dynamics into regular, sliding and escaping phases.
Lastly in \S\ref{sub:COORDRCS} we construct the coordinate system of \S\ref{sub:COORDASSUM} for the relay control example.

\subsection{A stochastically perturbed Filippov system and assumptions on the equations}
\label{sub:COORDASSUM}

For an $N$-dimensional Filippov system ($N \ge 2$) with a single switching manifold,
we suppose that we may choose our coordinate system such that 
the switching manifold coincides with $x_1 = 0$,
where $x_1 = e_1^{\sf T} \bx$ is the first component of the state variable, $\bx$.
We then write the Filippov system perturbed by noise as
\begin{equation}
d\bx(t) = \left\{ \begin{array}{lc}
\phi^{(L)}(\bx(t)) \;, & x_1(t) < 0 \\
\phi^{(R)}(\bx(t)) \;, & x_1(t) > 0
\end{array} \right\} \,dt + \sqrt{\ee} D \,d\bW(t) \;,
\label{eq:sde}
\end{equation}
where $\phi^{(L)}$ and $\phi^{(R)}$ are functions that are $C^2$ on the closure of their respective half-spaces,
$\bW(t)$ is a standard $N$-dimensional vector Brownian motion,
$0 < \ee \ll 1$ controls the noise amplitude, and
$D$ is an $N \times N$ matrix that specifies the relative strengths and correlations of the noise in different directions.
Throughout this paper it is convenient to separate the component of the noise in the $x_1$-direction
from the remaining directions, and we write
\begin{equation}
D D^{\sf T} = \left[ \begin{array}{c|c} \alpha & \beta^{\sf T} \\ \hline \beta & \gamma \end{array} \right] \;,
\label{eq:alphaBetaGamma}
\end{equation}
where $\alpha \in \mathbb{R}$, $\beta \in \mathbb{R}^{N-1}$ and $\gamma$ is an $(N-1) \times (N-1)$ matrix.
Many generalizations of (\ref{eq:sde}) are possible.
We anticipate that Filippov systems with multiple switching manifolds, nonsmooth switching manifolds,
coloured noise or multiplicative noise,
can be analyzed by extensions of the methods presented below.

We assume that when $\ee = 0$, (\ref{eq:sde}) has an attracting periodic orbit $\Gamma$ that includes at least one sliding segment.
We now reorient the coordinate axes relative to one such sliding segment
in order to analyze the stochastically perturbed dynamics
relating to this segment, and relating to the subsequent segment of $\Gamma$ that does
not intersect the switching manifold.
Our underlying strategy is to repeat this coordinate change and analysis
for all sliding segments in order to determine the overall effect of noise on $\Gamma$.

\begin{figure}[t!]
\begin{center}
\setlength{\unitlength}{1cm}
\begin{picture}(12,9)
\put(0,0){\includegraphics[height=9cm]{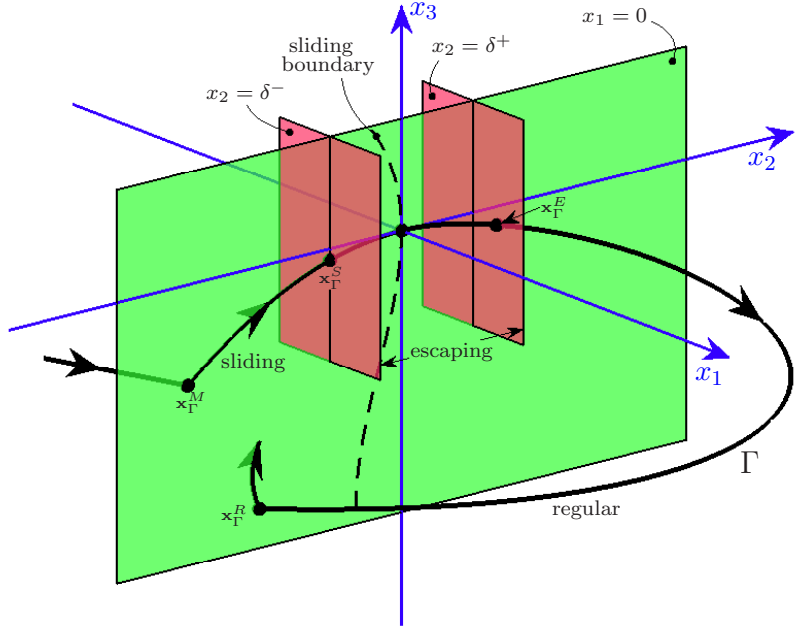}}	
\put(10.5,2.4){$\Gamma$}
\put(9.9,3.64){\small \color{blue} $x_1$}
\put(10.6,6.5){\small \color{blue} $x_2$}
\put(6.1,8.5){\small \color{blue} $x_3$}
\put(3,3.27){\tiny $\bx_{\Gamma}^M$}
\put(4.9,4.93){\tiny $\bx_{\Gamma}^S$}
\put(7.85,5.87){\tiny $\bx_{\Gamma}^E$}
\put(3.6,1.77){\tiny $\bx_{\Gamma}^R$}
\put(3.6,3.8){\scriptsize sliding}
\put(6.12,3.95){\scriptsize escaping}
\put(8,1.84){\scriptsize regular}
\put(3.4,7.38){\scriptsize $x_2 = \delta^-$}
\put(4.5,8){\scriptsize sliding}
\put(4.4,7.75){\scriptsize boundary}
\put(6.4,8){\scriptsize $x_2 = \delta^+$}
\put(8.4,8.4){\scriptsize $x_1 = 0$}
\end{picture}
\caption{
A schematic showing part of a periodic orbit, $\Gamma$, of (\ref{eq:sde}) with $\ee = 0$ that involves sliding
and conforms to assumptions (\ref{eq:a1})-(\ref{eq:a5}).
$\Gamma$ arrives at the switching manifold at a point, $\bx_{\Gamma}^M$,
slides along the switching manifold to the origin, $\bx = 0$,
then travels in the right-half space until returning to the switching manifold at a point, $\bx_{\Gamma}^R$.
To analyze stochastic dynamics,
we divide evolution from $\bx_{\Gamma}^M$ to $\bx_{\Gamma}^R$
into three distinct phases: sliding, escaping and regular dynamics.
\label{fig:phaseSchem}
}
\end{center}
\end{figure}

Without loss of generality, we may assume that the chosen sliding segment ends at the origin,
and that from the origin $\Gamma$ then enters the right half-space,
as shown in Fig.~\ref{fig:phaseSchem}.
Consequently, the right half-flow $\phi^{(R)}$ is tangent to the switching manifold at the origin:
\begin{equation}
e_1^{\sf T} \phi^{(R)}(0) = 0 \;.
\label{eq:a1}
\end{equation}
To ensure a non-degenerate scenario we also require
\begin{equation}
e_1^{\sf T} \phi^{(L)}(0) > 0 \;.
\label{eq:a2}
\end{equation}
For simplicity we choose the coordinate $x_2$ such that at the origin
$\Gamma$ is tangent to the $x_2$-axis and locally the value of $x_2$ increases with time on $\Gamma$.
Therefore
\begin{equation}
e_2^{\sf T} \phi^{(R)}(0) > 0 \;, \qquad
e_j^{\sf T} \phi^{(R)}(0) = 0 \;,~\forall j > 2 \;.
\label{eq:a3}
\end{equation}
In order to ensure $\Gamma$ enters the right half-space from the origin
in a non-degenerate fashion we require
\begin{equation}
e_1^{\sf T} \frac{\partial \phi^{(R)}}{\partial x_2}(0) > 0 \;.
\label{eq:a4}
\end{equation}
Lastly, if the system is at least three-dimensional we may choose the remaining axes 
so that the boundary of the stable sliding region is tangent to $x_2 = 0$ at the origin.
This requirement is indicated in Fig.~\ref{fig:phaseSchem}
and simplifies our expansions about the origin in \S\ref{sec:ESCAPE}.
Algebraically this requirement equates to
\begin{equation}
e_1^{\sf T} \frac{\partial \phi^{(R)}}{\partial x_j}(0) = 0 \;,~\forall j > 2 \;.
\label{eq:a5}
\end{equation}

\subsection{Three dynamical phases}
\label{sub:DIVISON}

The periodic orbit $\Gamma$ may have many sliding segments\removableFootnote{
Periodic orbits that involve transversal intersections with the switching manifold
can be treated as having consecutive excursions.
}.
The assumptions (\ref{eq:a1})-(\ref{eq:a5}) ensure that
our coordinate system is centred at the end point of an arbitrarily chosen sliding segment
in a convenient fashion.
Such a point corresponds to a transition from sliding dynamics to regular dynamics.
In the presence of noise, we prefer to treat this transition as a sequence of three phases:
a sliding phase, an escaping phase, and a regular phase.

We now define these phases precisely and introduce notation used in the remaining sections of the paper.
In order to treat escape from the vicinity of the switching manifold separately, we introduce small constants
$\delta^-$ and $\delta^+$ that satisfy
\begin{equation}
\delta^- < 0 < \delta^+ \;, \qquad
\left| \delta^{\pm} \right| \ll 1 \;.
\label{eq:deltaMinusdeltaPlus}
\end{equation}
For a sample solution to (\ref{eq:sde}) with (\ref{eq:a1})-(\ref{eq:a5}) and $\ee > 0$
that follows a path close to $\Gamma$,
we define the {\em sliding phase} as the part of the solution between
the point at which it arrives at the switching manifold and its first intersection with $x_2 = \delta^-$.
This is followed by an {\em escaping phase}
defined as the part of the solution between the end point of the sliding phase
and its first intersection with $x_2 = \delta^+$.
Lastly we refer to the subsequent part of the solution ending with its
next intersection with the switching manifold as a {\em regular phase}.
Throughout this paper we consider approximations to quantities such as transitional PDFs
and first passage times and locations.
Since these relate to finite intervals of time,
with the assumption that $\ee$ is sufficiently small,
a wild departure of a sample solution from close proximity to $\Gamma$
occurs sufficiently rarely that such large deviations may be ignored.

As indicated in Fig.~\ref{fig:phaseSchem}, for the periodic orbit $\Gamma$ we let $\bx_{\Gamma}^M$
denote the point at which the sliding segment starts\removableFootnote{
$M$ for manifold.
},
let $\bx_{\Gamma}^S$ denote the point at which the sliding segment intersects $x_2 = \delta^-$,
let $\bx_{\Gamma}^E$ denote the point at which $\Gamma$ intersects $x_2 = \delta^+$,
and let $\bx_{\Gamma}^R$ denote the next point at which $\Gamma$ intersects the switching manifold.
We let $t_{\Gamma}^S$, $t_{\Gamma}^E$ and $t_{\Gamma}^R$
denote the deterministic evolution times for the sliding, escaping and regular phases, respectively.



\subsection{The coordinate change for the relay control system}
\label{sub:COORDRCS}

Here we change the coordinates of the relay control example,
(\ref{eq:relayControlSystem2}) with (\ref{eq:ABCDvalues3d})-(\ref{eq:paramValues}),
so that it conforms to the general system (\ref{eq:sde}) with assumptions (\ref{eq:a1})-(\ref{eq:a5}).
The periodic orbit $\Gamma$, shown in Fig.~\ref{fig:ppSketch}, has two sliding segments.
We choose the end point of the upper sliding segment
(the sliding segment with $X_3 > 0$) as the centre of the new coordinates.
Since $\Gamma$ is symmetric
(specifically (\ref{eq:relayControlSystem2}) is unchanged under $\bX \mapsto -\bX$),
the results obtained for the three phases associated with this end point
may be applied directly to the remaining half of $\Gamma$.

To determine the location of the end point of the upper sliding segment,
we note that while $X_1 < 0$ trajectories rapidly contract to a one-dimensional weakly stable manifold.
The intersection of this stable manifold with the switching manifold
provides a suitable approximation to the starting point of the upper sliding segment.
Then from Filippov's solution for sliding motion we find
that upper sliding segment of $\Gamma$ ends at $\bX = (0,1,Z)^{\sf T}$, where $Z \approx 2.561$.
The relevant calculations for this derivation
and an exact expression for $Z$ are given in Appendix \ref{sec:DET}.

Further calculations reveal that the affine change of coordinates
\begin{equation}
\bx = P \bX + Q \;,
\label{eq:Xtox}
\end{equation}
where
\begin{equation}
P = \left[ \begin{array}{ccc}
1 & 0 & 0 \\
0 & 1 & 0 \\
0 & \frac{1}{Z + 2} & 1
\end{array} \right] \;, \qquad
Q = \left[ \begin{array}{c}
0 \\ -1 \\ -\frac{1}{Z + 2} - Z
\end{array} \right] \;,
\end{equation}
transforms (\ref{eq:relayControlSystem2}) with $\ee = 0$ to a system
satisfying (\ref{eq:a1})-(\ref{eq:a5}).
Specifically, under (\ref{eq:Xtox}) the system
(\ref{eq:relayControlSystem2}) with (\ref{eq:ABCDvalues3d})-(\ref{eq:paramValues}) becomes
\begin{equation}
d\bx(t) = \left\{ \begin{array}{lc}
\mathcal{A} \bx(t) + \mathcal{B}^{(L)} \;, & x_1(t) < 0 \\
\mathcal{A} \bx(t) + \mathcal{B}^{(R)} \;, & x_1(t) > 0
\end{array} \right\} \,dt + \sqrt{\ee} D \,d\bW(t) \;,
\label{eq:relayControlSystem5}
\end{equation}
where\removableFootnote{
Also
$\mathcal{B}^{(L)} = P B - P A P^{-1} Q$,
$\mathcal{B}^{(R)} = - P B - P A P^{-1} Q$,
$D = P B e_1^{\sf T}$.
}
\begin{equation}
\mathcal{A} = P A P^{-1} =
\left[ \begin{array}{ccc}
-2 \zeta \omega - \lambda & 1 & 0 \\
-2 \zeta \omega \lambda - \omega^2 & \frac{-1}{Z+2} & 1 \\
-\lambda \omega^2 - \frac{2 \zeta \omega \lambda + \omega^2}{Z+2} &
\frac{-1}{(Z+2)^2} &
\frac{1}{Z+2}
\end{array} \right] \;, \label{eq:calA}
\end{equation}
\begin{equation}
\mathcal{B}^{(L)} =
\left[ \begin{array}{c}
2 \\ Z-2 \\ \frac{2 Z}{Z+2}
\end{array} \right] \;, \qquad
\mathcal{B}^{(R)} =
\left[ \begin{array}{c}
0 \\ Z+2 \\ 0
\end{array} \right] \;, \qquad
D =
\left[ \begin{array}{ccc}
1 & 0 & 0 \\
-2 & 0 & 0 \\
\frac{Z}{Z+2} & 0 & 0
\end{array} \right] \;.
\label{eq:calBLBRD}
\end{equation}
The system (\ref{eq:relayControlSystem5}) with (\ref{eq:calA})-(\ref{eq:calBLBRD})
is used as an example to illustrate our methods in the next three sections.

\section{Regular stochastic dynamics}
\label{sec:EXCUR}
\setcounter{equation}{0}

Here we consider sample solutions to (\ref{eq:sde})-(\ref{eq:a5}) that start
from an initial point, $\bx_0$, on $x_2 = \delta^+$, to an intersection with the switching manifold, $x_1 = 0$.
For an arbitrary sample solution, 
we let $t^R$ denote the first passage time to $x_1 = 0$ and let $\bx^R$ denote
the corresponding arrival location of the solution.
When $\ee = 0$, these values are deterministic and we denote them by
$t_{\dd}^R$ and $\bx_{\dd}^R$ respectively.
Naturally the values of $t_{\dd}^R$ and $\bx_{\dd}^R$ depend on $\bx_0$,
but in this section it is convenient to ignore this dependency because
here we are not interested in variations in $\bx_0$.
Such variations are considered in \S\ref{sec:COMB}
where it is necessary to calculate exactly how deviations
in one phase of the dynamics influence dynamics in subsequent phases.
Note, if $\bx_0 = \bx_{\Gamma}^E$,
then $t_{\dd}^R = t_{\Gamma}^R$ and $\bx_{\dd}^R = \bx_{\Gamma}^R$, Fig.~\ref{fig:phaseSchem}.

We assume $\ee$ is small and $\bx_0$ is sufficiently far from $x_1 = 0$
such that $t^R \approx t_{\dd}^R$ and $\bx^R \approx \bx_{\dd}^R$, with high probability.
Indeed, we have $t^R - t_{\dd}^R = O(\sqrt{\ee})$
and $\bx^R - \bx_{\dd}^R = O(\sqrt{\ee})$.
Dynamics for this phase are governed purely by the right half-system of (\ref{eq:sde}):
\begin{equation}
d\bx(t) = \phi^{(R)}(\bx(t)) \,dt + \sqrt{\ee} D \,d\bW(t) \;.
\label{eq:sdeR}
\end{equation}
It suffices to use classical methods of analysis to study (\ref{eq:sdeR}).
We begin by obtaining $t^R$ and $\bx^R$ to $O(\sqrt{\ee})$
by sample path methods following \cite{FrWe12}.

\subsection{First order approximations}
\label{sub:REGVAR}

For the smooth stochastic differential equation, (\ref{eq:sdeR}),
we can expand $\bx(t)$ as a series involving powers of $\sqrt{\ee}$.
Specifically, by Theorem 2.2 of Chapter 2 of \cite{FrWe12} we can write
\begin{equation}
\bx(t) = \bx_{\dd}(t;\bx_0) + \sqrt{\ee} \,\bx^{(1)}(t) + o(\sqrt{\ee}) \;,
\label{eq:xExpFW}
\end{equation}
where $\bx_{\dd}$ denotes the solution to $\dot{\bx} = \phi^{(R)}(\bx)$, from $\bx_0$,
and $\bx^{(1)}$ satisfies
\begin{equation}
d\bx^{(1)}(t) = \rD_\bx \phi^{(R)}(\bx_{\dd}(t)) \,\bx^{(1)}(t) \,dt + D \,d\bW(t) \;,
\qquad \bx^{(1)}(0) = 0 \;,
\label{eq:OUFW}
\end{equation}
where $\rD_\bx \phi^{(R)}$ is the Jacobian of $\phi^{(R)}$.
Equation (\ref{eq:OUFW}) is a time-dependent Ornstein-Uhlenbeck process \cite{Sc10,Ga09}.
By using an integrating factor we obtain the explicit solution
\begin{equation}
\bx^{(1)}(t) = \int_0^t {\rm e}^{\int_s^t \rD_\bx \phi^{(R)}(\bx_{\dd}(\tilde{s})) \,d\tilde{s}} D \,d\bW(s) \;.
\end{equation}
Consequently $\bx^{(1)}(t)$ is a Gaussian random variable with zero mean and covariance matrix
\begin{equation}
K(t) = \int_0^t H(s,t) H(s,t)^{\sf T} \,ds \;,
\label{eq:XiFW}
\end{equation}
where
\begin{equation}
H(s,t) = {\rm e}^{\int_s^t \rD_\bx \phi^{(R)}(\bx_{\dd}(\tilde{s})) \,d\tilde{s}} D \;.
\label{eq:HFW}
\end{equation}

When $\ee = 0$, first passage to the switching manifold occurs at the point,
$\bx_{\dd}^R = \bx_{\dd}(t_{\dd}^R)$.
We assume that the deterministic solution intersects the switching manifold transversely
at this point, as is generically the case.
That is, we assume $e_1^{\sf T} v \ne 0$ where
\begin{equation}
v = \phi^{(R)}(\bx_{\dd}(t_{\dd}^R)) \;.
\label{eq:v}
\end{equation}
Then, for $\ee > 0$, by Theorem 2.3 of Chapter 2 of \cite{FrWe12},
the first passage statistics satisfy
\begin{eqnarray}
\mathbb{E}[t^R] &=& t_{\dd}^R + o(\sqrt{\ee}) \;,
\label{eq:meantex0} \\
{\rm Var}(t^R) &=& \frac{e_1^{\sf T} K(t_{\dd}^R) e_1}{(e_1^{\sf T} v)^2} \,\ee + o(\ee) \;,
\label{eq:vartex}
\end{eqnarray}
and
\begin{eqnarray}
\mathbb{E}[\bx^R] &=& \bx_{\dd}^R + o(\sqrt{\ee}) \;,
\label{eq:meanxex0} \\
{\rm Cov}(\bx^R) &=& \left( I - \frac{v e_1^{\sf T}}{e_1^{\sf T} v} \right)
K(t_{\dd}^R) \left( I - \frac{v e_1^{\sf T}}{e_1^{\sf T} v} \right)^{\sf T} \ee + o(\ee) \;.
\label{eq:covxex}
\end{eqnarray}

\begin{figure}[b!]
\begin{center}
\setlength{\unitlength}{1cm}
\begin{picture}(15,13)
\put(3.75,7){\includegraphics[height=6cm]{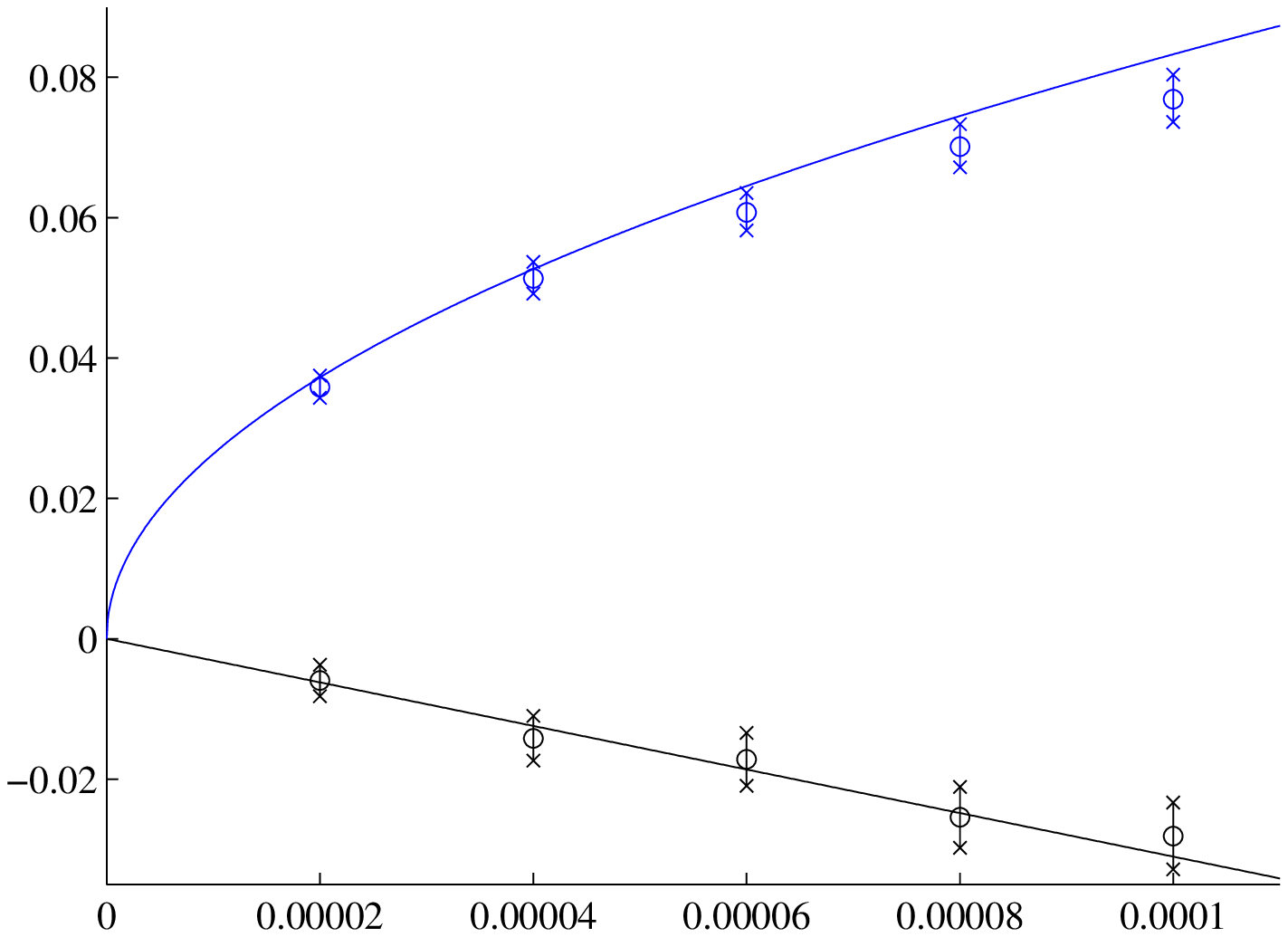}}
\put(0,0){\includegraphics[height=6cm]{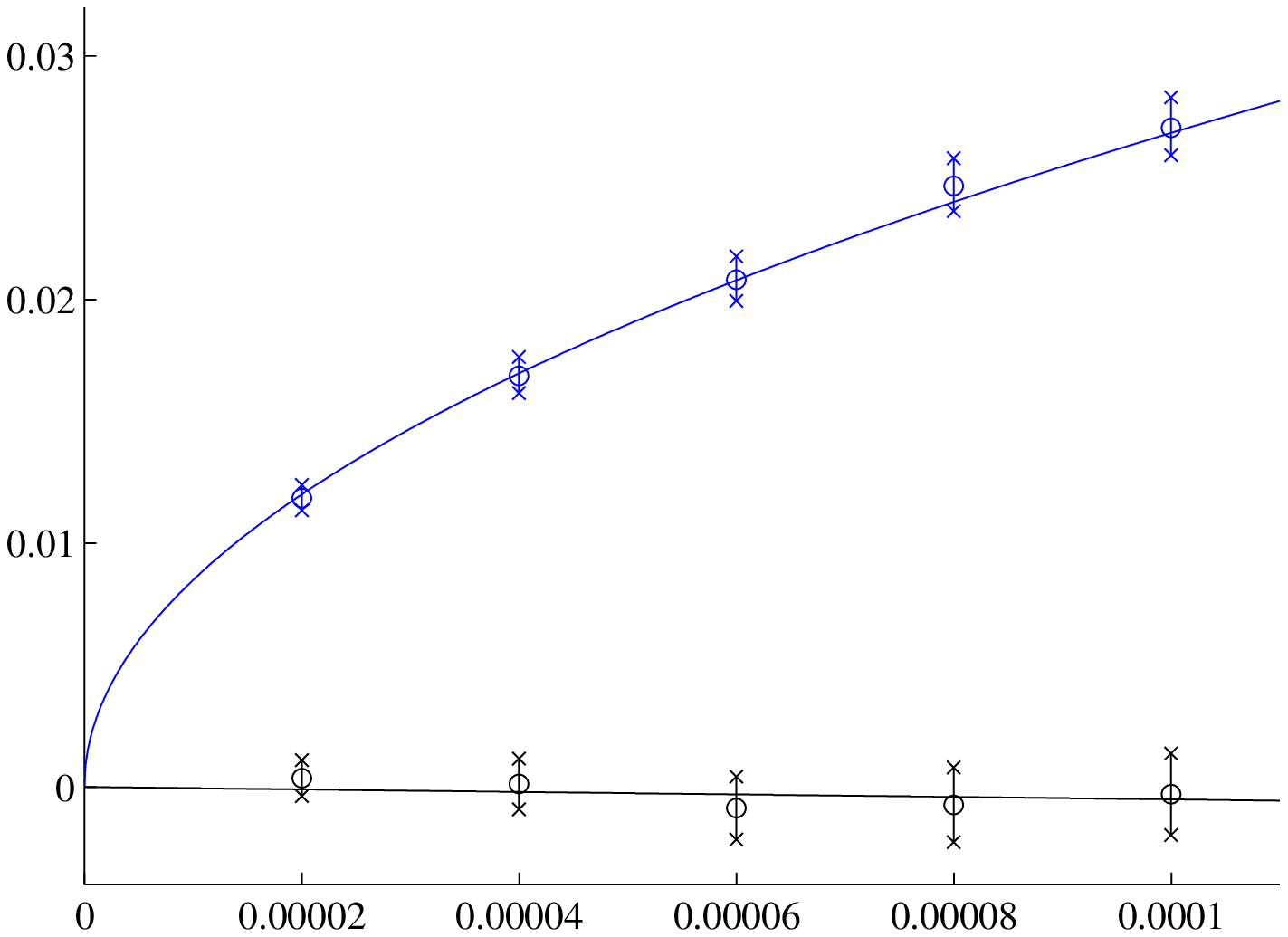}}
\put(7.5,0){\includegraphics[height=6cm]{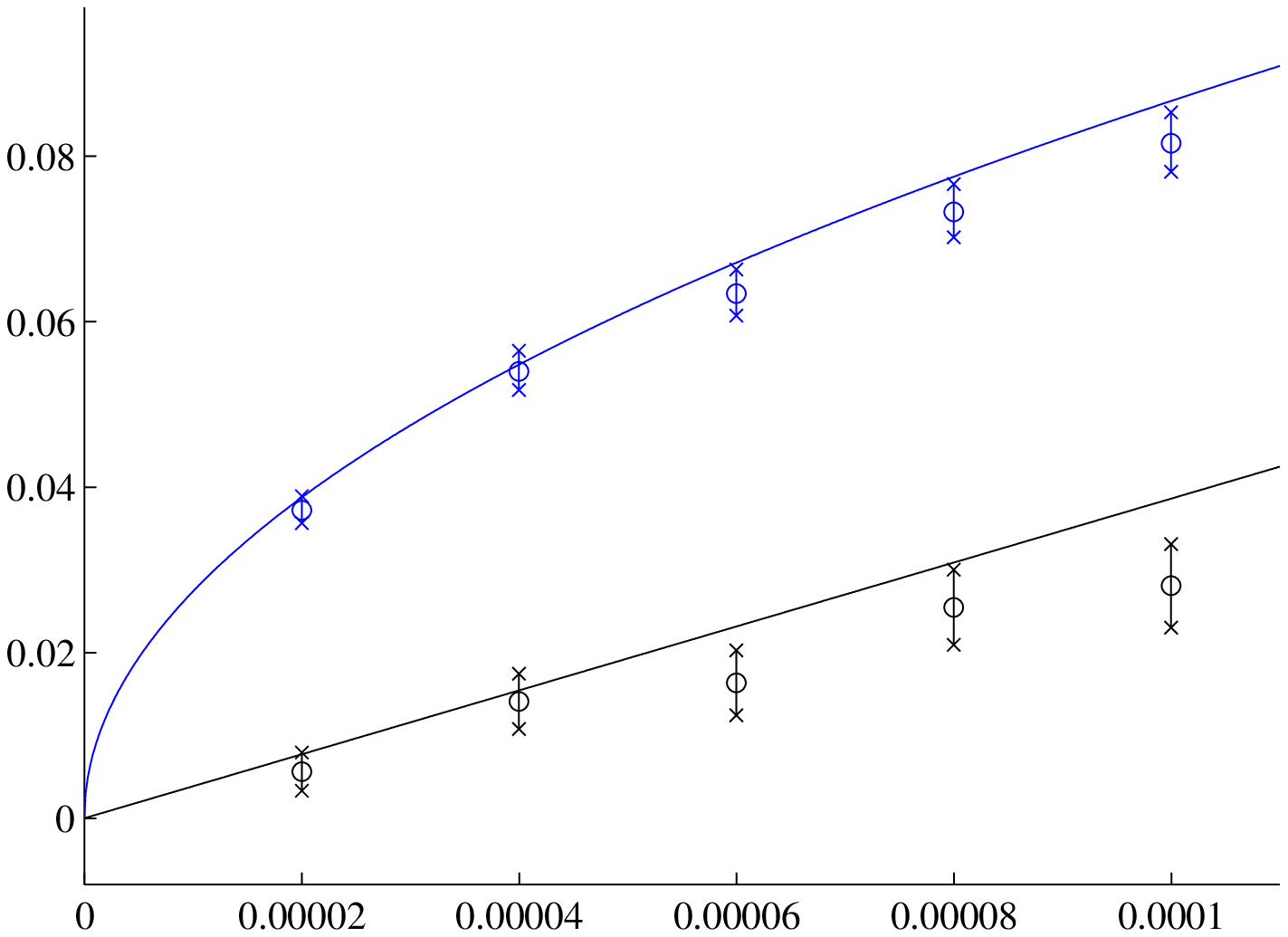}}
\put(4.75,12.8){\large \sf \bfseries A}
\put(1,5.8){\large \sf \bfseries B}
\put(8.5,5.8){\large \sf \bfseries C}
\put(7.55,7){\small $\ee$}
\put(8.8,8.6){\small ${\rm Diff} \left( t^R \right)$}
\put(6.5,12.05){\small \color{blue} ${\rm Std} \left( t^R \right)$}
\put(3.8,0){\small $\ee$}
\put(5,1.7){\small ${\rm Diff} \left( x_2^R \right)$}
\put(2.7,4.6){\small \color{blue} ${\rm Std} \left( x_2^R \right)$}
\put(11.3,0){\small $\ee$}
\put(12.5,1.6){\small ${\rm Diff} \left( x_3^R \right)$}
\put(10.2,4.6){\small \color{blue} ${\rm Std} \left( x_3^R \right)$}
\end{picture}
\caption{
First passage statistics for a regular phase of the dynamics for the relay control
system (\ref{eq:ABCDvalues3d})-(\ref{eq:relayControlSystem2}).
Theoretical results (shown as solid curves) are compared with Monte-Carlo simulations using the transformed system (\ref{eq:relayControlSystem5})-(\ref{eq:calBLBRD}).
In panel A, ${\rm Std}(t^R)$ and ${\rm Diff}(t^R)$ are approximated using (\ref{eq:vartex}) and (\ref{eq:meantex3}), respectively.
In panels B and C, ${\rm Std}(x_2^R)$ and ${\rm Std}(x_3^R)$ are approximated using (\ref{eq:covxex}),
and ${\rm Diff}(x_2^R)$ and ${\rm Diff}(x_3^R)$ are approximated using (\ref{eq:meanxjex2}).
We have used, $\bx_0 = \bx_{\Gamma}^E$, that is, the initial point is
the intersection of the deterministic periodic orbit, $\Gamma$, with $x_2 = \delta^+$,
(using $\delta^+ = 0.2$, as discussed in \S\ref{sec:COMB}).
The deterministic passage time and location
are $t_{\dd}^R \approx 4.263$ and $\bx_{\dd}^R \approx (0,-0.040,-4.770)$.
The data points were computed from $1000$ Monte-Carlo simulations for each value of $\ee$
using the Euler-Maruyama method with a fixed step size, $\Delta t = 0.00001$,
and $95\%$ confidence intervals are indicated as in Fig.~\ref{fig:manyPeriod}.
\label{fig:checkExcur}
}
\end{center}
\end{figure}

In Fig.~\ref{fig:checkExcur}, the formulas (\ref{eq:vartex}) and (\ref{eq:covxex}) are compared with Monte-Carlo simulations
for the relay control example, (\ref{eq:relayControlSystem2}).
The upper curve of panel A is the square root of the leading order term of (\ref{eq:vartex}).
We compute $K$ in this expression by numerically evaluating the integral (\ref{eq:XiFW}).
The upper curves of the lower two panels of Fig.~\ref{fig:checkExcur}
correspond to analogous lowest-order approximations using (\ref{eq:covxex}).

\subsection{Deviations in the mean}
\label{sub:REGMEAN}

From (\ref{eq:meantex0}) and (\ref{eq:meanxex0}) we see that to 
determine the lowest order nonzero terms of
${\rm Diff}(t^R) = \mathbb{E}[t^R] - t_{\Gamma}^R$ and 
${\rm Diff}(\bx^R) = \mathbb{E}[\bx^R] - \bx_{\Gamma}^R$,
more powerful methods are required.
To calculate $\mathbb{E}[t^R]$ to $O(\ee)$ we express this quantity
in terms of the transitional PDF (probability density function) $p^R(\bx,t;\bx_0)$ for (\ref{eq:sdeR}),
\begin{equation}
\mathbb{E}[t^R] = \int_0^\infty \int_{\mathbb{R}^N} p^R(\bx,t) \,d\bx \,dt \;.
\label{eq:meantex1}
\end{equation}
The Fokker-Planck equation for $p^R$ with corresponding boundary conditions, including an absorbing barrier 
at $x_1=0$, is considered below in (\ref{eq:FPEex})-(\ref{eq:FPEremainingBCs}).
A determination of $\mathbb{E}[\bx^R]$ requires using the joint PDF for the first passage time and location, $k(x_2,\ldots,x_N,t;\bx_0)$,
\begin{equation}
\mathbb{E}[x_j^R] = \int_0^\infty \int_{\mathbb{R}^{N-1}} x_j
k(x_2,\ldots,x_N,t) \,dx_2 \ldots \,dx_N \,dt \;.
\label{eq:meanxjex0}
\end{equation}
The joint PDF $k$ can be written in terms of $p^R$ through the probability current $J$ \cite{Sc10,Ga09} of (\ref{eq:sdeR}),
\begin{equation}
k(x_2,\ldots,x_N,t) = -e_1^{\sf T} J(0,x_2,\ldots,x_N,t) \;,
\end{equation}
where $J$ is given by
\begin{equation}
J(\bx,t) = \phi^{(R)} p^R - \frac{\ee}{2} \left[ \begin{array}{c}
\sum_{j=1}^N (D D^{\sf T})_{1j} \frac{\partial p^R}{\partial x_j} \\
\vdots \\
\sum_{j=1}^N (D D^{\sf T})_{Nj} \frac{\partial p^R}{\partial x_j}
\end{array} \right] \;.
\end{equation}
In view of the absorbing barrier at $x_1=0$, we have
\begin{equation}
k(x_2,\ldots,x_N,t) = \frac{\ee}{2} (D D^{\sf T})_{11} \frac{\partial p^R}{\partial x_1}(0,x_2,\ldots,x_N,t) \;,
\label{eq:k2}
\end{equation}
so that (\ref{eq:meanxjex0}) becomes
\begin{equation}
\mathbb{E}[x_j^R] = \frac{\ee}{2} (D D^{\sf T})_{11}
\int_0^\infty \int_{\mathbb{R}^{N-1}} x_j
\frac{\partial p^R}{\partial x_1}(0,x_2,\ldots,x_N,t)
\,dx_2 \ldots \,dx_N \,dt \;.
\label{eq:meanxjex1}
\end{equation}

The expressions for $\mathbb{E}[t^R]$ and $\mathbb{E}[x_j^R]$ involve $p^R$, which satisfies the following boundary value problem
\begin{eqnarray}
\frac{\partial p^R}{\partial t} &=& 
-\sum_{i=1}^N \frac{\partial}{\partial x_i} \left( \phi_i^{(R)}(\bx) p^R \right) +
\frac{\ee}{2} \sum_{i=1}^N \sum_{j=1}^N
(D D^{\sf T})_{ij} \frac{\partial^2 p^R}{\partial x_i \partial x_j} \;, \label{eq:FPEex} \\
p^R(\bx,0;\bx_0) &=& \delta(\bx - \bx_0) \;, \label{eq:FPEIC} \\
p^R(\bx,t;\bx_0) &=& 0 \;, {\rm ~whenever~} x_1=0 \;, \label{eq:FPEabsorbingBC} \\
p^R(\bx,t;\bx_0) &\to& 0 {\rm ~as~} ||\bx|| \to \infty \;, {\rm ~with~} x_1 \ge 0 \label{eq:FPEremainingBCs} \;.
\end{eqnarray}
Here (\ref{eq:FPEIC}) captures the initial condition,
(\ref{eq:FPEabsorbingBC}) is the absorbing boundary condition representing that the switching manifold is an absorbing barrier,
with (\ref{eq:FPEremainingBCs}) ensuring physically realistic behaviour.

The form of (\ref{eq:FPEex})-(\ref{eq:FPEremainingBCs}) suggests that the problem
can be solved using an asymptotic approach \cite{Za06,Ve05,KaPa03}.
The solution to (\ref{eq:FPEex})-(\ref{eq:FPEremainingBCs})
can be constructed using the solution to (\ref{eq:FPEex}) without the boundary conditions,
that is, the free-space solution $p_f(\bx,t;\bx_0)$.
With $\ee \ll 1$, $p_f$ is concentrated around the deterministic trajectory $\bx_{\dd}(t)$.
The PDF $p_f$ does not satisfy the boundary condition at $x_1=0$,
in particular near the point $\bx_{\dd}^R$.

A boundary layer analysis in terms of the local variable $z = \frac{x_1}{\ee}$
indicates that near the switching manifold
the boundary layer behaviour of the PDF is given by
the sum $p_{\ell}(z,x_2,\ldots,x_N,t) + p_f|_{x_1=0}$,
where $p_{\ell}$ is a local contribution to be determined.
The uniform solution is then obtained by matching the boundary layer solution to the outer free space solution,
$\lim_{z \to \infty}(p_{\ell} + p_f|_{x_1=0}) = \lim_{x_1 \to 0} p_f(\bx,t;\bx_0)$.
This implies $\lim_{z \to \infty} p_{\ell} = 0$, and hence
the uniform solution has the form $p^R = p_{\ell} + p_f$.
Furthermore, since $p_f$ decays exponentially away from its main concentration
around the deterministic trajectory, $\bx_{\dd}$,
it remains to find the local contribution $p_{\ell}$ that decays away from $\bx_{\dd}^R$.  
To this end it is appropriate to use the local variables,
\begin{equation}
z = \frac{x_1}{\ee} \;, \qquad
u_j = \frac{x_j - x_{{\dd},j}^R}{\sqrt{\ee}} \;, \forall j \ne 1 \;, \qquad
\tau = \frac{t - t_{\dd}^R}{\sqrt{\ee}} \;,
\label{eq:regularScaling}
\end{equation}
where the particular scaling for $\tau$ and $u_j$
is motivated by (\ref{eq:vartex}) and (\ref{eq:covxex}) respectively.
Then the local contribution can be written
\begin{eqnarray}
&& p_{\ell} \left( \ee z,\sqrt{\ee} u_2 + x_{{\dd},2}^R,\ldots,\sqrt{\ee} u_N +
x_{{\dd},N}^R,\sqrt{\ee} \tau + t_{\dd}^R \right)
= \ee^{-\frac{N}{2}} \mathcal{P}(z,u_2,\ldots,u_N,\tau) \nonumber \\
&&= \ee^{-\frac{N}{2}} \left( \mathcal{P}^{(0)}(z,u_2,\ldots,u_N,\tau)
+ \sqrt{\ee} \mathcal{P}^{(1)}(z,u_2,\ldots,u_N,\tau) + O(\ee) \right) \;,
\label{eq:pell}
\end{eqnarray}
with the expansion in powers of $\sqrt{\ee}$. 
The PDE and behaviour at infinity for $\mathcal{P}$ are given by
\begin{align}
\frac{1}{\sqrt{\ee}} \frac{\partial \mathcal{P}}{\partial \tau} &=
-\frac{1}{\ee} \phi_1^{(R)}(\bx_d^R)
\frac{\partial \mathcal{P}}{\partial z}
-\frac{1}{\sqrt{\ee}} \sum_{i=2}^N
\frac{\partial \phi_1^{(R)}}{\partial x_i}(\bx_d^R) u_i \frac{\partial \mathcal{P}}{\partial z}
-\frac{1}{\sqrt{\ee}} \sum_{i=2}^N
\phi_i^{(R)}(\bx_d^R) \frac{\partial \mathcal{P}}{\partial u_i} \nonumber \\
&+~\frac{1}{2 \ee} \left( D D^{\sf T} \right)_{1,1} \frac{\partial^2 \mathcal{P}}{\partial z^2}
+\frac{1}{\sqrt{\ee}} \sum_{i=2}^N
\left( D D^{\sf T} \right)_{i,1} \frac{\partial^2 \mathcal{P}}{\partial z \partial u_i} + O(\ee^0) \;, \label{eq:fpe3} \\
\mathcal{P} &\to 0 {\rm ~as~} z \to \infty {\rm ~or~} u_j \to \pm \infty \;, \forall j \ne 1 \;. \label{eq:remainingBCs}
\end{align}
The absorbing boundary condition for $p^R$ (\ref{eq:FPEabsorbingBC})
gives the condition for $\mathcal{P}$ at $z=0$,
\begin{eqnarray}
&& -p_f(0,\sqrt{\ee} u_2 + x_{{\dd},2}^R,\ldots,\sqrt{\ee} u_N + x_{{\dd},N}^R,\sqrt{\ee} \tau + t_{\dd}^R)
= \ee^{-\frac{N}{2}} \mathcal{P}(0,u_2,\ldots,u_N,\tau) \nonumber \\
&&= \ee^{-\frac{N}{2}} \left( f^{(0)}(u_2,\ldots,u_N,\tau) + \sqrt{\ee} f^{(1)}(u_2,\ldots,u_N,\tau) + O(\ee) \right)
\;. \label{eq:exBC0}
\end{eqnarray}
for some functions $f^{(i)}$.

\subsection{Regular dynamics for the relay control model}
\label{sub:REGRELAY}

Here we summarize the key steps in calculating the functions $f^{(i)}$ that appear in (\ref{eq:exBC0}),
and ultimately calculating $\mathbb{E}[t^R]$ and $\mathbb{E}[\bx^R]$ for the relay control example.
Details are deferred to Appendix \ref{sec:CALEX}.

As described in \S\ref{sub:COORDRCS},
in transformed coordinates the right half-system of the relay control example is
\begin{equation}
d\bx(t) = \left( \mathcal{A} \bx(t) + \mathcal{B}^{(R)} \right) dt
+ \sqrt{\ee} D \,d\bW(t) \;,
\label{eq:relayControlSystem5R}
\end{equation}
with (\ref{eq:calA}) and (\ref{eq:calBLBRD}).
The transitional PDF of (\ref{eq:relayControlSystem5R}) takes the form
$p^R(\bx,t) = p_f(\bx,t) + \ee^{-\frac{N}{2}} \mathcal{P}(z,u_2,\ldots,u_N,\tau)$
with initial and boundary conditions (\ref{eq:FPEIC})-(\ref{eq:FPEremainingBCs}).
The free-space contribution is given by
\begin{equation}
p_f(\bx,t) =
\frac{1}{(2 \pi \ee)^{\frac{3}{2}} \sqrt{\det(K(t))}}
{\rm e}^{-\frac{1}{2 \ee} (\bx - \bx_{\dd}(t))^{\sf T}
K(t)^{-1} (\bx - \bx_{\dd}(t))} \;,
\label{eq:pfs}
\end{equation}
with covariance matrix
\begin{equation}
K(t) = \int_0^t {\rm e}^{\mathcal{A} s}
D D^{\sf T}
{\rm e}^{\mathcal{A}^{\sf T} s} \,ds \;.
\end{equation}
The local contribution $\mathcal{P}$ satisfies (\ref{eq:fpe3})
with $\phi_1^{(R)}(\bx_{\Gamma}^R) = x_{\Gamma,2}^R$ 
and $\left( D D^{\sf T} \right)_{1,1} = 1$
(coefficients for higher order terms are given in Appendix \ref{sec:CALEX}).
Note that $t_{\Gamma}^R$ and $\bx_\Gamma^R$ are
obtained by solving (\ref{eq:relayControlSystem5R}) with $\ee = 0$
(refer to (\ref{eq:bxdetrelay})).
 
With $\mathcal{P}$ expanded as in (\ref{eq:pell}), the $O(1)$ equation is
\begin{equation}
\frac{1}{2} \mathcal{P}^{(0)}_{z z} - x_{\Gamma,2}^R \mathcal{P}^{(0)}_{z} = 0 \;.
\label{eq:calP0}
\end{equation}
Therefore
\begin{equation}
\mathcal{P}^{(0)} = -f^{(0)}(u_2,u_3,\tau) \,{\rm e}^{2 x_{\Gamma,2}^R z} \;,
\label{eq:calP02}
\end{equation}
where $f^{(0)}$ is determined from (\ref{eq:exBC0}) and (\ref{eq:pfs}) and given by (\ref{eq:f0}).
Note $x_{\Gamma,2}^R < 0$,
because $\phi_1^{(R)}(\bx_{\Gamma}^R) = x_{\Gamma,2}^R$ is the component of the vector field
orthogonal to the switching manifold evaluated at the deterministic passage location.
Higher order corrections to $\mathcal{P}$ satisfy
\begin{equation}
\frac{1}{2} \mathcal{P}^{(j)}_{z z} - x_{\Gamma,2}^R \mathcal{P}^{(j)}_{z} =
\mathcal{F} \left( \mathcal{P}^{(0)},\ldots,\mathcal{P}^{(j-1)} \right) \;,
\label{eq:calPj}
\end{equation}
where the right hand-side is a function of the lower order components and their derivatives.
By solving (\ref{eq:calPj}) with $j = 1$ using (\ref{eq:calP02}) and the boundary conditions
(\ref{eq:remainingBCs})-(\ref{eq:exBC0}), we obtain
\begin{equation}
\mathcal{P}^{(1)} = \left( -f^{(1)}(u_2,u_3,\tau) + g^{(1)}(u_2,u_3,\tau) z \right) {\rm e}^{2 x_{\Gamma,2}^R z} \;,
\label{eq:calP12}
\end{equation}
where $f^{(1)}$ is found from (\ref{eq:exBC0}) and (\ref{eq:f0})
and $g^{(1)}$ depends on $f^{(0)}$ and its derivatives (\ref{eq:g1}).

By evaluating (\ref{eq:meantex1}) with (\ref{eq:pfs}) and (\ref{eq:calP02}),
we obtain the following formula for the mean first passage time,
\begin{equation}
\mathbb{E}[t^R] = t_{\Gamma}^R +
\frac{1}{2 \left( \phi_1^{(R)}(\bx_{\Gamma}^R) \right)^2} \left( \dot{\kappa}_{11}(t_{\Gamma}^R)
+ \frac{\kappa_{11}(t_{\Gamma}^R) \ddot{x}_{{\dd},1}(t_{\Gamma}^R)}{\phi_1^{(R)}(\bx_{\Gamma}^R)} - 1 \right) \ee
+ O \left( \ee^{\frac{3}{2}} \right) \;,
\label{eq:meantex2}
\end{equation}
where $\kappa_{11}$ denotes the $(1,1)$-element of $K$.
As discussed in \S\ref{sub:COORDRCS},
the deterministic passage location $\bx_{\Gamma}^R$ is well approximated
by the intersection of the weakly stable manifold of the right half-space with the switching manifold.
Repeating (\ref{eq:xRint}), in transformed coordinates this intersection point is
\begin{equation}
\bx_{\rm int}^{(R)} =
\left[ 0 \;, -\frac{1}{\omega^2} \;, -2-\frac{2 \zeta}{\omega} - Z - \frac{1}{\omega^2(Z+2)} \right]^{\sf T} \;.
\label{eq:xRint2}
\end{equation}
Indeed, for the parameter values we have used numerical calculations reveal that
$|| \bx_{\Gamma}^R - \bx_{\rm int}^{(R)} || \approx 0.000087$.
By combining (\ref{eq:meantex2}) and (\ref{eq:xRint2}) we obtain
the useful approximation\removableFootnote{
A necessary calculation is:
\begin{eqnarray}
&\dot{\bx}_{\dd}(t_{\dd}^R)
= \mathcal{A} \bx_{\dd}(t_{\dd}^R) + \mathcal{B}^{(R)}
\approx \mathcal{A} \bx_{\rm int}^{(R)} + \mathcal{B}^{(R)}
= \left[ \begin{array}{c}
-\frac{1}{\omega^2} \\
-\frac{2 \zeta}{\omega} \\
-1 - \frac{2 \zeta}{\omega (Z+2)}
\end{array} \right]& \;, \\
&\ddot{\bx}_{\dd}(t_{\dd}^R)
= \mathcal{A} \dot{\bx}_{\dd}(t_{\dd}^R)
\approx \left[ \begin{array}{c}
\frac{\lambda}{\omega^2} \\
\frac{2 \zeta}{\omega} \\
\lambda + \frac{2 \zeta}{\omega (Z+2)}
\left( 1 + \lambda - \frac{1}{Z+2} \right)
\end{array} \right]& \;.
\end{eqnarray}
}
\begin{equation}
\mathbb{E}[t^R] \approx t_{\dd}^R +
\frac{\omega^4}{2} \left( \dot{\kappa}_{11}(t_{\dd}^R)
+ \lambda \kappa_{11}(t_{\dd}^R) - 1 \right) \ee
+ O \left( \ee^{\frac{3}{2}} \right) \;.
\label{eq:meantex3}
\end{equation}
As shown in Fig.~\ref{fig:checkExcur}-A, (\ref{eq:meantex3}) is consistent with Monte-Carlo simulations.

The mean values $\mathbb{E}[\bx_2^R]$ and $\mathbb{E}[\bx_3^R]$ cannot be expressed as concisely.
Numerically it is convenient to compute these values by evaluating $p_f$ directly, rather than calculating $f^{(1)}$,
as shown in Appendix \ref{sec:CALEX}.
The results, shown in panels B and C of Fig.~\ref{fig:checkExcur},
were computed by numerically evaluating (\ref{eq:meanxjex2}).

\section{Stochastically perturbed sliding dynamics}
\label{sec:SLIDE}
\setcounter{equation}{0}

Here we consider sample solutions to (\ref{eq:sde}) with assumptions (\ref{eq:a1})-(\ref{eq:a5})
from an initial point, $\bx_0$, on $x_1 = 0$, until an intersection $x_2 = \delta^-$.
For an arbitrary sample solution, 
we let $t^S$ and $\bx^S$ denote the first passage time and location to $x_2 = \delta^-$.
When $\ee = 0$, these values are deterministic and we denote them by
$t_{\dd}^S$ and $\bx_{\dd}^S$ respectively.

In this section it is convenient to write
\begin{equation}
\by = [ x_2, \ldots, x_N ]^{\sf T} \;, \qquad
\psi^{(L)} = [ \phi^{(L)}_2, \ldots, \phi^{(L)}_N ]^{\sf T} \;, \qquad
\psi^{(R)} = [ \phi^{(R)}_2, \ldots, \phi^{(R)}_N ]^{\sf T} \;,
\label{eq:slidingVariables}
\end{equation}
with which the piecewise-smooth stochastic differential equation (\ref{eq:sde}) may be written as
\begin{equation}
\left[ \begin{array}{c} dx_1(t) \\ d\by(t) \end{array} \right] =
\left\{ \begin{array}{lc}
\left[ \begin{array}{c} \phi^{(L)}_1(x_1(t),\by(t)) \\ \psi^{(L)}(x_1(t),\by(t)) \end{array} \right]
\;, & x_1(t) < 0 \\
\left[ \begin{array}{c} \phi^{(R)}_1(x_1(t),\by(t)) \\ \psi^{(R)}(x_1(t),\by(t)) \end{array} \right]
\;, & x_1(t) > 0
\end{array} \right\} \,dt + \sqrt{\ee} D \,d\bW(t) \;.
\label{eq:sde2}
\end{equation}
Since sample solutions remain near the switching manifold with high probability,
it is profitable to expand in $x_1$.
We rewrite (\ref{eq:sde2}) as
\begin{equation}
\left[ \begin{array}{c} dx_1(t) \\ d\by(t) \end{array} \right] =
\left\{ \begin{array}{lc}
\left[ \begin{array}{c}
a_L(\by(t)) + c_L(\by(t)) x_1(t) + O(x_1(t)^2) \\
b_L(\by(t)) + d_L(\by(t)) x_1(t) + O(x_1(t)^2)
\end{array} \right]
\;, & x_1(t) < 0 \\
\left[ \begin{array}{c}
-a_R(\by(t)) + c_R(\by(t)) x_1(t) + O(x_1(t)^2) \\
b_R(\by(t)) + d_R(\by(t)) x_1(t) + O(x_1(t)^2)
\end{array} \right]
\;, & x_1(t) > 0
\end{array} \right\} \,dt + \sqrt{\ee} D \,d\bW(t) \;.
\label{eq:sde3}
\end{equation}
Take care to note that $a_L$, $a_R$, $c_L$ and $c_R$ are scalars,
and $b_L$, $b_R$, $d_L$ and $d_R$ are $(N-1)$-dimensional vectors.

When $\ee = 0$ we use Filippov's convention to define a deterministic sliding solution \cite{Fi88,Fi60}.
On $x_1 = 0$, for values of $\by$ for which $a_L > 0$ and $a_R > 0$, we define
\begin{equation}
\left[ \begin{array}{c} \dot{x}_{{\dd},1} \\ \dot{\by}_{\dd} \end{array} \right] =
(1-\mu(\by_{\dd})) \left[ \begin{array}{c} a_L(\by_{\dd}) \\ b_L(\by_{\dd}) \end{array} \right] +
\mu(\by_{\dd}) \left[ \begin{array}{c} -a_R(\by_{\dd}) \\ b_R(\by_{\dd}) \end{array} \right] \;,
\end{equation}
where $\mu$ is given by the requirement $\dot{x}_{{\dd},1} = 0$.
That is
$\mu = \frac{a_L}{a_L+a_R}$,
and hence
\begin{equation}
\dot{\by}_{\dd} = \Omega \equiv \frac{a_L b_R + a_R b_L}{a_L + a_R} \;.
\label{eq:Omega}
\end{equation}
The sliding solution, $\by_{\dd}(t;\by_0)$,
satisfies $\by_{\dd}(0;\by_0) = \by_0$,
and $\dot{\by}_{\dd} = \Omega(\by_{\dd})$.

\subsection{Stochastic averaging}
\label{sub:AVERAGING}

The technique of stochastic averaging applies to stochastic systems with distinct time scales \cite{FrWe12,PaSt08,PaKo74,MoCu11}.
The underlying principle is to average fast variables in order
to obtain a simpler description of the behaviour of slow variables over a relatively long time frame.
Stochastic averaging has been an invaluable tool for understanding
periodically forced oscillators \cite{AnAs02,RoSp86,Da98},
and excitable systems \cite{BeGe06,Ba04b}.
Here we apply stochastic averaging
to determine the first passage statistics of stochastically perturbed sliding motion to
the plane $x_2 = \delta^-$.

From previous investigations \cite{SiKu13c} we know that $x_1(t) = O(\ee)$,
for stochastically perturbed sliding motion of (\ref{eq:sde3}).
This motivates the scaling
\begin{equation}
z = \frac{x_1}{\ee} \;,
\end{equation}
with which (\ref{eq:sde3}) may be written as
\begin{eqnarray}
dz(t) &=& \frac{1}{\ee} \left\{ \begin{array}{lc}
a_L(\by(t)) + \ee c_L(\by(t)) z(t) + O(\ee^2) \;, & z(t) < 0 \\
-a_R(\by(t)) + \ee c_R(\by(t)) z(t) + O(\ee^2) \;, & z(t) > 0
\end{array} \right\} \,dt +
\frac{1}{\sqrt{\ee}} \,e_1^{\sf T} D \,d\bW(t) \;, \label{eq:dhatx3} \\
\,d\by(t) &=& F(z(t),\by(t)) \,dt +
\sqrt{\ee} \left[ \begin{array}{c}
e_2^{\sf T} \\ \vdots \\ e_N^{\sf T}
\end{array} \right]
D \,d\bW(t) \;, \label{eq:dhaty2}
\end{eqnarray}
where
\begin{equation}
F(z,\by) = 
\frac{b_L+b_R}{2}
- \frac{b_L-b_R}{2} \,{\rm sgn}(z)
+ \ee \frac{d_L+d_R}{2} z
- \ee \frac{d_L-d_R}{2} z \,{\rm sgn}(z)
+ O(\ee^2) \;.
\label{eq:F}
\end{equation}
In view of the manner by which $\ee$ appears in (\ref{eq:dhatx3})-(\ref{eq:dhaty2}),
we may treat $z(t)$ and $\by(t)$ as fast and slow variables respectively.
Furthermore, the averaging approximation of (\ref{eq:dhaty2})
involves only terms in the drift coefficients that are of lower order than $\ee$
(i.e.~terms involving $a_L$, $a_R$, $b_L$ and $b_R$).
Below we show that with this approximation the mean of $\by$ coincides with the deterministic solution, $\by_{\dd}$.
A full computation of the noise-induced correction to the mean
due to terms of the next order is beyond the scope of this paper\removableFootnote{
I started to compute fourth order terms in our full asymptotic expansion of the PDF
(given to third order in \cite{SiKu14b}).
Despite using {\sc matlab} to do algebraic manipulations the calculations quickly became very complicated.
(I appear to get a PDE for $f^{(1)}$ that is a nonhomogeneous version
of the PDF for $f^{(0)}$ involving terms $y f^{(0)}$ and $y^3 f^{(0)}$.)
It seems possible to be able to follow the calculations all the way through and obtain the correction semi-explicitly,
but there are so many terms to keep track of,
and some of the expressions are much nastier than those at third order,
that the effort required is far greater than the value of the result itself.
}.
However, for the relay control example some of the coefficients $c_L$, $c_R$, $d_L$ and $d_R$, take relatively large values.
Here we suppose these coefficients are $O \left( \ee^{-\eta} \right)$, for some $\eta > 0$,
so that we can formally derive the correction via a straight-forward averaging approximation.

For any fixed $\by$ satisfying $a_L, a_R > 0$, as detailed in \cite{SiKu13c},
(\ref{eq:dhatx3}) has the quasi-steady-state density
\begin{equation}
p_{\rm qss}(z;\by) = 
\left( \frac{2 a_L a_R}{\alpha (a_L+a_R)}
- \frac{a_L^3 c_R + a_R^3 c_L}{a_L a_R (a_L+a_R)^2} \,\ee + O(\ee^2) \right)
\left\{ \begin{array}{lc}
{\rm e}^{\frac{1}{\alpha} \left( 2 a_L z
+ c_L z^2 \ee + O(\ee^2) \right)} \;, & z < 0 \\
{\rm e}^{-\frac{1}{\alpha} \left( 2 a_R z
- c_R z^2 \ee + O(\ee^2) \right)} \;, & z > 0
\end{array} \right. \;,
\label{eq:pqss}
\end{equation}
where $\alpha = (D D^{\sf T})_{11}$ as in (\ref{eq:alphaBetaGamma}).
Given $\by$, it is suitable to assume $z$ is distributed according to (\ref{eq:pqss}).
Averaging $F$ over $p_{\rm qss}$ yields\removableFootnote{
From (\ref{eq:pqss}) we readily obtain
\begin{eqnarray}
\mathbb{E} \left[ {\rm sgn}(z) \,\big|\, \by \right] &=&
\frac{a_L-a_R}{a_L+a_R}
+ \frac{a_L^2 c_R - a_R^2 c_L}{a_L a_R (a_L+a_R)^2} \alpha \ee + O(\ee^2) \\
\mathbb{E} \left[ z \,\big|\, \by \right] &=&
\frac{a_L-a_R}{2 a_L a_R} \alpha + O(\ee) \\
\mathbb{E} \left[ z \,{\rm sgn}(z) \,\big|\, \by \right] &=&
\frac{a_L^2+a_R^2}{2 a_L a_R (a_L+a_R)} \alpha + O(\ee) \;.
\end{eqnarray}
}
\begin{eqnarray}
\overline{F}(\by) &\equiv& \mathbb{E} \left[ F(z,\by) \,\big|\, \by \right] \nonumber \\
&=& \int_{-\infty}^\infty
F(z,\by) p_{\rm qss}(z;\by) \,dz \nonumber \\
&=& \Omega(\by) + \Lambda(\by) \alpha \ee + O(\ee^2) \;,
\label{eq:meanFcal}
\end{eqnarray}
where $\Omega$ is given by (\ref{eq:Omega}) and
\begin{equation}
\Lambda = \frac{(a_L^2 d_R - a_R^2 d_L)(a_L+a_R)
- (a_L^2 c_R - a_R^2 c_L)(b_L-b_R)}
{2 a_L a_R (a_L+a_R)^2} \;.
\label{eq:Lambda}
\end{equation}
The quantity $\Lambda$ was obtained in \cite{SiKu13c}
in the case that $\phi^{(L)}$ and $\phi^{(R)}$ are independent of $\by$.
The averaged equation of (\ref{eq:dhaty2}) is
$d\overline{\by}(t) = \overline{F}(\overline{\by}(t)) dt$.
Therefore we have the ODE
\begin{equation}
\frac{d\overline{\by}}{dt} = \Omega(\overline{\by})
+ \Lambda(\overline{\by}) \alpha \ee + O(\ee^2) \;.
\label{eq:dyBar}
\end{equation}
Note, in the limit $\ee \to 0$, (\ref{eq:dyBar}) is equal to (\ref{eq:Omega}).
Formally there exists a sequence of stochastic solutions to (\ref{eq:sde2})
that converges to $\by_{\dd}(t)$ as $\ee \to 0$ \cite{BuOu09}.

Since (\ref{eq:dyBar}) is the averaged equation
we can use it to obtain the mean of the first passage statistics for small $\ee$.
The mean of the first passage time, $t^S$, is approximated by
\begin{equation}
e_1^{\sf T} \overline{\by} \left( \mathbb{E} \left[ t^S \right] \right) \approx \delta^- \;,
\label{eq:tBarDef}
\end{equation}
and the mean of the $\by$-component of the first passage location is approximated by
\begin{equation}
\mathbb{E} \left[ \by^S \right] \approx \overline{\by} \left( \mathbb{E} \left[ t^S \right] \right) \;.
\label{eq:byBarDef}
\end{equation}
We expect these approximations to be exact to at least $O(\ee)$, and write
\begin{eqnarray}
\overline{\by}(t) &=& \by_{\dd}(t) + \overline{\by}^{(1)}(t) \ee + O(\ee^2) \;,
\label{eq:yBarSeries} \\
\mathbb{E}[t^S] &=& t_{\dd}^S + t^{S,1} \ee + O(\ee^2) \;,
\label{eq:tSlMeanSeries} \\
\mathbb{E}[\by^S] &=& \by_{\dd}^S + \by^{S,1} \ee + O(\ee^2) \;.
\label{eq:bySlMeanSeries}
\end{eqnarray}
By substituting (\ref{eq:yBarSeries}) into (\ref{eq:dyBar}) we obtain
\begin{equation}
\overline{\by}^{(1)}(t) = \alpha
\int_0^t {\rm e}^{\int_s^t (\rD_\by \Omega)(\by_{\dd}(u)) \,du}
\Lambda(\by_{\dd}(s)) \,ds \;,
\label{eq:yBar1}
\end{equation}
which indicates the leading-order deviation of $\overline{\by}(t)$ from
the deterministic value $\by_{\dd}(t)$.
We then express the leading-order deviations of $\mathbb{E}[t^S]$ and $\mathbb{E}[\by^S]$ in terms of $\overline{\by}^{(1)}(t_{\dd}^S)$.
First, from (\ref{eq:tBarDef}) and (\ref{eq:tSlMeanSeries}),
\begin{equation}
t^{S,1} =
-\frac{e_1^{\sf T} \overline{\by}^{(1)}(t_{\dd}^S)}
{e_1^{\sf T} \Omega(\by_{\dd}(t_{\dd}^S))} \;.
\label{eq:meantsl}
\end{equation}
Second, from (\ref{eq:byBarDef}) and (\ref{eq:bySlMeanSeries}),
\begin{equation}
\by^{S,1} = \Omega(t_{\dd}^S) t^{S,1}
+ \overline{\by}^{(1)}(t_{\dd}^S) \;.
\label{eq:meanbysl1}
\end{equation}
The $x_1$-value of the mean first passage location is found using (\ref{eq:pqss}):
\begin{equation}
\mathbb{E} \left[ x_1^S \right] = \ee \int_{-\infty}^{\infty} z p_{\rm qss}(z;\by) \,dz =
\frac{a_L-a_R}{2 a_L a_R} \bigg|_{\by = \by_{\dd}^S} \alpha \ee + O(\ee^2) \;.
\label{eq:meanx1sl}
\end{equation}

\subsection{Linear diffusion approximation}

Here we calculate deviations in $t^S$ and $\bx^S$ from their respective mean values.
Our approach is to use a linear diffusion approximation 
to obtain a stochastic differential equation for the difference between $\by$ and its averaged value,
and analyze first passage to $x_2 = \delta^-$.

We write the slow-fast system (\ref{eq:dhatx3})-(\ref{eq:dhaty2}) as
\begin{eqnarray}
dz(t) &=& \frac{1}{\ee} \left( \left\{ \begin{array}{lc}
a_L(\by(t)) \;, & z(t) < 0 \\
-a_R(\by(t)) \;, & z(t) > 0
\end{array} \right\} + O(\ee) \right) \,dt +
\frac{1}{\sqrt{\ee}} \,e_1^{\sf T} D \,d\bW(t) \;, \label{eq:dhatx4} \\
d\by(t) &=& \big( F_0(z(t),\by(t)) + O(\ee) \big) \,dt +
\sqrt{\ee} \left[ \begin{array}{c}
e_2^{\sf T} \\ \vdots \\ e_N^{\sf T}
\end{array} \right]
D \,d\bW(t) \;, \label{eq:dhaty3}
\end{eqnarray}
where
\begin{equation}
F_0(z,\by) = \frac{b_L(\by)+b_R(\by)}{2} - \frac{b_L(\by)-b_R(\by)}{2} \,{\rm sgn}(z) \;,
\label{eq:F0}
\end{equation}
constitutes the leading order component of $F$ (\ref{eq:F}).
From (\ref{eq:meanFcal}), the averaged value of $F_0$ is
$\mathbb{E} \left[ F_0(z,\by) | \by \right] = \Omega(\by)$.
In view of (\ref{eq:Omega}), the averaged value of $\by(t)$ is $\by_{\dd}(t)$,
and for this reason we define
\begin{equation}
\hat{\by}(t) = \by(t) - \by_{\dd}(t) \;.
\label{eq:yhat}
\end{equation}

In \cite{SiKu14b} we performed an asymptotic expansion of the Fokker-Planck equation
for (\ref{eq:dhatx4})-(\ref{eq:dhaty3}).
Assuming the validity of this expansion, we showed that as $\ee \to 0$ the distribution of
\begin{equation}
\bY(t) = \frac{\hat{\by}(t)}{\sqrt{\ee}} \;,
\label{eq:checkby}
\end{equation}
converges weakly to that of
\begin{equation}
d\bY(t) = (\rD_\by \Omega)(\by_{\dd}(t)) \bY(t) \,dt
+ M(\by_{\dd}(t)) \,d\bW(t) \;,
\label{eq:linearDiffusionApprox2}
\end{equation}
where
\begin{equation}
M(\by) = \left[ -\frac{b_L(\by)-b_R(\by)}{a_L(\by)+a_R(\by)} \,\bigg|\, I \right] D \;,
\label{eq:M}
\end{equation}
and $I$ is the $(N-1)$-dimensional identity matrix.
In order to obtain statistics for $t^S$ and $\bx^S$ by applying
standard first passage theory to (\ref{eq:linearDiffusionApprox2}),
we require strong convergence.
For this reason we use the method of averaging to derive a linear diffusion approximation,
which provides strong convergence \cite{Ki03,Ar03}, and compare it to (\ref{eq:linearDiffusionApprox2}).

Expanding (\ref{eq:dhaty3}) about $\by = \by_{\dd}$ produces\removableFootnote{
We start with (here $r = \frac{t}{\ee}$)
\begin{align}
\frac{1}{\ee} d\hat{\by}(r) &= \frac{1}{\ee} d\by(r) - \frac{1}{\ee} d\by_{\dd}(\ee r) \nonumber \\
&= F_0(z(r),\by(r)) \,dr + \left[ \begin{array}{c}
e_2^{\sf T} \\ \vdots \\ e_N^{\sf T}
\end{array} \right]
D \,d\bW(r) - \Omega(\by_{\dd}(\ee r)) \,dr + O(\ee) \;.
\end{align}
We then substitute
\begin{equation}
F_0(z(r),\by(r)) = F_0(z(r),\by_{\dd}(\ee r))
+ \left( \rD_\by F_0 \right)(z(t),\by_{\dd}(\ee r)) \left( \by(r) - \by_{\dd}(\ee r) \right)
+ O \left( \big| \by - \by_{\dd} \big|^2 \right) \;,
\end{equation}
and note that the error term is $O(\ee)$.
Finally we add and subtract the term
$(\rD_\by \Omega)(\by_{\dd}(\ee r)) \hat{\by}(r)$
and claim that
$\left( \left( \rD_\by F_0 \right)(z(t),\by_{\dd}(\ee r))
- (\rD_\by \Omega)(\by_{\dd}(\ee r)) \right) \hat{\by}(r) = O \left( \sqrt{\ee} \right)$,
although I don't know how to justify this.
}
\begin{equation}
d\hat{\by}(t) = (\rD_\by \Omega)(\by_{\dd}(t)) \hat{\by}(t) \,dt +
\big( F_0(z(t),\by_{\dd}(t)) - \Omega(\by_{\dd}(t)) \big) \,dt
+ \sqrt{\ee} \left[ \begin{array}{c}
e_2^{\sf T} \\ \vdots \\ e_N^{\sf T}
\end{array} \right]
D \,d\bW(t) + O \left( \ee^{\frac{3}{2}} \right) \;.
\label{eq:linearDiffusionApprox0}
\end{equation}
The fast variable $z(t)$ appears in only the middle term of (\ref{eq:linearDiffusionApprox0}).
The essence of the linear diffusion approximation is to replace this term with an equivalent diffusion
\cite{FrWe12,PaSt08,MoCu11,Kh66b}.
Such a computation is beyond the scope of this paper in the general situation
that the noise terms of (\ref{eq:dhatx4}) and (\ref{eq:dhaty3}) are correlated.
In this case the middle term of (\ref{eq:linearDiffusionApprox0})
and the noise term of (\ref{eq:linearDiffusionApprox0}) are not independent
and it seems necessary to study the occupation times of $z(t)$ on
either side of zero, as in \cite{SiKu14b}\removableFootnote{
To describe short time dynamics with correlated noise in the simplest case
it seems that we require knowledge of the joint PDF of $W(r)$ and the positive occupation time of $z(t)$.
On the other hand it is not clear that a ``diffusion'' approximation is even appropriate,
that is it is possible that $\hat{\by}$ cannot be accurately described by an equation of the form
$d \hat{\by} = a(\hat{\by}) \,dt + b(\hat{\by}) \,d\bV(t)$.
Over long times we know that the distribution of $\hat{\by}$ is roughly Gaussian \cite{SiKu14b}.
}.
With the correlation matrix of (\ref{eq:dhatx4})-(\ref{eq:dhaty3})
partitioned as in (\ref{eq:alphaBetaGamma}),
if $\beta = 0$ then the noise terms of (\ref{eq:dhatx4}) and (\ref{eq:dhaty3}) are uncorrelated.
In this case (\ref{eq:linearDiffusionApprox0}) admits the linear diffusion approximation
\begin{equation}
d\hat{\by}(t) = (\rD_\by \Omega)(\by_{\dd}(t)) \hat{\by}(t) \,dt
+ \sigma(\by_{\dd}(t)) \sqrt{\alpha \ee} \,dW(t)
+ \sqrt{\ee} \tilde{D} \,d\bV(t) \;,
\label{eq:linearDiffusionApprox}
\end{equation}
where $W(t)$ is a one-dimensional Brownian motion,
$\bV(t)$ is an $(N-1)$-dimensional Brownian motion independent to $W(t)$,
$\tilde{D} \tilde{D}^{\sf T} = \gamma$, and
\begin{equation}
\sigma \sigma^{\sf T} = \frac{(b_L-b_R) (b_L-b_R)^{\sf T}}{(a_L+a_R)^2} \;.
\label{eq:sigma2}
\end{equation}
The formula (\ref{eq:sigma2}) is derived in Appendix \ref{sec:SIGMA}
by using an explicit expression for the transitional PDF of the leading order
truncation of (\ref{eq:dhatx4}).
The approximation (\ref{eq:linearDiffusionApprox}) is called ``linear''
because the drift term of (\ref{eq:linearDiffusionApprox}) is linear.
The validity of (\ref{eq:linearDiffusionApprox})
requires that the drift term of (\ref{eq:dhatx4}) is Lipschitz in $\by$
and is not influenced by the fact that this drift term is discontinuous in $z$.

The correlation matrices for the noise terms in (\ref{eq:linearDiffusionApprox})
are $\sigma \sigma^{\sf T} \alpha$ and $\gamma$, respectively.
In comparison, the correlation matrix for (\ref{eq:linearDiffusionApprox2}) is
\begin{equation}
M M^{\sf T} = \sigma \sigma^{\sf T} \alpha - \sigma \beta^{\sf T} - \beta \sigma^{\sf T} + \gamma \;.
\label{eq:MMT}
\end{equation}
and so (\ref{eq:linearDiffusionApprox}) is equivalent to (\ref{eq:linearDiffusionApprox2})
when $\beta = 0$.
Given this agreement we conjecture that (\ref{eq:linearDiffusionApprox2}) has strong convergence for any $D$,
with which we may apply standard first passage theory to (\ref{eq:linearDiffusionApprox2})
and obtain statistics for $t^S$ and $\bx^S$.
Indeed the first passage theory
provides a good match to the results of Monte-Carlo simulations of the relay control system
for various choices of $D$ as discussed in the next section. 



Equation (\ref{eq:linearDiffusionApprox2}) is a time-dependent Ornstein-Uhlenbeck process \cite{Sc10,Ga09},
and thus the PDF for $\bY(t)$ is the Gaussian
\begin{equation}
p^S(\bY,t) = \frac{1}{(2 \pi \ee)^{\frac{N-1}{2}} \sqrt{\det(\Theta(t))}}
{\rm e}^{-\frac{1}{2} \bY^{\sf T} \Theta(t)^{-1} \bY} \;,
\label{eq:pycheck}
\end{equation}
where
\begin{equation}
\Theta(t) = \int_0^t
{\rm e}^{\int_s^t (\rD_\by \Omega)(\by_{\dd}(\tilde{s})) \,d\tilde{s}}
M(\by_{\dd}(s)) M(\by_{\dd}(s))^{\sf T}
{\rm e}^{\int_s^t (\rD_\by \Omega)(\by_{\dd}(\tilde{s}))^{\sf T} \,d\tilde{s}} \,ds \;.
\label{eq:Theta}
\end{equation}
Then the leading order terms of ${\rm Var} \left( t^S \right)$ and ${\rm Cov} \left( \by^S \right)$
are found via standard first passage theory \cite{FrWe12} 
(employed in \S\ref{sec:EXCUR}) applied to (\ref{eq:linearDiffusionApprox2}):
\begin{eqnarray}
{\rm Var} \left( t^S \right) &=& \frac{\theta_{11}(t_{\dd}^S)}
{e_1^{\sf T} \Omega(\by_{\dd}(t_{\dd}^S))} \ee + O(\ee^2) \;,
\label{eq:vartsl} \\
{\rm Cov} \left( \by^S \right) &=&
\left( I - \frac{\Omega e_1^{\sf T}}{e_1^{\sf T} \Omega} \right)
\Theta
\left( I - \frac{\Omega e_1^{\sf T}}{e_1^{\sf T} \Omega} \right)^{\sf T}
\bigg|_{\by = \by_{\dd}(t_{\dd}^S)} \ee + O(\ee^2) \;,
\label{eq:covbysl}
\end{eqnarray}
where $\theta_{11}$ is the top left entry of $\Theta$.
Lastly, since $x_1$ operates on a fast time-scale relative to $\by$,
to leading order, $x_1^S$ is uncorrelated to $\by^S$.
From (\ref{eq:pqss}), we have
\begin{equation}
{\rm Var}(x_1^S) = 
\frac{a_L^2 + a_R^2}{4 a_L^2 a_R^2} \bigg|_{\by = \by_{\dd}(t_{\dd}^S)} \ee^2
+ O(\ee^3) \;.
\label{eq:varx1sl}
\end{equation}

\subsection{Summary and comparison to numerical simulations}

\begin{figure}[b!]
\begin{center}
\setlength{\unitlength}{1cm}
\begin{picture}(15,13)
\put(3.75,7){\includegraphics[height=6cm]{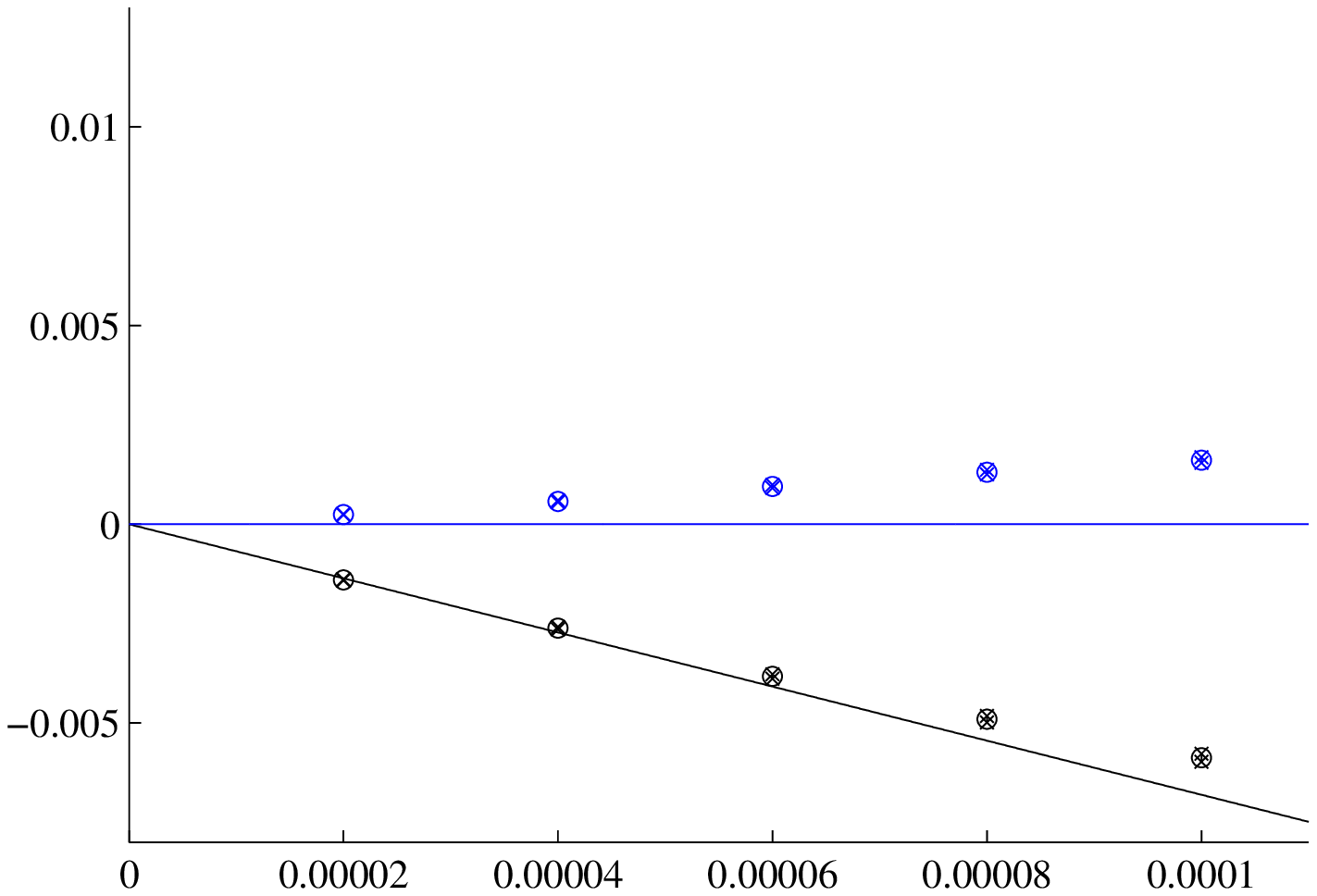}}
\put(0,0){\includegraphics[height=6cm]{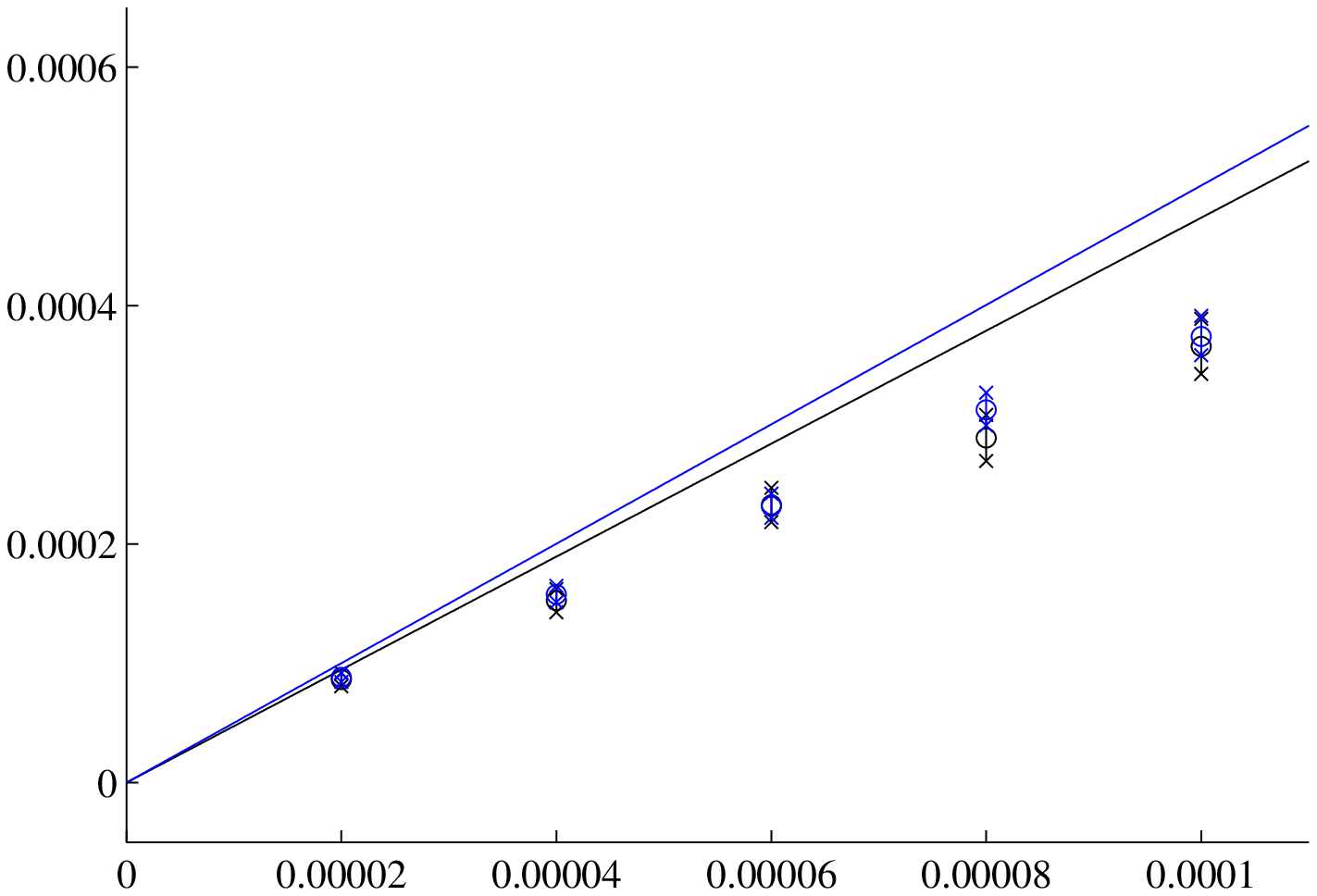}}
\put(7.5,0){\includegraphics[height=6cm]{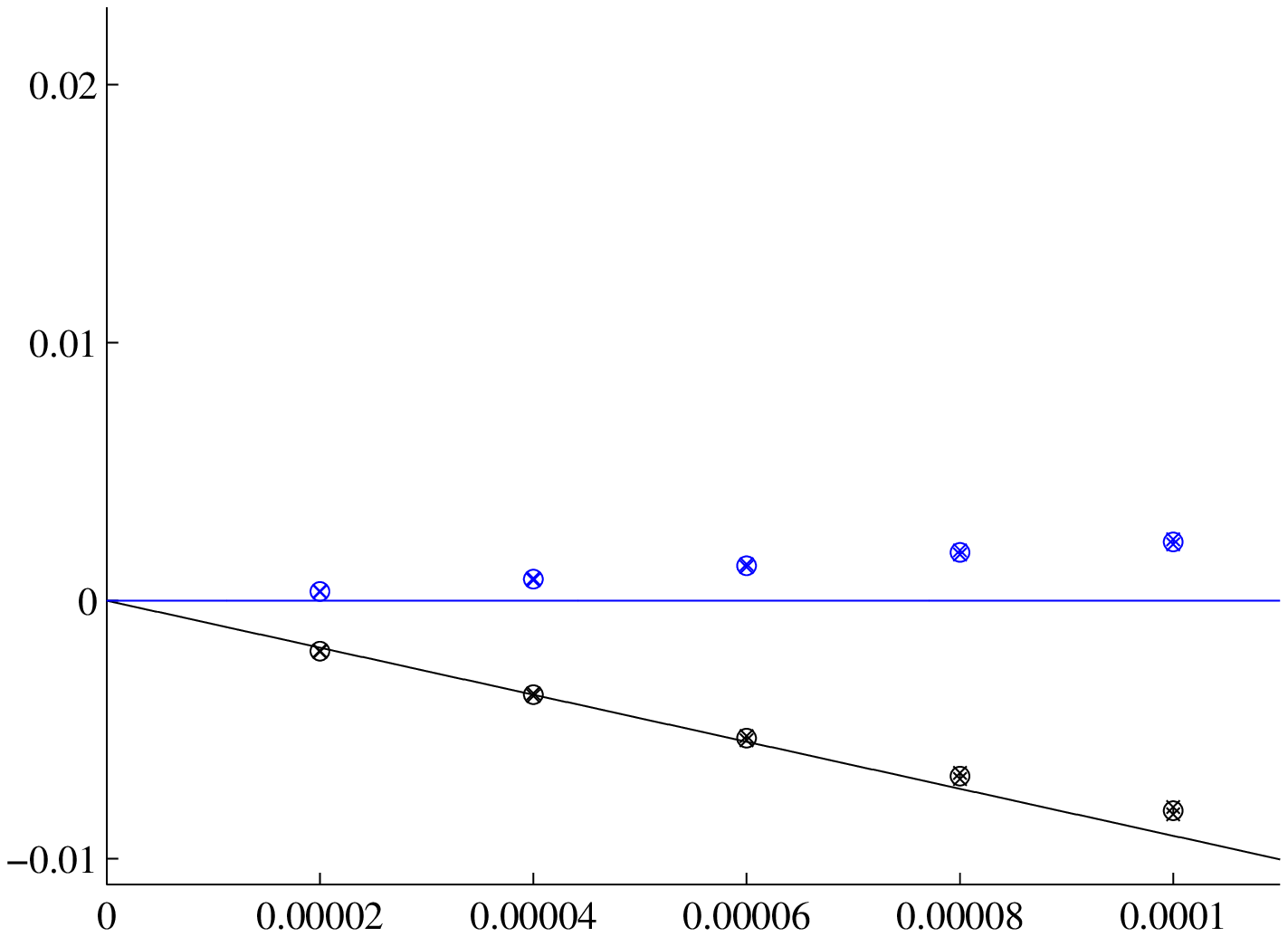}}
\put(4.75,12.8){\large \sf \bfseries A}
\put(1,5.8){\large \sf \bfseries B}
\put(8.5,5.8){\large \sf \bfseries C}
\put(7.55,7){\small $\ee$}
\put(9.5,8.6){\small ${\rm Diff} \left( t^S \right)$}
\put(6.5,10){\small \color{blue} ${\rm Std} \left( t^S \right)$}
\put(3.8,0){\small $\ee$}
\put(5,2.4){\small ${\rm Diff} \left( x_1^S \right)$}
\put(2.7,3.2){\small \color{blue} ${\rm Std} \left( x_1^S \right)$}
\put(11.3,0){\small $\ee$}
\put(12.5,1.66){\small ${\rm Diff} \left( x_3^S \right)$}
\put(10.2,2.74){\small \color{blue} ${\rm Std} \left( x_3^S \right)$}
\end{picture}
\caption{
A comparison of first passage statistics for stochastically perturbed sliding motion for
the relay control example, (\ref{eq:relayControlSystem2}) with
(\ref{eq:ABCDvalues3d})-(\ref{eq:paramValues}), with theoretical results.
Here $\bx_0 = \bx_{\Gamma}^M$,
and we have used $\delta^- = 0.1$, as discussed in \S\ref{sec:COMB}.
The deterministic passage time and location
are $t_{\dd}^S \approx 1.032$ and $\bx_{\dd}^S \approx (0,-0.1,-0.000685)$.
The data points were computed from $1000$ Monte-Carlo simulations for each value of $\ee$
using the Euler-Maruyama method with a fixed step size, $\Delta t = 0.00001$.
The circles, bars and crosses indicate mean values and 95\% confidence intervals, as in Fig.~\ref{fig:manyPeriod}.
The solid curves are the theoretical predictions.
${\rm Std}(t^S)$ is approximated by (\ref{eq:vartsl}),
${\rm Std}(x_1^S)$ by (\ref{eq:varx1sl}),
and ${\rm Std}(x_1^S)$ by (\ref{eq:covbysl}).
${\rm Diff}(t^S)$ is approximated by (\ref{eq:meantsl}),
${\rm Diff}(x_1^S)$ by (\ref{eq:meanx1sl})
and ${\rm Diff}(x_3^S)$ by (\ref{eq:byBarDef}) and (\ref{eq:meanbysl1}).
\label{fig:checkSlide}
}
\end{center}
\end{figure}

\begin{figure}[b!]
\begin{center}
\setlength{\unitlength}{1cm}
\begin{picture}(15,13)
\put(3.75,7){\includegraphics[height=6cm]{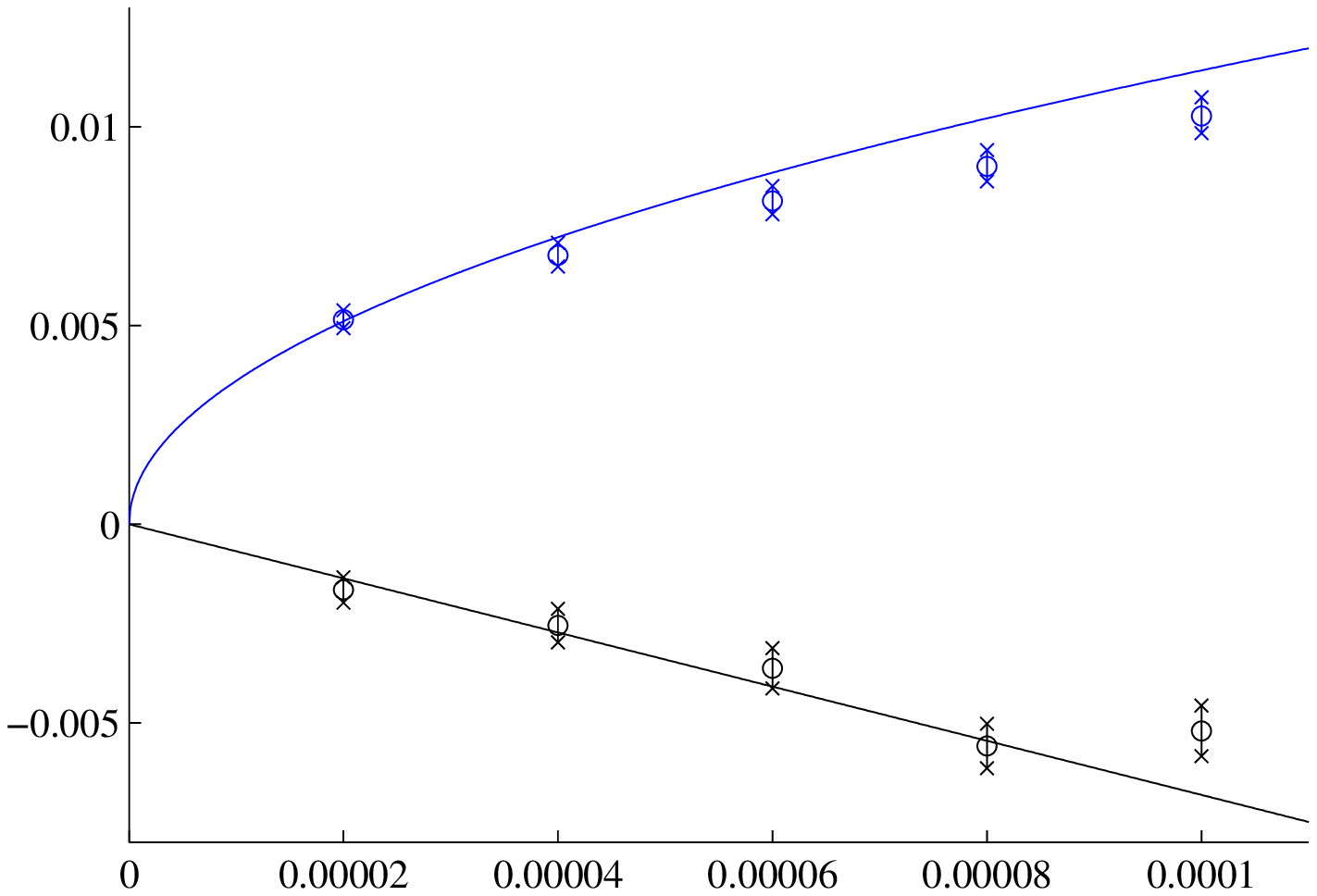}}
\put(0,0){\includegraphics[height=6cm]{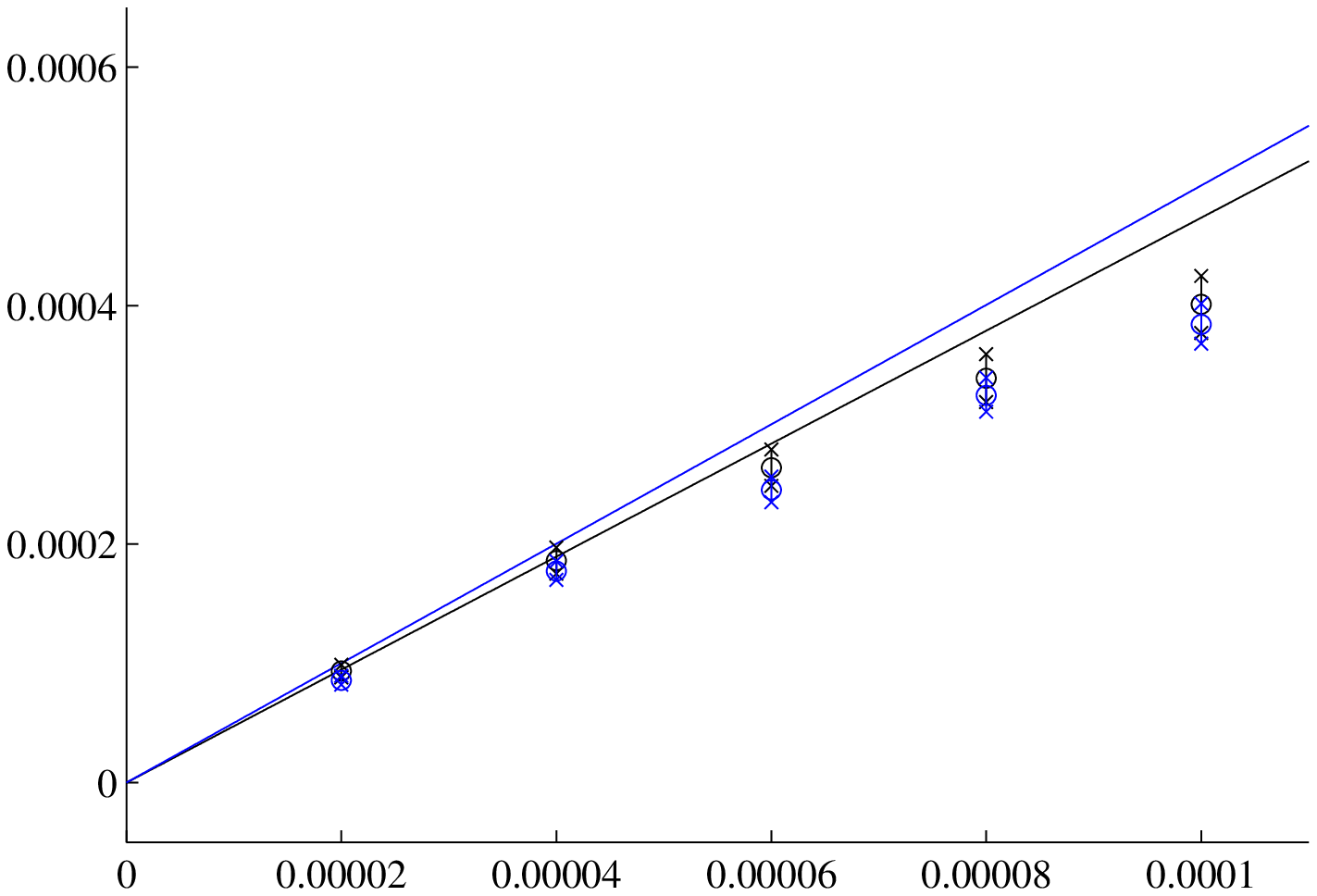}}
\put(7.5,0){\includegraphics[height=6cm]{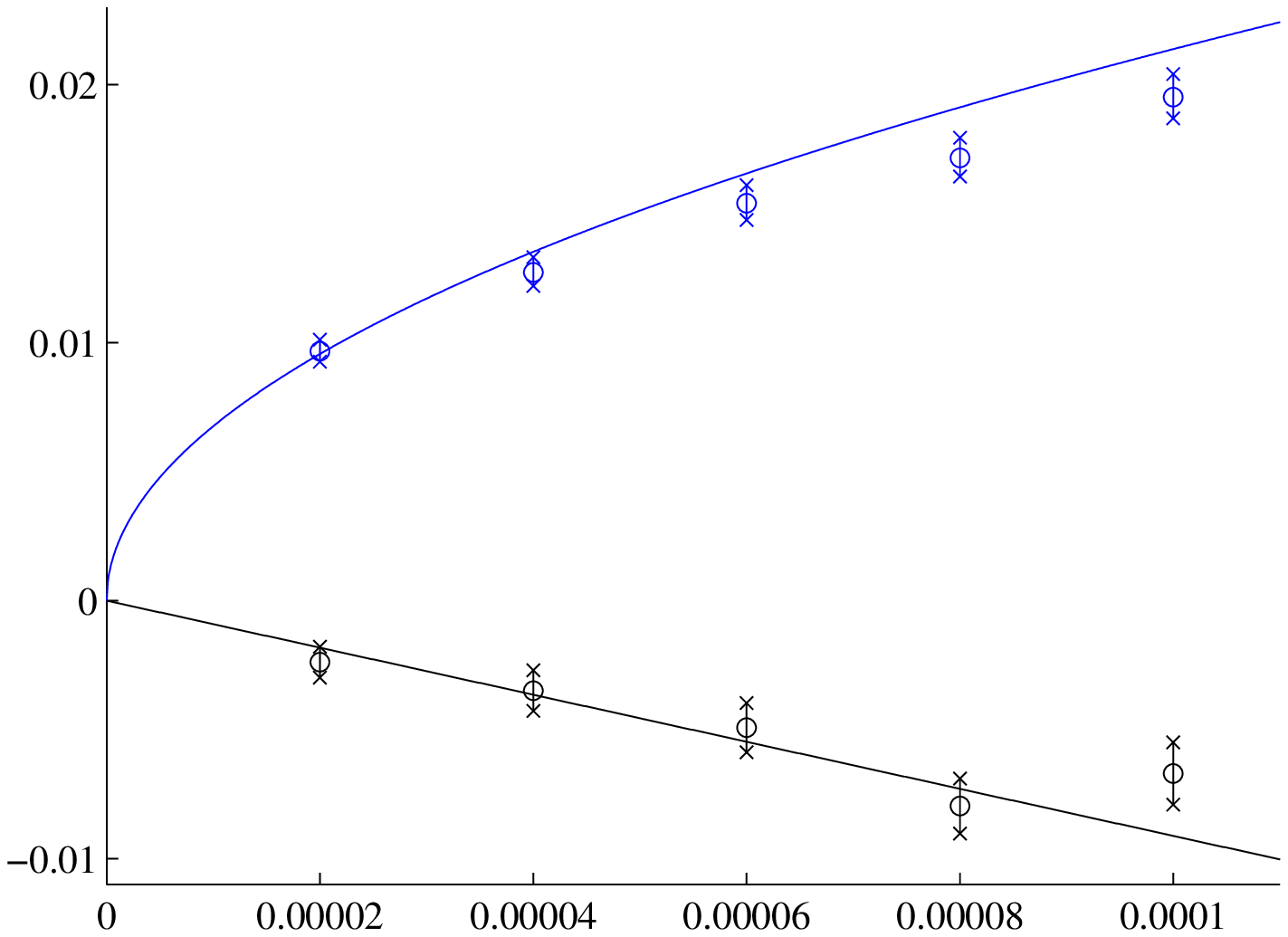}}
\put(4.75,12.8){\large \sf \bfseries A}
\put(1,5.8){\large \sf \bfseries B}
\put(8.5,5.8){\large \sf \bfseries C}
\put(7.55,7){\small $\ee$}
\put(6.5,8.4){\small ${\rm Diff} \left( t^S \right)$}
\put(8.5,12.2){\small \color{blue} ${\rm Std} \left( t^S \right)$}
\put(3.8,0){\small $\ee$}
\put(5,2.5){\small ${\rm Diff} \left( x_1^S \right)$}
\put(2.7,3.3){\small \color{blue} ${\rm Std} \left( x_1^S \right)$}
\put(11.3,0){\small $\ee$}
\put(12.5,1.66){\small ${\rm Diff} \left( x_3^S \right)$}
\put(10.9,5.3){\small \color{blue} ${\rm Std} \left( x_3^S \right)$}
\end{picture}
\caption{
First passage statistics of the sliding phase for the relay control example
as in Fig.~\ref{fig:checkSlide}, except for this figure we have used
$D = e_1 e_1^{\sf T}$ in place of $D = B e_1^{\sf T}$ in (\ref{eq:ABCDvalues3d}).
\label{fig:checkSlide2}
}
\end{center}
\end{figure}

Fig.~\ref{fig:checkSlide} compares the above theoretical results
with Monte-Carlo simulations of the relay control example
from $\bx_0 = \bx_{\Gamma}^M$ to $x_2 = \delta^-$.
For panel A, ${\rm Diff}(t^S)$ is approximated by $\ee t^{S,1}$ using (\ref{eq:meantsl}),
and ${\rm Std}(t^S)$ is approximated using (\ref{eq:vartsl}).
For panel B, ${\rm Diff}(x_1^S)$ is approximated by (\ref{eq:meanx1sl}),
and ${\rm Std}(x_1^S)$ by (\ref{eq:varx1sl}).
Lastly for panel C, ${\rm Diff}(x_3^S)$ is approximated
with (\ref{eq:byBarDef}) and (\ref{eq:meanbysl1}),
and the standard deviation with (\ref{eq:covbysl}).

Notice that the approximations to the standard deviations of $t^S$ and $x_3^S$ are zero.
This is because, remarkably, for our example the matrix $M M^{\sf T}$ is identically zero, and thus so is $\Theta$.
To see this why this is the case,
we notice that the left and right half-systems of (\ref{eq:relayControlSystem5})
are identical up to a constant vector, that is
\begin{equation}
a_L + a_R = e_1^{\sf T}
\left( \mathcal{B}^{(L)} - \mathcal{B}^{(R)} \right) \;, \qquad
b_L - b_R =
\left[ \begin{array}{c} e_2^{\sf T} \\ \vdots \\ e_N^{\sf T} \end{array} \right]
\left( \mathcal{B}^{(L)} - \mathcal{B}^{(R)} \right) \;,
\label{eq:aLplusaRbLplusbR}
\end{equation}
are constant.
Moreover, in view of (\ref{eq:calBLBRD}) we can write
\begin{equation}
D = \frac{1}{2} \left( \mathcal{B}^{(L)} - \mathcal{B}^{(R)} \right) e_1^{\sf T} \;.
\label{eq:calD2}
\end{equation}
By substituting (\ref{eq:sigma2}), (\ref{eq:aLplusaRbLplusbR}) and (\ref{eq:calD2}) into (\ref{eq:MMT}),
we immediately obtain $M M^{\sf T} = 0$.

This can also be explained geometrically.
The drift term of the stochastic differential equation (\ref{eq:sde2})
is piecewise, and in general points in different unrelated directions for $x_1 < 0$ and $x_1 > 0$.
When we rewrite (\ref{eq:sde2}) in terms of the variables
$x_1$ and $\hat{\by} = \by - \by_{\dd}(t)$
(representing deviations from the deterministic solution),
to lowest order the drift for $x_1 > 0$ is a scalar multiple of the drift for $x_1 < 0$.
For our example (\ref{eq:relayControlSystem2}),
both left and right drift vectors are multiples of $B$
and the noise term in (\ref{eq:relayControlSystem2})
is one-dimensional and also a multiple $B$.
Consequently, to leading order deviations occur along a line in the direction of $B$.
Deviations in any direction orthogonal to the switching manifold are $O(\ee)$,
hence in this case overall deviations are $O(\ee)$.
For this reason, for our example, the $O(\sqrt{\ee})$ term
representing deviations in $y_i$ is zero.

In Fig.~\ref{fig:checkSlide2} we repeat the numerical comparison
using $D = e_1 e_1^{\sf T}$ in place of $D = B e_1^{\sf T}$ in (\ref{eq:ABCDvalues3d})\removableFootnote{
See {\sc drawSlide3b.m}.
}.
Now $M M^{\sf T}$ is nonzero and
${\rm Std}(t^S)$ and ${\rm Std}(x_3^S)$ are $O(\sqrt{\ee})$.
Again the theoretical calculations are consistent with the numerical simulations.
We have found that the first passage predictions match the results of Monte-Carlo simulations
for other choices of $D$, including those with $\beta \ne 0$\removableFootnote{
I performed the experiment using, $D = [1,0,1]^{\sf T} e_1^{\sf T}$,
for which the noise terms are correlated,
and found that the first passage theory predictions still seem to work, see {\sc drawSlide4.m}.
}.

\section{Escaping analysis}
\label{sec:ESCAPE}
\setcounter{equation}{0}

In this section we study (\ref{eq:sde}) near the origin
and in the range $\delta^- \le x_2 \le \delta^+$.
As in the previous section it is convenient to write
$\by = [ x_2, \ldots, x_N ]^{\sf T}$ and, repeating (\ref{eq:sde3}), write (\ref{eq:sde}) as
\begin{equation}
\left[ \begin{array}{c} dx_1(t) \\ d\by(t) \end{array} \right] =
\left\{ \begin{array}{lc}
\left[ \begin{array}{c}
a_L(\by(t)) + c_L(\by(t)) x_1(t) + O(x_1(t)^2) \\
b_L(\by(t)) + d_L(\by(t)) x_1(t) + O(x_1(t)^2)
\end{array} \right]
\;, & x_1(t) < 0 \\
\left[ \begin{array}{c}
-a_R(\by(t)) + c_R(\by(t)) x_1(t) + O(x_1(t)^2) \\
b_R(\by(t)) + d_R(\by(t)) x_1(t) + O(x_1(t)^2)
\end{array} \right]
\;, & x_1(t) > 0
\end{array} \right\} \,dt + \sqrt{\ee} D \,d\bW(t) \;.
\label{eq:sde4}
\end{equation}
Here we expand the coefficients in (\ref{eq:sde4}) about $\by = 0$.
In view of assumptions (\ref{eq:a1})-(\ref{eq:a5}), we can write\removableFootnote{
Specifically,
(i) (\ref{eq:a2}) implies the first equation and $a_L > 0$,
(ii) (\ref{eq:a1}) and (\ref{eq:a5}) together imply the second equation,
(iii) (\ref{eq:a4}) implies $\frac{\partial a_R}{\partial x_2} > 0$,
(iv) (\ref{eq:a3}) implies the third equation and $b_{R1} > 0$, and 
(v) (\ref{eq:a3}) implies the fourth equation.
}
\begin{equation}
\begin{gathered}
a_L(\by) = a_L + O(||\by||) \;, \qquad
a_R(\by) = \frac{\partial a_R}{\partial x_2} x_2 + O(||\by||^2) \;, \\
b_{R1}(\by) = b_{R1} + O(||\by||) \;, \qquad
b_{Ri}(\by) = \sum_{j=2}^N \frac{\partial b_{Ri}}{\partial x_j} x_j + O(||\by||^2) \;,
\forall i \ne 1 \;, \\
c_R(\by) = c_R + O(||\by||) \;, \qquad
d_R(\by) = d_R + O(||\by||) \;,
\end{gathered}
\end{equation}
where, on the right hand sides, and in the remainder of this section, the coefficients are evaluated at $\by = 0$, and
\begin{equation}
a_L > 0 \;, \qquad
\frac{\partial a_R}{\partial x_2} < 0 \;, \qquad
b_{R1} > 0 \;.
\end{equation}
To study dynamics near the origin asymptotically in $\ee$, we scale space and time.
Consider the general scaling\removableFootnote{
Note, \LaTeX~can't render $\mathcal{Y}$ bold!
}
\begin{equation}
\mathcal{X}_1 = \frac{x_1}{\ee^{\lambda_1}} \;, \qquad
\mathcal{X}_i = \frac{x_i}{\ee^{\lambda_2}} \;, \forall i \ne 1 \;, \qquad
\mathcal{T} = \frac{t}{\ee^{\lambda_3}} \;,
\label{eq:tildes}
\end{equation}
where $\lambda_1, \lambda_2, \lambda_3 > 0$.
By substituting (\ref{eq:tildes}) into (\ref{eq:sde4}), for $\mathcal{X}_1 > 0$ we obtain\removableFootnote{
Originally I had this in terms of a PDE, but there are several reasons why I prefer the SDE.
\begin{enumerate}
\item
Scaling in the PDE appears more rigorous but neither way is properly rigorous.
To validate the asymptotic expansion one could try to apply the maximum principle,
but this is very difficult and well beyond the scope of this paper.
\item
In the context of the SDE the scaling seems easier to understand:
The values of $\lambda_1$, $\lambda_2$, $\lambda_3$ are such that the following three terms are leading order:
(i) noise in $\mathcal{X}_1$,
(ii) drift in $\mathcal{X}_1$, 
and (iii) drift in $\mathcal{X}_2$.
\item
By doing the scaling in terms of a PDE,
one ends up with the desired BVP (\ref{eq:fpeu})-(\ref{eq:icu}),
but an extra leap is required to extract back the corresponding SDE (\ref{eq:du}).
Strictly speaking this SDE is not really necessary here,
but it does significantly add to our understanding of what is going on.
\item
If the scaling is done in the context of the elliptic BVP corresponding to escape from a
neighbourhood of the origin,
then since the PDE is independent of time,
the analysis gives us $\lambda_1$ and $\lambda_2$, but not $\lambda_3$.
\end{enumerate}
},
\begin{equation}
\begin{split}
d\mathcal{X}_1(\mathcal{T}) &=
\left( -\ee^{\lambda_2+\lambda_3-\lambda_1}
{\textstyle \frac{\partial a_R}{\partial x_2}} \mathcal{X}_2(\mathcal{T})
+ \ee^{\lambda_3} c_R \mathcal{X}_1(\mathcal{T}) +
O \left( \ee^{\lambda_3-\lambda_1+2{\rm min}(\lambda_1,\lambda_2)} \right) \right) \,d\mathcal{T} \\
&+ \ee^{\frac{\lambda_3+1}{2} - \lambda_1} e_1^{\sf T} D \,d\bW(\mathcal{T}) \;, \\
d\mathcal{X}_2(\mathcal{T}) &= \left( \ee^{\lambda_3-\lambda_2} b_{R1} +
O \left( \ee^{\lambda_3-\lambda_2+{\rm min}(\lambda_1,\lambda_2)} \right) \right) \,d\mathcal{T}
+ \ee^{\frac{\lambda_3+1}{2} - \lambda_2} e_2^{\sf T} D \,d\bW(\mathcal{T}) \;, \\
d\mathcal{X}_i(\mathcal{T}) &=
O \left( \ee^{\lambda_3-\lambda_2+{\rm min}(\lambda_1,\lambda_2)} \right) \,d\mathcal{T}
+ \ee^{\frac{\lambda_3+1}{2} - \lambda_2} e_{i+1}^{\sf T} D \,d\bW(\mathcal{T}) \;,
\quad \forall i \ge 3 \;,
\end{split}
\label{eq:sdetildeR}
\end{equation}
and for $\mathcal{X}_1 < 0$,
\begin{equation}
\begin{split}
d\mathcal{X}_1(\mathcal{T}) &= \left( \ee^{\lambda_3-\lambda_1} a_L +
O \left( \ee^{\lambda_3-\lambda_1+{\rm min}(\lambda_1,\lambda_2)} \right) \right) \,d\mathcal{T}
+ \ee^{\frac{\lambda_3+1}{2} - \lambda_1} e_1^{\sf T} D \,d\bW(\mathcal{T}) \;, \\
d\mathcal{X}_i(\mathcal{T}) &= \left( \ee^{\lambda_3-\lambda_2} b_{Li} +
O \left( \ee^{\lambda_3-\lambda_2+{\rm min}(\lambda_1,\lambda_2)} \right) \right) \,d\mathcal{T}
+ \ee^{\frac{\lambda_3+1}{2} - \lambda_2} e_{i+1}^{\sf T} D \,d\bW(\mathcal{T}) \;,
\quad \forall i \ge 2 \;.
\end{split}
\label{eq:sdetildeL}
\end{equation}

To identify the appropriate choice for each $\lambda_j$, $j=1,2,3$,
one would normally consider the asymptotic behaviour as $\ee \to 0$
of the Fokker-Planck equation for the joint probability density of $\mathcal{X}_i$ for all $i$.
We do not provide such calculations here,
and instead consider the stochastic differential equation directly for the sake of brevity.
Both approaches lead to the same conclusions.

First note that the ratio of the drift in the $\mathcal{X}_1$-direction for $\mathcal{X}_1 > 0$,
to the drift in the $\mathcal{X}_1$-direction for $\mathcal{X}_1 < 0$, approaches zero as $\ee \to 0$,
and for $\mathcal{X}_1 < 0$ this drift is directed to the right.
Therefore we expect sample solutions to be located almost entirely in the right half-space.
We then choose $\lambda_1$, $\lambda_2$ and $\lambda_3$ such that three
terms in (\ref{eq:sdetildeR}) are $O(1)$ and all other terms are of higher order.
This gives two possibilities.
First, we may have $\lambda_1 = \lambda_2 = \lambda_3 = 1$,
but this proves unhelpful because in this case the drift in $\mathcal{X}_1$ is a higher order term
and so this scaling does not capture dynamics escaping a neighbourhood of the switching manifold,
which is what we are trying to describe.
This suggests we need to look on a longer time-scale, that is, we should choose $\lambda_3 < 1$.
Indeed the second possibility is:
\begin{equation}
\lambda_1 = \frac{2}{3} \;, \qquad
\lambda_2 = \frac{1}{3} \;, \qquad
\lambda_3 = \frac{1}{3} \;.
\label{eq:chirholambda}
\end{equation}
Then for $\mathcal{X}_1 < 0$,
\begin{equation}
d\mathcal{X}_1(\mathcal{T}) = \frac{1}{\ee^{\frac{1}{3}}} \,a_L \,d\mathcal{T} + O(\ee^0) \;.
\label{eq:dXEscapeLeft}
\end{equation}
Consideration of the corresponding Fokker-Planck equation for $\mathcal{X}_1 < 0$
leads to the conclusion that the probability that $\mathcal{X}_1 < 0$ is negligible as $\ee \to 0$.
Consequently we consider the dynamics for $\mathcal{X}_1> 0$ only,
providing an appropriate boundary condition at $\mathcal{X}_1 =0$ below in (\ref{eq:bc0u}).

With (\ref{eq:sdetildeR}) and (\ref{eq:chirholambda}), for $\mathcal{X}_1 > 0$
\begin{equation}
\begin{split}
d\mathcal{X}_1(\mathcal{T}) &= -\frac{\partial a_R}{\partial x_2} \mathcal{X}_2(\mathcal{T}) \,d\mathcal{T}
+ e_1^{\sf T} D \,d\bW(\mathcal{T}) 
+ O \left( \ee^{\frac{1}{3}} \right) \;, \\ 
d\mathcal{X}_2(\mathcal{T}) &= b_{R1} \,d\mathcal{T}
+ O \left( \ee^{\frac{1}{3}} \right) \;, \\ 
d\mathcal{X}_i(\mathcal{T}) &= O \left( \ee^{\frac{1}{3}} \right) \;,
\quad \forall i \ge 3 \;.
\end{split}
\label{eq:sdetildeR2}
\end{equation}
As $\ee \to 0$, (\ref{eq:sdetildeR2}) approaches
\begin{equation}
\begin{split}
d\mathcal{X}_1(\mathcal{T}) &= -\frac{\partial a_R}{\partial x_2} \mathcal{X}_2(\mathcal{T}) \,d\mathcal{T}
+ \sqrt{\alpha} \,dW(\mathcal{T}) \;, \\
d\mathcal{X}_2(\mathcal{T}) &= b_{R1} \,d\mathcal{T} \;, \\
d\mathcal{X}_i(\mathcal{T}) &= 0 \;, \quad \forall i \ge 3 \;,
\end{split}
\label{eq:sdetildeR3}
\end{equation}
where $W(\mathcal{T})$ is a scalar Brownian motion,
and $\alpha = (D D^{\sf T})_{11}$ (\ref{eq:alphaBetaGamma}).
We now perform an additional scaling to simplify (\ref{eq:sdetildeR3}).
Note that $\mathcal{X}_2(\mathcal{T})$, as governed by the limiting equation (\ref{eq:sdetildeR3}), is deterministic.
Let $\mathcal{T}_0$ be the time at which $\mathcal{X}_2 = 0$.
Then, with
\begin{equation}
u = \frac{\left| \frac{\partial a_R}{\partial x_2} \right|^{\frac{1}{3}} b_{R1}^{\frac{1}{3}}}
{\alpha^{\frac{4}{3}}} \,\mathcal{X}_1 \;, \qquad
s = \frac{\left| \frac{\partial a_R}{\partial x_2} \right|^{\frac{2}{3}} b_{R1}^{\frac{2}{3}}}
{\alpha^{\frac{2}{3}}} \left( \mathcal{T} - \mathcal{T}_0 \right) \;,
\label{eq:us}
\end{equation}
(\ref{eq:sdetildeR3}) reduces to\removableFootnote{
The appropriate scaling of $\mathcal{X}_2$ is
$\tilde{Y} = \frac{\left| \frac{\partial a_R}{\partial x_2} \right|^{\frac{2}{3}}}
{b_{R1}^{\frac{1}{3}} \alpha^{\frac{2}{3}}} \,\mathcal{X}_2$,
with which we have $\tilde{Y}(s) = s$.
}
\begin{equation}
du(s) = s \,ds + dW(s) \;.
\label{eq:du}
\end{equation}

For the purposes of describing escaping dynamics,
we determine the transitional PDF for (\ref{eq:du}), call it $p^E(u,s)$,
by writing it as the solution to a boundary value problem.
The PDF satisfies the Fokker-Planck equation
\begin{equation}
p^E_s = -s p^E_u + \frac{1}{2} p^E_{uu} \;.
\label{eq:fpeu}
\end{equation}
To produce meaningful solutions we impose the following boundary condition at infinity:
\begin{equation}
p^E(u,s) \to 0 {\rm ~as~} u \to \infty \;.
\end{equation}
At $u=0$ we use the boundary condition
\begin{equation}
s p^E(0,s) - \frac{1}{2} p^E_u(0,s) = 0 \;,
\label{eq:bc0u}
\end{equation}
which may be justified in two ways.
First, the requirement that $u > 0$ with probability $1$ is equivalent to
$\frac{\partial}{\partial s} \int_0^\infty p^E(u,s) \,du = 0$,
and applying this identity to (\ref{eq:fpeu}) produces (\ref{eq:bc0u}).
Second, by (\ref{eq:dXEscapeLeft}) dynamics for $\mathcal{X}_1 < 0$
has drift directed to the right, which in the limit $\ee \to 0$ is infinitely large.
Therefore sample solutions to (\ref{eq:du}) that reach $u=0$ are reflected back to the right;
indeed (\ref{eq:bc0u}) is a reflecting boundary condition \cite{Sc10,Ga09}.
Also we suppose $u(s_0) = u_0$ at some initial time $s_0$, which corresponds to the initial condition
\begin{equation}
p^E(u,s_0)
= \delta(u - u_0) \;.
\label{eq:icu}
\end{equation}

\begin{figure}[t!]
\begin{center}
\setlength{\unitlength}{1cm}
\begin{picture}(7.2,6.2)
\put(0,.2){\includegraphics[height=6cm]{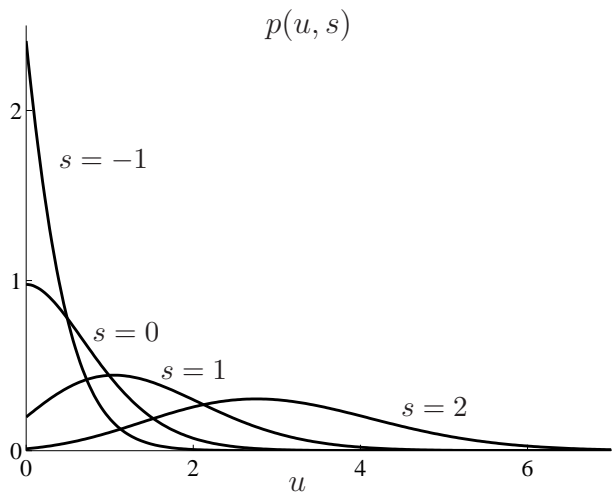}}
\put(3.7,0){$u$}
\put(3.4,6.1){$p(u,s)$}
\put(.65,4.3){\small $s = -1$}
\put(1.1,2){\small $s = 0$}
\put(2,1.5){\small $s = 1$}
\put(5.2,1){\small $s = 2$}
\end{picture}
\caption{
The PDF, (\ref{eq:Kn00}), for four different values of the time, $s$.
The quantity, $u$, represents distance from $x_1=0$ (\ref{eq:us}).
\label{fig:Kn00}
}
\end{center}
\end{figure}

An explicit expression for the solution to (\ref{eq:fpeu})-(\ref{eq:icu}) is derived in \cite{KnYa01}
(by scaling $p^E$ in order to remove the
time-dependency in the coefficients of (\ref{eq:fpeu}) and taking Laplace transforms)
but takes a rather complicated form.
For our purposes it is useful to take $s_0 \to -\infty$,
because $s_0 \propto \frac{\delta^-}{\ee^{\frac{1}{3}}} \to -\infty$, as $\ee \to 0$,
when taking an initial point with $x_2 = \delta^-$ in the original equation (\ref{eq:sde4}).
Also, it is reasonable to take $u_0 \to 0$
because, as shown in \S\ref{sec:SLIDE}, $x_1(t)$ is a fast variable
and with high probability repeatedly intersects the switching manifold
as the sample solution approaches the escaping phase.
As shown by Knessl \cite{Kn00},
the solution to (\ref{eq:fpeu})-(\ref{eq:icu}) with $s_0 \to -\infty$ and $u_0 \to 0$
is given by
\begin{equation}
p^E(u,s) = 2^{\frac{2}{3}} {\rm e}^{-\frac{s^3}{6}} {\rm e}^{u s} Y(u,s) \;,
\label{eq:Kn00}
\end{equation}
where $Y$ is given by the inverse Laplace transform\removableFootnote{
The integral is over some Bromwich contour
(here the imaginary axis is a suitable contour
because all singularities of $\frac{{\rm Ai} \left( 2^{\frac{1}{3}} (u+\nu) \right)}
{{\rm Ai} \left( 2^{\frac{1}{3}} \nu \right)^2}$ have negative real part,
and this is what we use for numerical evaluation of $Y(u,s)$).
The integral converges because $\frac{{\rm Ai} \left( 2^{\frac{1}{3}} (u+{\rm i}\beta) \right)}
{{\rm Ai} \left( 2^{\frac{1}{3}} {\rm i}\beta \right)^2} \to 0$ sufficiently quickly as $\beta \to \infty$.
However, since $\frac{{\rm Ai} \left( 2^{\frac{1}{3}} (u+\nu) \right)}
{{\rm Ai} \left( 2^{\frac{1}{3}} \nu \right)^2} \to \infty$ as $\nu \to \infty$,
this function cannot be the Laplace transform of anything!
Thus, oddly, the Laplace transform of $Y(u,s)$ is equal to something other than
$\frac{{\rm Ai} \left( 2^{\frac{1}{3}} (u+\nu) \right)}
{{\rm Ai} \left( 2^{\frac{1}{3}} \nu \right)^2}$.
}
\begin{equation}
Y(u,s) = \frac{1}{2 \pi {\rm i}} \int_{\rm Br}
\frac{{\rm Ai} \left( 2^{\frac{1}{3}} (u+\nu) \right)}
{{\rm Ai} \left( 2^{\frac{1}{3}} \nu \right)^2} \,{\rm e}^{\nu s} \,d\nu \;,
\end{equation}
and ${\rm Ai}(u)$ is the Airy function\removableFootnote{
$Ai(u) = \frac{1}{\pi} \int_0^\infty {\rm cos} \left( \frac{1}{3} s^3 + u s \right) \,ds$
}.
The PDF (\ref{eq:Kn00}) is shown in Fig.~\ref{fig:Kn00}.

\begin{figure}[t!]
\begin{center}
\setlength{\unitlength}{1cm}
\begin{picture}(15,13)
\put(3.75,7){\includegraphics[height=6cm]{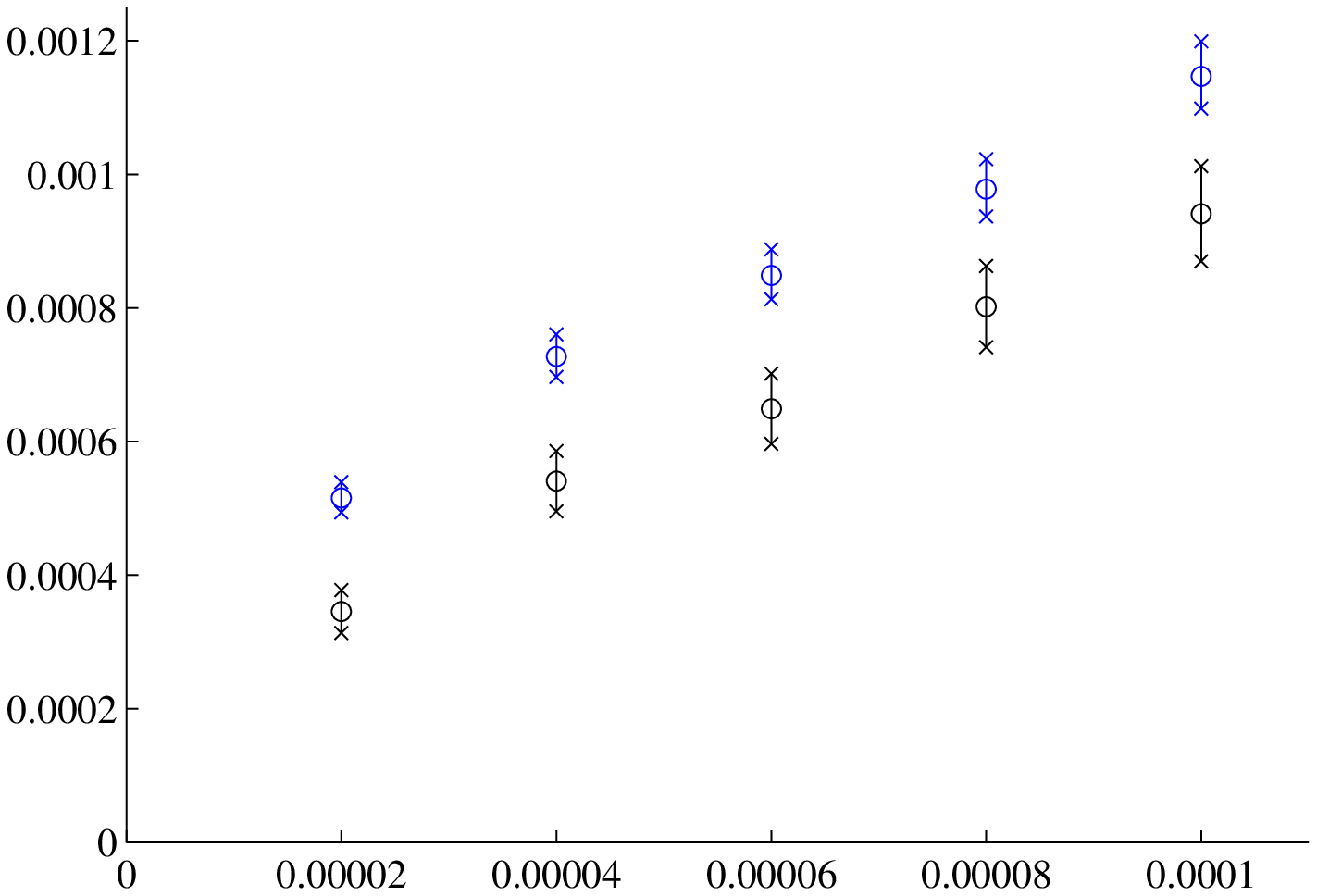}}
\put(0,0){\includegraphics[height=6cm]{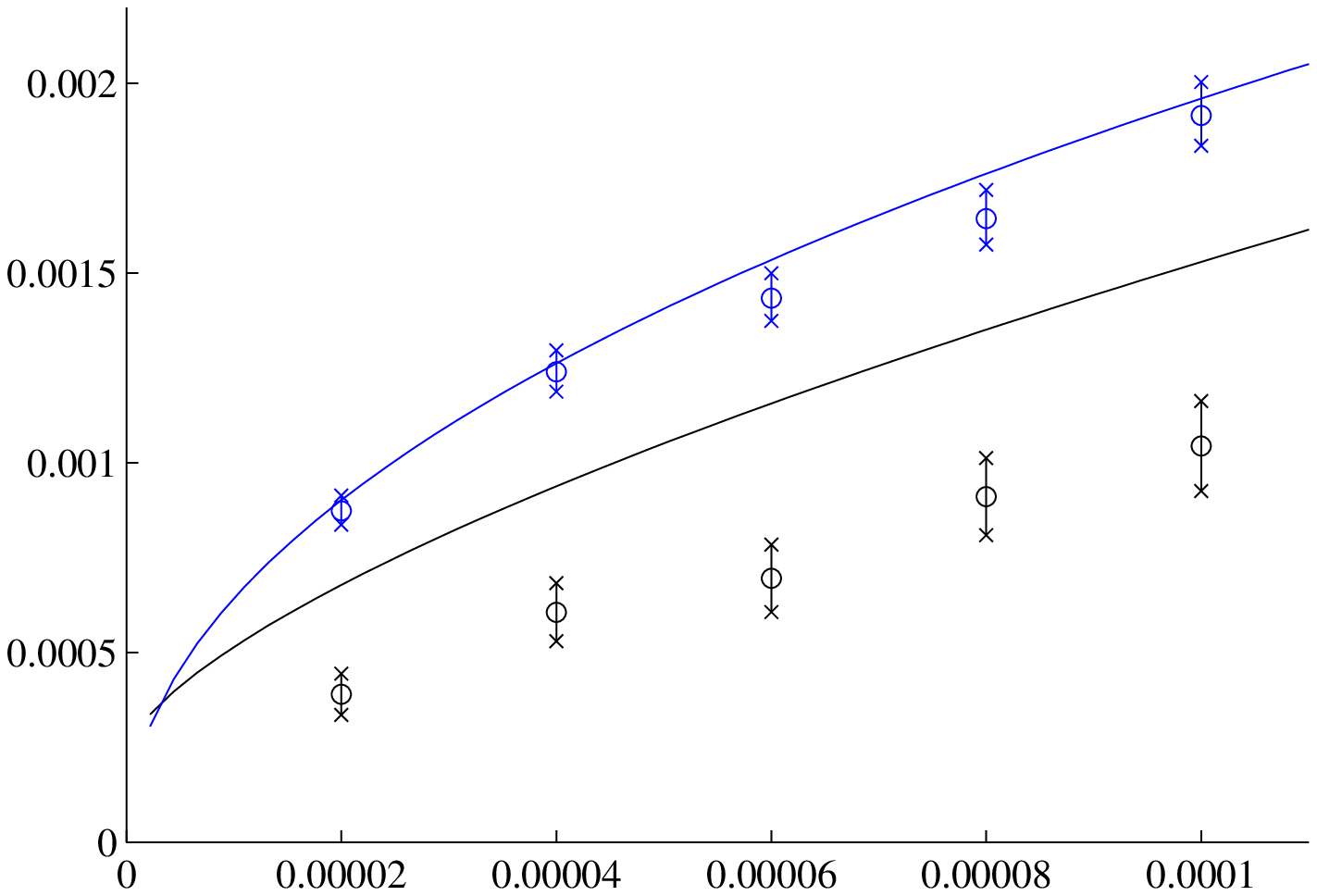}}
\put(7.5,0){\includegraphics[height=6cm]{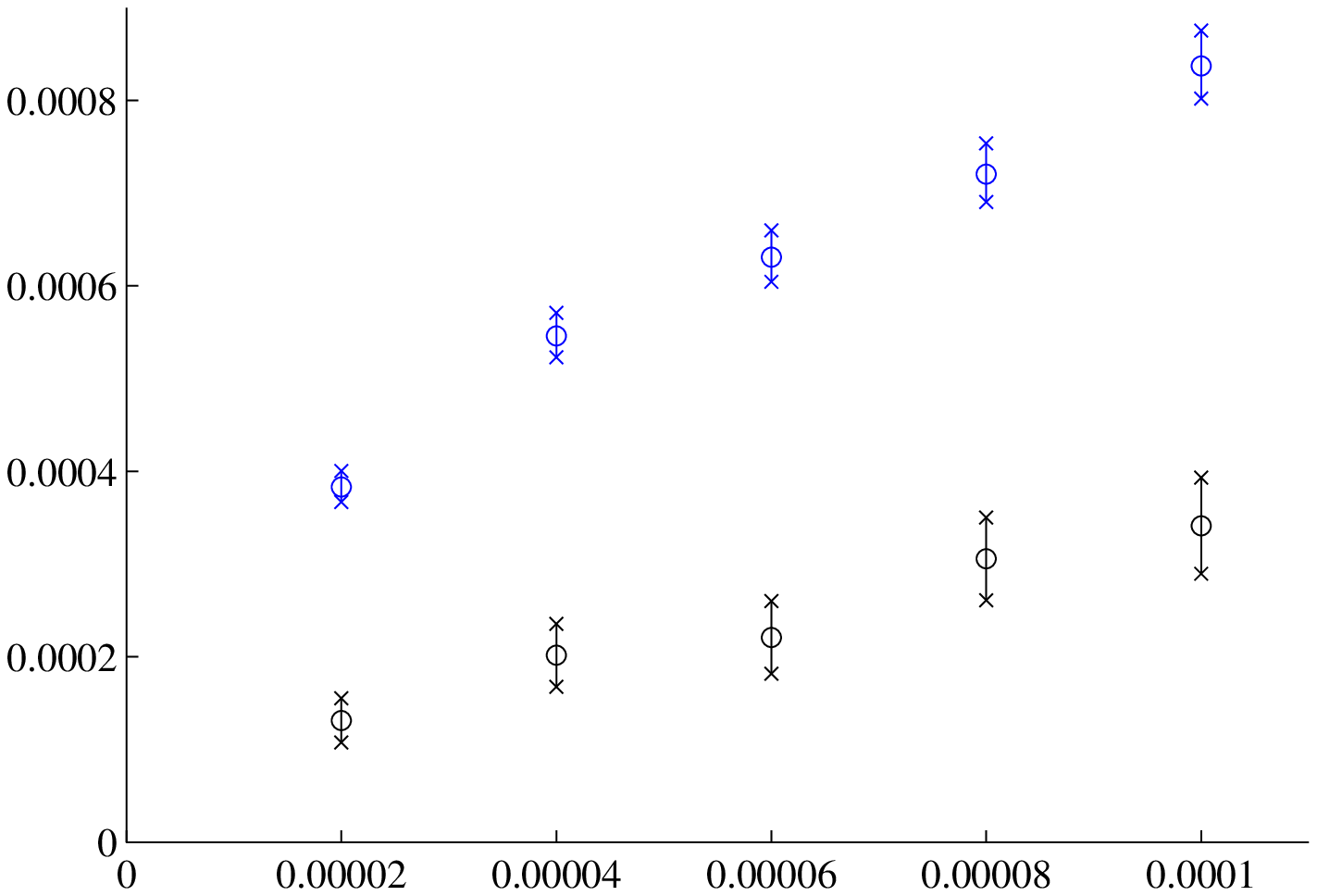}}
\put(4.75,12.8){\large \sf \bfseries A}
\put(1,5.8){\large \sf \bfseries B}
\put(8.5,5.8){\large \sf \bfseries C}
\put(7.55,7){\small $\ee$}
\put(9.5,10.15){\small ${\rm Diff} \left( t^E \right)$}
\put(6.3,11){\small \color{blue} ${\rm Std} \left( t^E \right)$}
\put(3.8,0){\small $\ee$}
\put(5,2.1){\small ${\rm Diff} \left( x_1^E \right)$}
\put(2.7,4.2){\small \color{blue} ${\rm Std} \left( x_1^E \right)$}
\put(11.3,0){\small $\ee$}
\put(12.5,1.65){\small ${\rm Diff} \left( x_3^E \right)$}
\put(10.9,4.6){\small \color{blue} ${\rm Std} \left( x_3^E \right)$}
\end{picture}
\caption{
First passage statistics of an escaping phase for
the relay control example, (\ref{eq:relayControlSystem2}) with
(\ref{eq:ABCDvalues3d})-(\ref{eq:paramValues}).
Here $\bx_0 = \bx_{\Gamma}^S$,
and we have used (\ref{eq:specificdeltas}) for the values of $\delta^-$ and $\delta^+$.
The deterministic passage time and location
are $t_{\dd}^E \approx 0.0668$ and $\bx_{\dd}^E \approx (0.00415,0.2,-0.000642)$.
The data points were computed from $1000$ Monte-Carlo simulations for each value of $\ee$
using the Euler-Maruyama method with a fixed step size, $\Delta t = 0.00001$.
We have included $95\%$ confidence intervals using the same conventions as in Fig.~\ref{fig:manyPeriod}.
The solid curves in panel B are theoretical approximations
obtained by evaluating (\ref{eq:Kn00}).
\label{fig:checkEscape}
}
\end{center}
\end{figure}

Fig.~\ref{fig:checkEscape} shows the result of Monte-Carlo simulations of the relay control example
(\ref{eq:relayControlSystem2}) with (\ref{eq:ABCDvalues3d})-(\ref{eq:paramValues})
for first passage from $\bx = \bx_{\dd}^S$ to $x_2 = \delta^+$
using $\delta^- = -0.1$ and $\delta^+ = 0.2$.
The curves in panel B were obtained by using (\ref{eq:Kn00}).
If smaller values of $\delta^-$ and $\delta^+$ are used,
then the approximation of ${\rm Diff}(x_1^E)$ for small $\ee$ improves
because the approximation is applied over a smaller region.
Since $\delta^-$ and $\delta^+$ are small,
the values in Fig.~\ref{fig:checkEscape} are significantly smaller
than the analogous values of the previous two sections.
In the next section we find that,
in agreement with this observation,
the escaping phase does not have a significant effect on the statistics of the oscillation time, $t_{\rm osc}$,
because the escaping phase corresponds to a relatively short time-frame.
Moreover, the leading-order description of the escaping dynamics derived in this section
does not provide us with a way to accurately approximate the data in panels A and C.

\section{Combining the results}
\label{sec:COMB}
\setcounter{equation}{0}

Here we use the results of Sections \ref{sec:EXCUR}--\ref{sec:ESCAPE}
to construct approximations to ${\rm Diff}(t_{\rm osc})$ and ${\rm Std}(t_{\rm osc})$
for the relay control system (\ref{eq:ABCDvalues3d})-(\ref{eq:relayControlSystem2}).
We employ approximations at various stages of the construction
in order to obtain results that can be interpreted in terms 
of the values of the parameters and geometric features of the system.
The first passage times and locations discussed in sections \ref{sec:EXCUR}--\ref{sec:ESCAPE}
are stochastic quantities that depend significantly on one another.
We were able to ignore this interdependence in these sections
because in each case our attention was restricted to an individual phase.
In this section, however, it is necessary to consider the dependence carefully.

As evident from Fig.~\ref{fig:ppSketch},
on either side of the switching manifold, $\Gamma$ rapidly approaches a slow manifold.
For this reason, each point at which a regular phase ends, denoted $\bx^R$,
is practically independent of the point at which the regular phase starts, denoted $\bx^E$.
Consequently, we ignore the distribution of $\bx^E$ in the computation of the distribution of $\bx^R$.
For systems for which such a simplification is not possible,
one could use the results of the previous three sections
to numerically evaluate the statistics of the oscillation time via an iterative procedure.
In such a procedure one would consecutively apply
the distributions of the various phases of the dynamics in order,
rather than a stochastic simulation of many realizations of the equations.

The final approximations are affected by the values of $\delta^-$ and $\delta^+$.
In particular, the results for escaping are based on a series expansion of the system about $\bx = 0$
that is applied for values of $x_2$ over the range $\delta^- < x_2 < \delta^+$.
For this reason, errors in approximations relating to escaping increase with the magnitude of $\delta^-$ and $\delta^+$,
and therefore we need $\delta^-$ and $\delta^+$ to be small.
However, the results for sliding are singular in the limit $\delta^- \to 0$,
and for this reason the accuracy of approximations relating to sliding decrease as $\delta^-$ approaches zero.
Also, the results for the regular phase assume that initial points on $x_2 = \delta^+$
are sufficiently far from $x_1 = 0$ so that a sample solution from an initial point is highly
unlikely to reach $x_1 = 0$ before undergoing a large excursion with $x_1 > 0$.
Hence the values of $\delta^-$ and $\delta^+$ cannot be too small.
For simplicity we take $\delta^-$ and $\delta^+$ independent of $\ee$.
In view of the above points, and based on using $\ee \le 0.0001$, throughout this paper we have used
\begin{equation}
\delta^- = -0.1 \;, \qquad \delta^+ = 0.2 \;.
\label{eq:specificdeltas}
\end{equation}
The final approximations are not substantially altered by using other values of $\delta^-$ and $\delta^+$
that are the same order of magnitude as the values in (\ref{eq:specificdeltas}).

Next we introduce additional notation, 
for which it is helpful refer to Fig.~\ref{fig:phaseSchem}.
From an initial point $\bx^M$ that lies on the switching manifold and is near $\bx_{\Gamma}^M$,
we let $t^S$ and $\bx^S$ denote the first passage time and location to $x_2 = \delta^-$.
We let
\begin{eqnarray}
{\rm Diff} \left( t^S \big| \bx^M \right) &\equiv& \mathbb{E} \left[ t^S \big| \bx^M \right]
- t_{\dd}^S \left( \bx^M \right) \;,
\label{eq:Difftsl} \\
{\rm Diff} \left( \bx^S \big| \bx^M \right) &\equiv& \mathbb{E} \left[ \bx^S \big| \bx^M \right]
- \bx_{\dd}^S \left( \bx^M \right) \;,
\label{eq:Diffbxsl}
\end{eqnarray}
denote the differences between their means and deterministic values.
Below we evaluate (\ref{eq:Difftsl}) and (\ref{eq:Diffbxsl}) at $\bx_{\Gamma}^M$ by using
(\ref{eq:tSlMeanSeries}) to compute $\mathbb{E} \left[ t^S \big| \bx_{\Gamma}^M \right]$, and
(\ref{eq:bySlMeanSeries}) and (\ref{eq:meanx1sl})
to compute $\mathbb{E} \left[ \bx^S \big| \bx_{\Gamma}^M \right]$.
Also, we evaluate ${\rm Std} \left( t^S \big| \bx_{\Gamma}^M \right)$ with (\ref{eq:vartsl}),
and ${\rm Cov} \left( \bx^S \big| \bx_{\Gamma}^M \right)$ with (\ref{eq:covbysl}) and (\ref{eq:varx1sl}).

Similarly, from an initial point $\bx^E$ that lies on $x_2 = \delta^+$ and near $\bx_{\Gamma}^E$,
we let $t^R$ and $\bx^R$ denote the first passage time and location to the switching manifold.
Below we evaluate
\begin{eqnarray}
{\rm Diff} \left( t^R \big| \bx^E \right) &\equiv& \mathbb{E} \left[ t^R \big| \bx^E \right]
- t_{\dd}^R \left( \bx^E \right) \;,
\label{eq:Difftex} \\
{\rm Diff} \left( \bx^R \big| \bx^E \right) &\equiv& \mathbb{E} \left[ \bx^R \big| \bx^E \right]
- \bx_{\dd}^R \left( \bx^E \right) \;,
\label{eq:Diffbxex}
\end{eqnarray}
at $\bx_{\Gamma}^E$ by using (\ref{eq:meantex3}) and (\ref{eq:meanxjex2}), respectively.
${\rm Std} \left( t^R \big| \bx_{\Gamma}^E \right)$ and ${\rm Cov} \left( \bx^R \big| \bx_{\Gamma}^E \right)$
are given by (\ref{eq:vartex}) and (\ref{eq:covxex}), respectively.

Calculations relating the escaping phase do not enter into our final approximations
because escaping phases occur over significantly shorter time-frames than sliding and regular phases,
as discussed at the end of \S\ref{sec:ESCAPE}.

\subsection{An approximation to ${\rm Diff}(t_{\rm osc})$}
\label{sub:DIFF}

From a sample solution to (\ref{eq:ABCDvalues3d})-(\ref{eq:relayControlSystem2})
computed over a length of time that is substantially greater than the period of $\Gamma$,
we can identify first passage locations,
$\bx^M$, $\bx^S$, $\bx^E$ and $\bx^R$,
and first passage times,
$t^S$, $t^E$ and $t^R$
corresponding to the beginning and end of sliding, escaping and regular phases.
Since (\ref{eq:ABCDvalues3d})-(\ref{eq:relayControlSystem2}) exhibits a simple symmetry about $x_1 = 0$,
the distributions of $\bx^M$ and $\bx^R$ are symmetric.
Also, as discussed above, each $\bx^R$ is practically independent of the previous point $\bx^E$.
Consequently it is suitable to use the approximation
\begin{equation}
{\rm Diff} \left( \bx^M \right) \approx -{\rm Diff} \left( \bx^R \big| \bx_{\Gamma}^E \right) \;.
\label{eq:Diffbxpr2}
\end{equation}

\begin{figure}[b!]
\begin{center}
\setlength{\unitlength}{1cm}
\begin{picture}(15,6)
\put(0,0){\includegraphics[height=6cm]{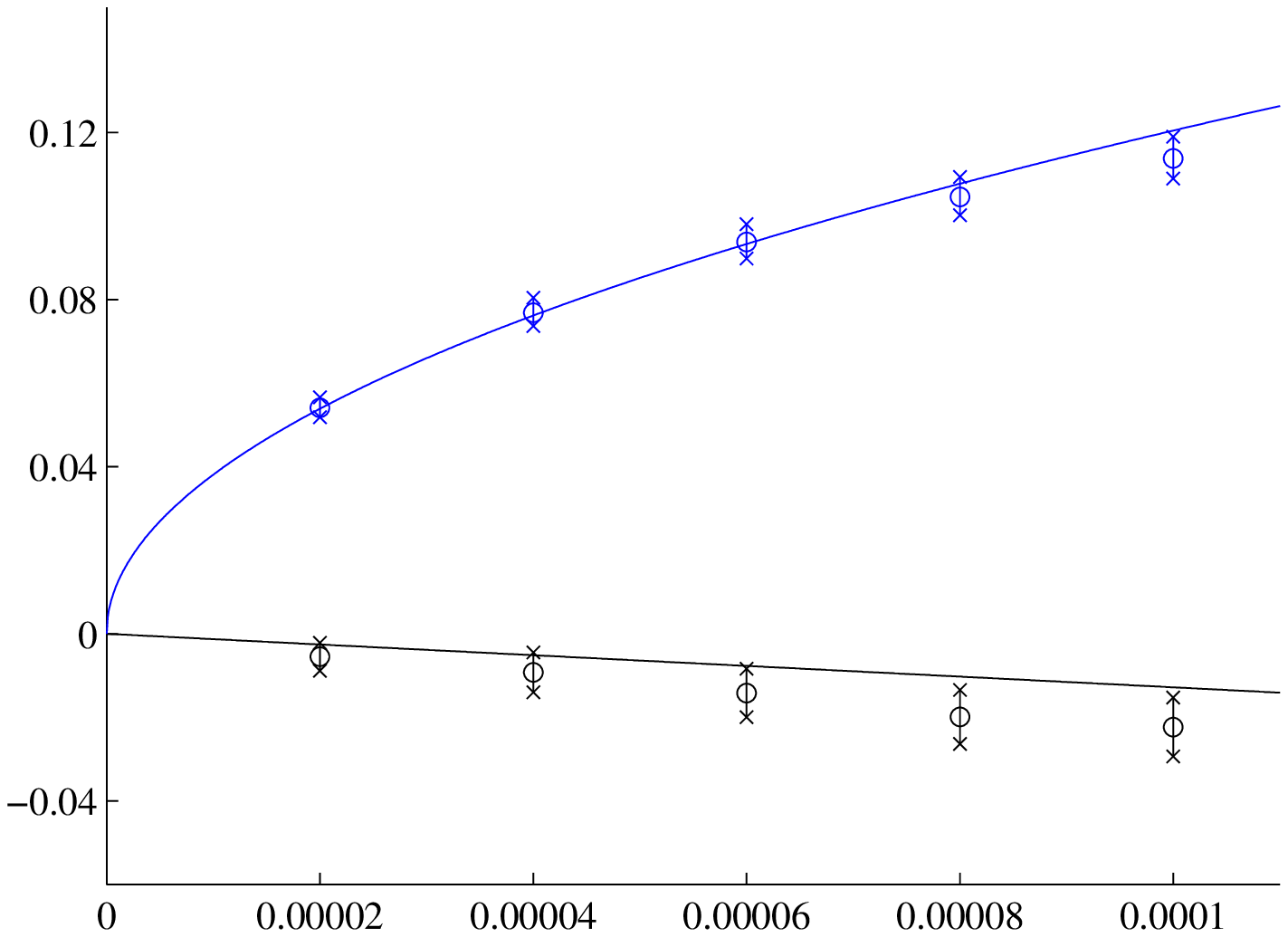}}
\put(7.5,0){\includegraphics[height=6cm]{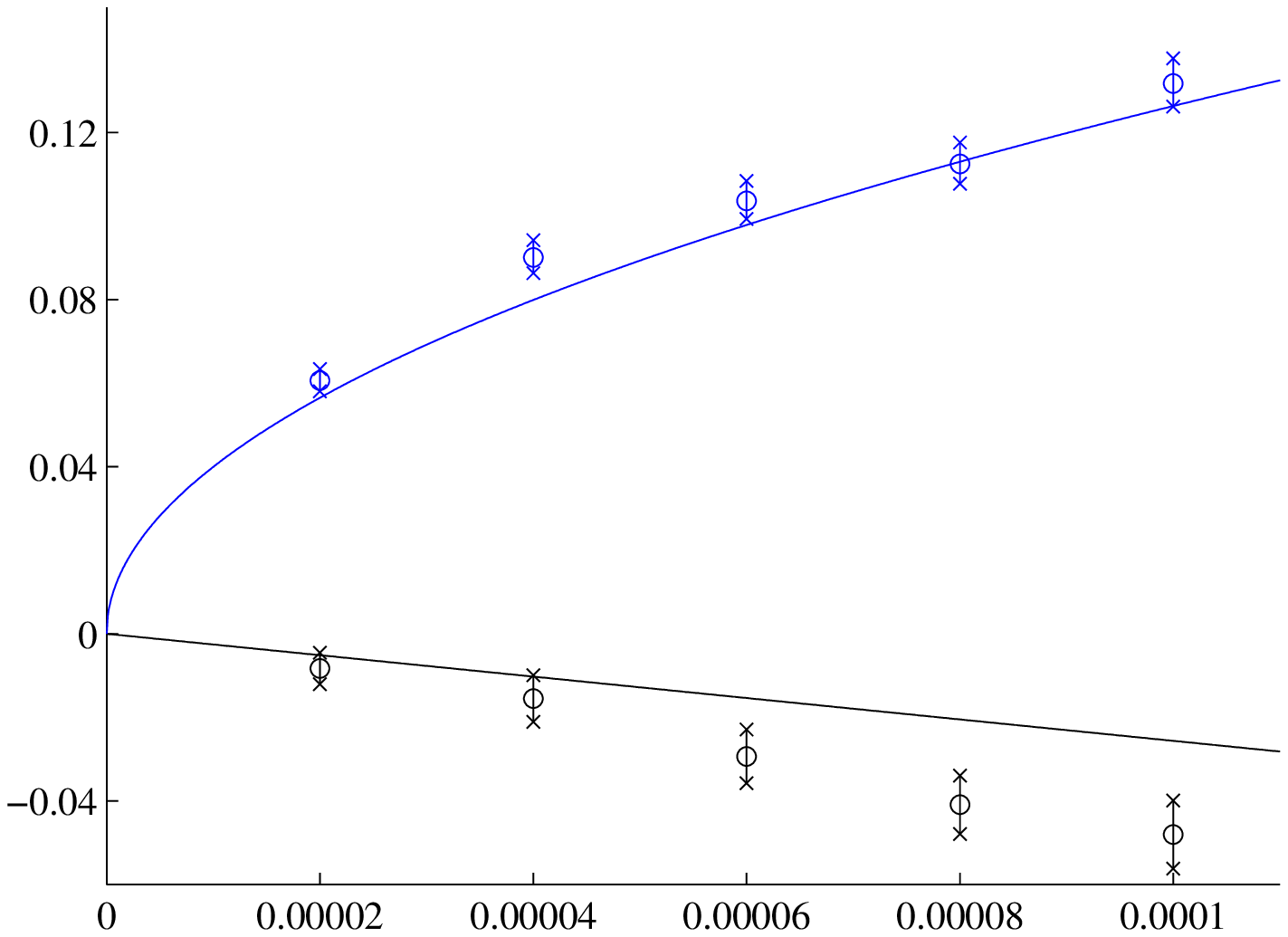}}
\put(1,5.8){\large \sf \bfseries A}
\put(8.5,5.8){\large \sf \bfseries B}
\put(3.8,0){\small $\ee$}
\put(5,2.3){\small ${\rm Diff} \left( t_{\frac{1}{2} {\rm osc}} \right)$}
\put(2.7,4.9){\small \color{blue} ${\rm Std} \left( t_{\frac{1}{2} {\rm osc}} \right)$}
\put(11.3,0){\small $\ee$}
\put(12.5,2){\small ${\rm Diff} \left( t_{\rm osc} \right)$}
\put(10.9,5.15){\small \color{blue} ${\rm Std} \left( t_{\rm osc} \right)$}
\end{picture}
\caption{
A comparison of Monte-Carlo simulations with the theoretical approximations derived in the text
for oscillation times of the system (\ref{eq:paramValues})-(\ref{eq:relayControlSystem2}).
As in Fig.~\ref{fig:manyPeriod}, the data points were computed by solving the system
using the Euler-Maruyama method with $\Delta t = 0.00001$.
For panel A, 500 oscillations were computed yielding 1000 values of $t_{\frac{1}{2} {\rm osc}}$
-- the time taken to return to the switching manifold after one large excursion.
Twice as many oscillations were computed for panel B.
The circles, bars and crosses indicate mean values and 95\% confidence intervals, as in Fig.~\ref{fig:manyPeriod}.
The solid straight lines are the theoretical approximations (\ref{eq:Diffthalfosc2}) and (\ref{eq:Difftosc2}).
The solid curves are the theoretical approximations (\ref{eq:Varthalfosc2}) and (\ref{eq:Vartosc2}).
\label{fig:checkOsc}
}
\end{center}
\end{figure}

We can compute ${\rm Diff} \left( t^S \right)$
by evaluating the following expression that is derived in Appendix \ref{sec:TSL}
via a Taylor series expansion
\begin{equation}
{\rm Diff} \left( t^S \right) = {\rm Diff} \left( t^S \big| \bx_{\Gamma}^M \right) 
+ \rD_{\bx} t_{\dd}^S \left( \bx_{\Gamma}^M \right)^{\sf T} {\rm Diff} \left( \bx^M \right)
+ \sum_{i=1}^N \sum_{j=1}^N \rD_{\bx}^2 t_{\dd}^S \left( \bx_{\Gamma}^M \right)_{i,j}
{\rm Cov} \left( x_{\Gamma}^M \right)_{i,j} + O \left( \ee^{\frac{3}{2}} \right) \;.
\label{eq:Difftsl2}
\end{equation}
Note, $t_{\dd}^S$ is a function of the point $\bx^M$ at which the sliding phase begins.
In (\ref{eq:Difftsl2}), $t_{\dd}^S$ and its derivatives are evaluated at $\bx^M = \bx_{\Gamma}^M$.
Each term in (\ref{eq:Difftsl2}) is $O(\ee)$,
but, for our example, components of the vector ${\rm Diff} \left( \bx^M \right)$ are of much larger magnitude than
elements of the matrix ${\rm Cov} \left( x_{\Gamma}^M \right)$.
For this reason we use the approximation
\begin{equation}
{\rm Diff} \left( t^S \right) \approx {\rm Diff} \left( t^S \big| \bx_{\Gamma}^M \right) 
+ \rD_{\bx} t_{\dd}^S \left( \bx_{\Gamma}^M \right)^{\sf T} {\rm Diff} \left( \bx^M \right) \;,
\label{eq:Difftsl3}
\end{equation}
which is evaluated using (\ref{eq:Diffbxpr2}).
Similarly we use
\begin{align}
{\rm Diff} \left( t^E \right) &\approx {\rm Diff} \left( t^E \big| \bx_{\Gamma}^S \right) 
+ \rD_{\bx} t_{\dd}^E \left( \bx_{\Gamma}^S \right)^{\sf T} {\rm Diff} \left( \bx^S \right) \;, \\
{\rm Diff} \left( t^R \right) &\approx {\rm Diff} \left( t^R \big| \bx_{\Gamma}^E \right)
+ \rD_{\bx} t_{\dd}^R \left( \bx_{\Gamma}^E \right)^{\sf T} {\rm Diff} \left( \bx^E \right) \;, \\
{\rm Diff} \left( \bx^S \right) &\approx {\rm Diff} \left( \bx^S \big| \bx_{\Gamma}^M \right)
+ \rD_{\bx} \bx_{\dd}^S \left( \bx_{\Gamma}^M \right) {\rm Diff} \left( \bx^M \right) \;, \\
{\rm Diff} \left( \bx^E \right) &\approx {\rm Diff} \left( \bx^E \big| \bx_{\Gamma}^S \right)
+ \rD_{\bx} \bx_{\dd}^E \left( \bx_{\Gamma}^S \right) {\rm Diff} \left( \bx^S \right) \;.
\label{eq:DiffMany}
\end{align}
The difference for the time of half an oscillation is given simply by\removableFootnote{
Despite weak dependencies between the $t^S$, $t^E$ and $t^R$, this expression is exact
because it concerns mean values.
Compare equation (\ref{eq:Varthalfosc3}).
}
\begin{equation}
{\rm Diff} \left( t_{\frac{1}{2} {\rm osc}} \right) = {\rm Diff} \left( t^S \right)
+ {\rm Diff} \left( t^E \right) + {\rm Diff} \left( t^R \right) \;.
\label{eq:Diffthalfosc}
\end{equation}
Substituting (\ref{eq:Difftsl3})-(\ref{eq:DiffMany}) into (\ref{eq:Diffthalfosc}) and expanding brackets
yields an approximation for ${\rm Diff}(t_{\frac{1}{2} {\rm osc}})$
that is a sum of nine terms\removableFootnote{
Specifically
\begin{eqnarray}
{\rm Diff} \left( t_{\frac{1}{2} {\rm osc}} \right)
&=& {\rm Diff} \left( t^S \big| \bx_{\Gamma}^M \right)
+ \rD_{\bx} t_{\dd}^S \left( \bx_{\Gamma}^M \right)^{\sf T} {\rm Diff} \left( \bx^M \right)
+ {\rm Diff} \left( t^E \big| \bx_{\Gamma}^S \right) \nonumber \\
&&+~\rD_{\bx} t_{\dd}^E \left( \bx_{\Gamma}^S \right)^{\sf T} {\rm Diff} \left( \bx^S \big| \bx_{\Gamma}^M \right)
+ \rD_{\bx} t_{\dd}^E \left( \bx_{\Gamma}^S \right)^{\sf T} \rD_{\bx} \bx_{\dd}^S \left( \bx_{\Gamma}^M \right)
{\rm Diff} \left( \bx^M \right) \nonumber \\
&&+~{\rm Diff} \left( t^R \big| \bx_{\Gamma}^E \right)
+ \rD_{\bx} t_{\dd}^R \left( \bx_{\Gamma}^E \right)^{\sf T} {\rm Diff} \left( \bx^E \big| \bx_{\Gamma}^S \right) \nonumber \\
&&+~\rD_{\bx} t_{\dd}^R \left( \bx_{\Gamma}^E \right)^{\sf T} \rD_{\bx} \bx_{\dd}^E \left( \bx_{\Gamma}^S \right)
{\rm Diff} \left( \bx^S \big| \bx_{\Gamma}^M \right) \nonumber \\
&&+~\rD_{\bx} t_{\dd}^R \left( \bx_{\Gamma}^E \right)^{\sf T} \rD_{\bx} \bx_{\dd}^E \left( \bx_{\Gamma}^S \right)
\rD_{\bx} \bx_{\dd}^S \left( \bx_{\Gamma}^M \right) {\rm Diff} \left( \bx^M \right) \;.
\end{eqnarray}
}.
Monte-Carlo simulations reveal that for our example
three of these terms have significantly larger values than the remaining six terms,
and for simplicity we approximate ${\rm Diff}(t_{\frac{1}{2} {\rm osc}})$
using the three largest terms:
\begin{equation}
{\rm Diff} \left( t_{\frac{1}{2} {\rm osc}} \right)
\approx {\rm Diff} \left( t^R \big| \bx_{\Gamma}^E \right)
+ \rD_{\bx} t_{\dd}^S \left( \bx_{\Gamma}^M \right)^{\sf T} {\rm Diff} \left( \bx^M \right)
+ {\rm Diff} \left( t^S \big| \bx_{\Gamma}^M \right) \;.
\label{eq:Diffthalfosc2}
\end{equation}
The first term in (\ref{eq:Diffthalfosc2}) represents
the additional time that regular phases take, on average, due to the presence of noise.
This term is negative-valued and large because by (\ref{eq:meantex3})
its magnitude is proportional to $\omega^4$, and we have used $\omega = 5$.
The second term in (\ref{eq:Diffthalfosc2}) represents
the average additional time that sliding phases take
due to noise causing the points $\bx^M$, at which sliding phases start,
to be deviated from the deterministic value $\bx_{\Gamma}^M$ in a particular direction, on average.
Indeed, as evident in Fig.~\ref{fig:checkExcur}-C,
small noise induces a large positive shift in the $x_3$-component of the average value of $\bx^M$.
Finally the third term of (\ref{eq:Diffthalfosc2}) represents
the additional time that sliding phases take, on average, due to the noise.
In view of (\ref{eq:tSlMeanSeries}), (\ref{eq:yBar1}) and (\ref{eq:meantsl}),
this term is proportional to $\Lambda$ (\ref{eq:Lambda}).
The third term of (\ref{eq:Diffthalfosc2}) is large,
but not as large as the first term because $\Lambda$ is proportional to $\omega^2$.
The approximation (\ref{eq:Diffthalfosc2}) is compared with Monte-Carlo simulations in Fig.~\ref{fig:checkOsc}-A.
Also
\begin{equation}
{\rm Diff}(t_{\rm osc}) = 2 {\rm Diff}(t_{\frac{1}{2} {\rm osc}}) \;,
\label{eq:Difftosc2}
\end{equation}
is used for Fig.~\ref{fig:checkOsc}-B.

\subsection{An approximation to ${\rm Std}(t_{\rm osc})$}

Taylor expanding $\bx^M$ about its deterministic value $\bx_{\Gamma}^M$,
(mentioned in \S\ref{sub:DIFF}), also leads to the formula
\begin{equation}
{\rm Var} \left( t^S \right) = {\rm Var} \left( t^S \big| \bx_{\Gamma}^M \right)
+ \rD_{\bx} t_{\dd}^S \left( \bx_{\Gamma}^M \right)^{\sf T}
{\rm Cov} \left( \bx^M \right)
\rD_{\bx} t_{\dd}^S \left( \bx_{\Gamma}^M \right) + O \left( \ee^{\frac{3}{2}} \right) \;,
\label{eq:combvartsl}
\end{equation}
which expresses ${\rm Var} \left( t^S \right)$ in terms of elements
that we can evaluate using equations derived in earlier sections.
A derivation of (\ref{eq:combvartsl}) is given in Appendix \ref{sec:TSL}.
Via similar calculations we obtain
\begin{eqnarray}
{\rm Var} \left( t^E \right) &=& {\rm Var} \left( t^E \big| \bx_{\Gamma}^S \right)
+ \rD_{\bx} t_{\dd}^E  \left( \bx_{\Gamma}^S \right)^{\sf T}
{\rm Cov} \left( \bx^S \right)
\rD_{\bx} t_{\dd}^E  \left( \bx_{\Gamma}^S \right) + O \left( \ee^{\frac{3}{2}} \right) \;, \\
{\rm Var} \left( t^R \right) &=& {\rm Var} \left( t^R \big| \bx_{\Gamma}^E \right)
+ \rD_{\bx} t_{\dd}^R  \left( \bx_{\Gamma}^E \right)^{\sf T}
{\rm Cov} \left( \bx^E \right)
\rD_{\bx} t_{\dd}^R  \left( \bx_{\Gamma}^E \right) + O \left( \ee^{\frac{3}{2}} \right) \;, \\
{\rm Cov} \left( \bx^S \right) &=& {\rm Cov} \left( \bx^S \big| \bx_{\Gamma}^M \right)
+ \rD_{\bx} \bx_{\dd}^S  \left( \bx_{\Gamma}^M \right)
{\rm Cov} \left( \bx^M \right)
\rD_{\bx} \bx_{\dd}^S  \left( \bx_{\Gamma}^M \right)^{\sf T} + O \left( \ee^{\frac{3}{2}} \right) \;, \\
{\rm Cov} \left( \bx^E \right) &=& {\rm Cov} \left( \bx^E \big| \bx_{\Gamma}^S \right)
+ \rD_{\bx} \bx_{\dd}^E  \left( \bx_{\Gamma}^S \right)
{\rm Cov} \left( \bx^S \right)
\rD_{\bx} \bx_{\dd}^E  \left( \bx_{\Gamma}^S \right)^{\sf T} + O \left( \ee^{\frac{3}{2}} \right) \;.
\label{eq:VarMany}
\end{eqnarray}
We use the approximation
\begin{equation}
{\rm Var} \left( t_{\frac{1}{2} {\rm osc}} \right) \approx {\rm Var} \left( t^S \right)
+ {\rm Var} \left( t^E \right) + {\rm Var} \left( t^R \right) \;,
\label{eq:Varthalfosc3}
\end{equation}
because, for our example, the value of each $t^E$ is practically independent
of the preceding value of $t^S$, and the value of each $t^R$ is practically
independent of the preceding value of $t^E$.
This is due to strong attraction to $x_1=0$ for the duration of the sliding phase,
which is inherent in stochastically perturbed sliding motion
and causes the $x_1$-value of $\bx^S$ (which is the primary influence on the value of $t^E$)
to have a negligible correlation to $t^S$.
In an analogous fashion to the calculations in \S\ref{sub:DIFF},
by substituting (\ref{eq:combvartsl})-(\ref{eq:VarMany}) into (\ref{eq:Varthalfosc3})
and expanding brackets, we produce an approximation to ${\rm Var}(t_{\frac{1}{2} {\rm osc}})$
that is a sum of nine terms\removableFootnote{
Specifically,
\begin{eqnarray}
{\rm Var} \left( t_{\frac{1}{2} {\rm osc}} \right)
&=& {\rm Var} \left( t^S \big| \bx_{\Gamma}^M \right)
+ \rD_{\bx} t_{\dd}^S \left( \bx_{\Gamma}^M \right)^{\sf T}
{\rm Cov} \left( \bx^M \right)
\rD_{\bx} t_{\dd}^S \left( \bx_{\Gamma}^M \right)
+ {\rm Var} \left( t^E \big| \bx_{\Gamma}^S \right) \nonumber \\
&&+~\rD_{\bx} t_{\dd}^E \left( \bx_{\Gamma}^S \right)^{\sf T}
{\rm Cov} \left( \bx^S \big| \bx_{\Gamma}^M \right)
\rD_{\bx} t_{\dd}^E \left( \bx_{\Gamma}^S \right) \nonumber \\
&&+~\rD_{\bx} t_{\dd}^E \left( \bx_{\Gamma}^S \right)^{\sf T}
\rD_{\bx} \bx_{\dd}^S \left( \bx_{\Gamma}^M \right)
{\rm Cov} \left( \bx^M \right)
\rD_{\bx} \bx_{\dd}^S \left( \bx_{\Gamma}^M \right)^{\sf T}
\rD_{\bx} t_{\dd}^E \left( \bx_{\Gamma}^S \right) \nonumber \\
&&+~{\rm Var} \left( t^R \big| \bx_{\Gamma}^E \right)
+ \rD_{\bx} t_{\dd}^R \left( \bx_{\Gamma}^E \right)^{\sf T}
{\rm Cov}\left( \bx^E \big| \bx_{\Gamma}^S \right)
\rD_{\bx} t_{\dd}^R \left( \bx_{\Gamma}^E \right) \nonumber \\
&&+~\rD_{\bx} t_{\dd}^R \left( \bx_{\Gamma}^E \right)^{\sf T}
\rD_{\bx} \bx_{\dd}^E \left( \bx_{\Gamma}^S \right)
{\rm Cov} \left( \bx^S \big| \bx_{\Gamma}^M \right)
\rD_{\bx} \bx_{\dd}^E \left( \bx_{\Gamma}^S \right)^{\sf T}
\rD_{\bx} t_{\dd}^R \left( \bx_{\Gamma}^E \right) \nonumber \\
&&+~\rD_{\bx} t_{\dd}^R \left( \bx_{\Gamma}^E \right)^{\sf T}
\rD_{\bx} \bx_{\dd}^E \left( \bx_{\Gamma}^S \right)
\rD_{\bx} \bx_{\dd}^S \left( \bx_{\Gamma}^M \right)
{\rm Cov} \left( \bx^M \right) \nonumber \\
&& \qquad \rD_{\bx} \bx_{\dd}^S \left( \bx_{\Gamma}^M \right)^{\sf T}
\rD_{\bx} \bx_{\dd}^E \left( \bx_{\Gamma}^S \right)^{\sf T}
\rD_{\bx} t_{\dd}^R \left( \bx_{\Gamma}^E \right) \;.
\end{eqnarray}
}.
Monte-Carlo simulations reveal that three of these terms dominate.
By dropping the other six terms we generate the approximation
\begin{align}
&{\rm Var} \left( t_{\frac{1}{2} {\rm osc}} \right)
\approx \rD_{\bx} t_{\dd}^S \left( \bx_{\Gamma}^M \right)^{\sf T}
{\rm Cov} \left( \bx^M \right)
\rD_{\bx} t_{\dd}^S \left( \bx_{\Gamma}^M \right)
+ {\rm Var} \left( t^R \big| \bx_{\Gamma}^E \right) \nonumber \\
&+~\rD_{\bx} t_{\dd}^R \left( \bx_{\Gamma}^E \right)^{\sf T}
\rD_{\bx} \bx_{\dd}^E \left( \bx_{\Gamma}^S \right)
\rD_{\bx} \bx_{\dd}^S \left( \bx_{\Gamma}^M \right)
{\rm Cov} \left( \bx^M \right) \rD_{\bx} \bx_{\dd}^S \left( \bx_{\Gamma}^M \right)^{\sf T}
\rD_{\bx} \bx_{\dd}^E \left( \bx_{\Gamma}^S \right)^{\sf T}
\rD_{\bx} t_{\dd}^R \left( \bx_{\Gamma}^E \right) \;.
\label{eq:Varthalfosc2}
\end{align}
The three terms in (\ref{eq:Varthalfosc2}) can be interpreted geometrically.
Noise creates variability in the values of $\bx^M$.
The variance that this induces in the values of the sliding times $t^S$
is represented by the first term of (\ref{eq:Varthalfosc2}).
Variability in $\bx^M$ is also responsible for variance in the time of regular phases;
this is represented by the third term of (\ref{eq:Varthalfosc2}).
The second term of (\ref{eq:Varthalfosc2})
simply represents the variance in the time of the regular phases
given that regular phases start at the deterministic location $\bx_{\Gamma}^E$.
The approximation (\ref{eq:Varthalfosc2}) is compared with Monte-Carlo simulations in Fig.~\ref{fig:checkOsc}-A.

Lastly we determine ${\rm Var} \left( t_{\rm osc} \right)$ from ${\rm Var}(t_{\frac{1}{2} {\rm osc}})$.
Each oscillation time $t_{\rm osc}$ is the sum of two consecutive half oscillation times,
call them $t_{\frac{1}{2} {\rm osc},1}$ and $t_{\frac{1}{2} {\rm osc},2}$.
For our example, the value of each $t_{\frac{1}{2} {\rm osc},2}$
depends on heavily on the value of $t_{\frac{1}{2} {\rm osc},1}$.
This is because if $t_{\frac{1}{2} {\rm osc},1}$ is, say,
less than its deterministic value $\frac{t_{{\rm osc},\Gamma}}{2}$,
then the point at which the first half oscillation ends, $\bx^R$,
is likely to be skewed in a particular direction from $\bx_{\Gamma}^R$.
The second half oscillation begins at this end point
which affects the value of $t_{\frac{1}{2} {\rm osc},2}$.

To treat this difficulty, we define
\begin{equation}
\varrho = \frac{d \mathbb{E} \left[ t_{\frac{1}{2} {\rm osc},2} \big| t_{\frac{1}{2} {\rm osc},1} \right]}
{d t_{\frac{1}{2} {\rm osc},1}}
\Bigg|_{t_{\frac{1}{2} {\rm osc},1} \,=\, t_{\frac{1}{2} {\rm osc},\Gamma}} \;,
\label{eq:varrho}
\end{equation}
which measures the rate at which the mean value of $t_{\frac{1}{2} {\rm osc},2}$,
given $t_{\frac{1}{2} {\rm osc},1}$,
changes with $t_{\frac{1}{2} {\rm osc},1}$.
If $t_{\frac{1}{2} {\rm osc},1}$ and $t_{\frac{1}{2} {\rm osc},2}$ were independent
then we would have $\varrho = 0$.
Via straight-forward calculations based upon conditioning over the value of $\bx^R$\removableFootnote{
Let us derive this formula in a slightly more general context.
Consider a stochastic trajectory that travels first to a surface, $\Sigma_1$,
then on to another surface, $\Sigma_2$.
Let $t_1$ and $t_2$ be the respective passage times
and $\bx_{\rm mid}$ the intersection point on $\Sigma_1$.
We assume $t_1$ and $\bx_{\rm mid}$
are correlated normally distributed random variables with known statistics:
\begin{eqnarray}
t_1 &\sim& N(t_{{\dd},1},\ee t_{1,{\rm var}}) \;, \\
\bx_{\rm mid} &\sim& N(\bx_{\rm mid,det},\ee \bx_{\rm mid,cov}) \;.
\end{eqnarray}
We also assume that, for any given $\bx_{\rm mid}$,
$t_2$ is a normally distributed random variable with known statistics, and we write
\begin{equation}
t_2 ~\big|~ \bx_{\rm mid} \sim N(t_{{\dd},2}(\bx_{\rm mid}), \ee t_{2,{\rm var}}) \;,
\end{equation}
where in this exposition we ignore the dependence of
$t_{2,{\rm var}}$ on $\bx_{\rm mid}$ because this provides a final contribution
that is of an order higher than $\ee$.
Our goal is to derive ${\rm Var}(t_1 + t_2)$.
This non-trivial because $t_1$ and $\bx_{\rm mid}$ are correlated,
and therefore $t_1$ and $t_2$ are correlated.
Let $t_{\rm osc} = t_1 + t_2$.
Then
\begin{eqnarray}
p(t_{\rm osc}) &=& \int_{-\infty}^\infty
p(t_2 = t_{\rm osc} - t_1 | t_1) p(t_1) \,dt_1 \nonumber \\
&=& \int_{-\infty}^\infty
\frac{1}{\sqrt{2 \pi \ee t_{2,{\rm var}}}}
{\rm e}^{\frac{-(t_{\rm osc} - t_1 - \mathbb{E}[t_2|t_1])^2}{2 \ee t_{2,{\rm var}}}}
\frac{1}{\sqrt{2 \pi \ee t_{1,{\rm var}}}}
{\rm e}^{\frac{-(t_1 - t_{{\dd},1})^2}{2 \ee t_{1,{\rm var}}}} \,dt_1 \;.
\label{eq:derivptosc}
\end{eqnarray}
Expanding the function $\mathbb{E}[t_2|t_1]$ as a Taylor series
centred at $t_1 = t_{{\dd},1}$ produces
\begin{equation}
\mathbb{E}[t_2|t_1] = t_{{\dd},2} + \varrho (t_1 - t_{{\dd},1})
+ O(|t_1 - t_{{\dd},1}|^2) \;,
\end{equation}
where we let
\begin{equation}
\varrho = \frac{d \mathbb{E}[t_2|t_{{\dd},1}]}{d t_1} \;.
\end{equation}
Consequently,
$t_{\rm osc} - t_1 - \mathbb{E}[t_2|t_1]
= t_{\rm osc} - t_{{\rm osc},\Gamma} - (1+\varrho)(t_1 - t_{{\dd},1})$,
and by substituting this into (\ref{eq:derivptosc}) we arrive at
\begin{equation}
p(t_{\rm osc}) = \frac{1}{2 \pi \ee \sqrt{t_{1,{\rm var}} t_{2,{\rm var}}}}
\int_{-\infty}^\infty
{\rm e}^{\frac{-1}{2 \ee t_{1,{\rm var}} t_{2,{\rm var}}}
\left( t_{1,{\rm var}} \left( t_{\rm osc} - t_{{\rm osc},\Gamma} - (1+\varrho)(t_1 - t_{{\dd},1}) \right)^2
+ t_{2,{\rm var}} (t_1 - t_{{\dd},1})^2 \right)} \,dt_1 \;.
\end{equation}
By completing the square in the exponent,
we can simplify this to
\begin{equation}
p(t_{\rm osc}) = \frac{1}{\sqrt{2 \pi \ee \left( t_{2,{\rm var}} + (1+\varrho)^2 t_{1,{\rm var}} \right)}}
{\rm e}^{\frac{-(t_{\rm osc} - t_{{\rm osc},\Gamma})^2}
{2 \ee \left( t_{2,{\rm var}} + (1+\varrho)^2 t_{1,{\rm var}} \right)}} \;.
\end{equation}
Therefore,
\begin{equation}
{\rm Var}(t_1 + t_2) = \ee \left( t_{2,{\rm var}} + (1+\varrho)^2 t_{1,{\rm var}} \right) \;.
\end{equation}
},
it can be shown that
\begin{equation}
{\rm Var} \left( t_{\rm osc} \right) = {\rm Var} \left( t_{\frac{1}{2} {\rm osc},2} \right)
+ (1+\varrho)^2 {\rm Var} \left( t_{\frac{1}{2} {\rm osc},1} \right) + O \left( \ee^{\frac{3}{2}} \right) \;.
\end{equation}
Due to the symmetry of the relay control system,
${\rm Var}(t_{\frac{1}{2} {\rm osc},1}) = {\rm Var}(t_{\frac{1}{2} {\rm osc},2})$,
and therefore
\begin{equation}
{\rm Var} \left( t_{\rm osc} \right) =
\left( 1 + (1+\varrho)^2 \right) {\rm Var} \left( t_{\frac{1}{2} {\rm osc}} \right)
+ O \left( \ee^{\frac{3}{2}} \right) \;.
\label{eq:Vartosc2}
\end{equation}
For our example, $\varrho \approx -0.68$.
We obtained this value numerically by first computing the adjusted mean value of $\bx^R$
given that $t_{\frac{1}{2} {\rm osc},1} = \frac{t_{{\rm osc},\Gamma}}{2} + \Delta t$
(using $\Delta t = 0.0001$),
then numerically solving the system with $\ee = 0$ for half an oscillation from $\bx^R$
in order to obtain the expected value of $t_{\frac{1}{2} {\rm osc},2}$,
and lastly using a first order finite difference approximation to evaluate (\ref{eq:varrho}).
The approximation (\ref{eq:Vartosc2}) is compared with Monte-Carlo simulations in Fig.~\ref{fig:checkOsc}-B.

\section{Conclusions}
\label{sec:CONC}
\setcounter{equation}{0}

In this paper we have quantitatively analyzed the effect of noise on periodic orbits
of Filippov systems that involve segments of sliding motion.
Our results apply to the general $N$-dimensional stochastic differential equation, (\ref{eq:sde}),
which is formed by adding white Gaussian noise of amplitude $\sqrt{\ee}$
to a Filippov system with a single switching manifold.
We assume that in the absence of noise, i.e.~with $\ee = 0$, (\ref{eq:sde})
has an attracting periodic orbit $\Gamma$ of period $t_{{\rm osc},\Gamma}$.
For small $\ee > 0$, sample solutions to (\ref{eq:sde}) are likely to follow paths near $\Gamma$.
From such solutions we can identify oscillation times, $t_{\rm osc} \approx t_{{\rm osc},\Gamma}$,
defined by measuring the time taken between appropriate returns to the switching manifold ($x_1 = 0$).
In order to determine the statistics of $t_{\rm osc}$ for small $\ee > 0$,
we split the stochastic dynamics into three phases: regular, sliding and escaping, see Fig.~\ref{fig:phaseSchem},
and analyzed each phase separately.

\subsubsection*{Regular dynamics}

From an initial point $\bx_0$ in the right half-space,
we let $t^R$ and $\bx^R$ denote the time and location for first passage to $x_1 = 0$.
To derive the mean values of these quantities to $O(\ee)$,
we searched for an asymptotic solution to the Fokker-Planck equation of (\ref{eq:sde})
with an absorbing boundary condition at $x_1 = 0$, (\ref{eq:FPEex})-(\ref{eq:FPEremainingBCs}),
by introducing a boundary layer near $x_1 = 0$
and expanding about the deterministic passage time and location (\ref{eq:regularScaling}).
With the solution expanded in the form (\ref{eq:pell}),
only the first term of the local PDF $\mathcal{P}^{(0)}$
is required to obtain $\mathbb{E} \left[ t^R \big| \bx_0 \right]$ to $O(\ee)$.
To determine $\mathbb{E} \left[ \bx^R \big| \bx_0 \right]$ to $O(\ee)$,
we also require the second term, $\mathcal{P}^{(1)}$.
We computed $\mathcal{P}^{(1)}$ by numerically evaluating integrals, see Appendix \ref{sec:CALEX}.
Standard calculations based on a sample path methodology are
sufficient to determine ${\rm Var} \left( t^R \big| \bx_0 \right)$ and
${\rm Cov} \left( \bx^R \big| \bx_0 \right)$ to $O \left( \sqrt{\ee} \right)$.

\subsubsection*{Sliding dynamics}

In \S\ref{sec:SLIDE} we analyzed stochastically perturbed sliding motion.
We let $t^S$ and $\bx^S$ denote the time and location for
the first passage of (\ref{eq:sde}) to $x_2 = \delta^-$ from an initial point $\bx_0$
that lies on the switching manifold.
We assumed that the deterministic solution from $\bx_0$ to $x_2 = \delta^-$
is contained entirely within the interior of a stable sliding region.
Stochastic dynamics of (\ref{eq:sde}) in the $x_1$-direction, i.e.~orthogonal to the switching manifold,
occurs on an $O(\ee)$ time-scale and for this reason it is suitable to employ stochastic averaging
to analyze the overall dynamics from $\bx_0$ to $x_2 = \delta^-$.
To leading order, the averaged solution is identical to Filippov's solution of the deterministic equations.
We estimated ${\rm Diff} \left( t^S \big| \bx_0 \right)$ and
${\rm Diff} \left( \bx^S \big| \bx_0 \right)$ to $O(\ee)$
by including terms of the next order in the averaging calculation.
The key quantity affecting the magnitude of these differences is $\Lambda$ (\ref{eq:Lambda})
which denotes the $O(\ee)$ component of the average drift in directions parallel to the switching manifold.

We obtained the leading order terms of ${\rm Var} \left( t^S \big| \bx_0 \right)$ 
and ${\rm Cov} \left( \bx^S \big| \bx_0 \right)$ 
through the use of a linear diffusion approximation derived via averaging.
In particular we found that deviations of the first passage location $\bx^S$
orthogonal to the switching manifold are $O(\ee)$,
whereas deviations in a direction parallel to the switching manifold are $O(\sqrt{\ee})$,
as evident in panels B and C of Fig.~\ref{fig:checkSlide2}.
This is because the discontinuity in the equations along $x_1 = 0$
inhibits deviations in the $x_1$-direction.
Furthermore, for the relay control system with noise added purely to the control response,
the leading order terms of ${\rm Std} \left( t^S \big| \bx_0 \right)$
and ${\rm Std} \left( x_j^S \big| \bx_0 \right)$ for $j > 2$ vanish
because the noise effectively acts only in the $x_1$-direction.
Consequently these standard deviations are $O(\ee)$, Fig.~\ref{fig:checkSlide}.

\subsubsection*{Escaping dynamics}

We defined escaping dynamics as sections of solutions that lie within the strip, $\delta^- < x_2 < \delta^+$,
where $\delta^-$ and $\delta^+$ are suitably small, (\ref{eq:deltaMinusdeltaPlus}).
As shown in \S\ref{sec:ESCAPE},
the spatial and time scales for escaping dynamics are
$x_1 = O(\ee^{\frac{2}{3}})$,
$x_j = O(\ee^{\frac{1}{3}})$ for $j > 1$, and
$t = O(\ee^{\frac{1}{3}})$.
We derived the leading order component of the transitional PDF for (\ref{eq:sde}) for an escaping phase by assuming $x_1 > 0$,
imposing a reflecting boundary condition at $x_1 = 0$, and solving the corresponding Fokker-Planck equation.
The result is Knessl's solution (\ref{eq:Kn00}).
However, escaping phases make up only a small fraction of dynamics
over a full oscillation and have little effect on $t_{\rm osc}$.
Indeed our final approximations of
${\rm Diff} \left( t_{\rm osc} \right)$ and ${\rm Std} \left( t_{\rm osc} \right)$ in \S\ref{sec:COMB}
do not involve calculations relating to escaping.

\subsubsection*{The statistics of $t_{\rm osc}$ for relay control}

In \S\ref{sec:COMB} we combined the results
to approximate ${\rm Diff} \left( t_{\rm osc} \right)$ and ${\rm Std} \left( t_{\rm osc} \right)$
for the relay control system (\ref{eq:ABCDvalues3d})-(\ref{eq:relayControlSystem2}).
Fig.~\ref{fig:checkOsc} reveals that the approximations (\ref{eq:Diffthalfosc2}), (\ref{eq:Difftosc2}), (\ref{eq:Varthalfosc2}) and (\ref{eq:Vartosc2})
fit the results of Monte-Carlo simulations reasonably well.
In view of the complexity in evaluating these approximations,
a geometric understanding of the terms in these equations is arguably more useful than the approximations themselves.
Here we use the results to obtain four reasons why
the noise significantly reduces the average oscillation time for the relay control example
at relatively small values of $\ee$, as seen in Fig.~\ref{fig:manyPeriod}.

${\rm Diff}(t_{\frac{1}{2} {\rm osc}})$ is approximated by (\ref{eq:Diffthalfosc2})
as a sum of three terms that we have ordered by decreasing magnitude.
The first term, ${\rm Diff} \left( t^R \big| \bx_{\Gamma}^E \right)$,
is negative and represents the difference created by the noise causing solutions
to return to the switching manifold earlier, on average, than in the absence of noise.
By (\ref{eq:meantex2}), this term is proportional to the square of the inverse of the velocity of $\Gamma$ at $\bx^R$ in the $x_1$-direction.
For the relay control system the velocity is $-\frac{1}{\omega^2}$,
where $\omega = 5$, and for this reason the first term is relatively large.
The second term of (\ref{eq:Diffthalfosc2}),
$\rD_{\bx} t_{\dd}^S \left( \bx_{\Gamma}^M \right)^{\sf T} {\rm Diff} \left( \bx^M \right)$,
represents the difference created by sliding phases taking, on average, less time than they would without noise
due to sliding phases starting at points $\bx^M$ that are, on average,
deviated from $\bx_{\Gamma}^M$ in a particular direction.
The description of these two terms provides one reason why small noise significantly decreases $t_{\rm osc}$:
{\em Since $\Gamma$ slowly approaches the switching manifold along a path
that has a sharp angle relative to the switching manifold,
small noise tends to push solutions onto the switching manifold early
and at points deviated from $\bx_{\Gamma}^R$}.
Loosely speaking, the noise causes solutions to ``cut the corner'' at $\bx_{\Gamma}^R$.

The third term of (\ref{eq:Diffthalfosc2}), ${\rm Diff} \left( t^S \big| \bx_{\Gamma}^M \right)$,
is proportional to $\Lambda$ (\ref{eq:Lambda}).
For the relay control system, $\Lambda$ involves terms in the first column of $A$, such as $\omega^2$, which is relatively large.
This suggests that the noise-induced effect observed in Fig.~\ref{fig:checkExcur} is due in part to this term
which we interpret as the result of {\em noise pushing solutions slightly off the switching manifold,
and causing the nature of the vector field away from the switching manifold to influence dynamics}.


Equation (\ref{eq:Varthalfosc2}) approximates ${\rm Var}(t_{\frac{1}{2} {\rm osc}})$.
The second term of (\ref{eq:Varthalfosc2}) represents the variance in the times of the regular phases,
and the first [resp.~third] term of (\ref{eq:Varthalfosc2}) represents the variance in $t^S$ [resp.~$t^R$]
due to the variability in the points $\bx^M$.
Therefore deviations in $t_{\rm osc}$ are due primarily to
the variability in the first passage statistics of the regular phase.
For the relay control example, the deviations in these statistics are not as large
as one might expect because away from $x_1 = 0$ solutions rapidly contract onto a slow manifold\removableFootnote{
Using $D = e_1 e_1^{\sf T}$ in place of $D = B e_1^{\sf T}$ in (\ref{eq:ABCDvalues3d})
gives similar theoretical values for ${\rm Var}(t^R)$
and ${\rm Cov}(\bx^R)$.
}.
{\em Hence the slow-fast nature of the system
inhibits large deviations in $t_{\rm osc}$}.
For this reason ${\rm Std}(t_{\rm osc})$ is relatively small;
this constitutes a third reason for the nature of Fig.~\ref{fig:manyPeriod}.

Finally, in (\ref{eq:relayControlSystem2}) {\em the noise is added purely to the control response
causing the leading order contribution of the noise during stochastically perturbed sliding motion to vanish}.
Specifically, $M M^{\sf T} = 0$ where $M$ is the diffusion matrix in (\ref{eq:linearDiffusionApprox2})
-- the averaging approximation to the difference between stochastic solutions and the deterministic solution in the sliding phase.
Hence, again, for our particular example, ${\rm Std} \left( t_{\rm osc} \right)$ is less than
we would expect it to be in general.

\subsubsection*{Issues and future work}

The many approximations in our calculations combine to form discrepancies
between numerical results, obtained by Monte-Carlo simulations, and theoretical results, as evident in Fig.~\ref{fig:checkOsc}.
For instance, we used only the three largest terms in our expressions for
the statistics of $t_{\rm osc}$ and $t_{\frac{1}{2} {\rm osc}}$ for Fig.~\ref{fig:checkOsc}.
Our calculations for each of the three phases involve expansions in $\ee$
and approximations are obtained by truncating these expansions.
Consequently the accuracy of the approximations decreases with increasing values of $\ee$.
Calculations regarding escaping involve the assumption $x_2 = O(\ee^{\frac{1}{3}})$.
Thus, strictly speaking, the values $\delta^-$ and $\delta^+$ should be $O(\ee^{\frac{1}{3}})$,
but for simplicity we have set them as constants (\ref{eq:specificdeltas}).
Another source of error is that for the relay control example
the distance of the point $\bx_{\Gamma}^M$, at which the sliding phase begins,
to the boundary of the stable sliding region is approximately equal to $\frac{1}{\omega^2}$, see Appendix \ref{sec:DET}.
With $\omega = 5$ this distance is relatively small
causing inaccuracy in the analysis for the sliding phase\removableFootnote{
With $\ee = 0.0001$, $6.4\%$ of $\bx^M$ values were actually located outside the stable sliding region.
(Naturally this percentage increases with $\ee$ and decreases with $\omega$.)
}
because the calculations are singular in the limit that the distance of $\Gamma$ to the sliding boundary goes to zero.
Also, we have not attempted to compute the $O(\ee)$ term of ${\rm Std} \left( t_{\rm osc} \right)$
necessary to fairly compare this value to ${\rm Diff} \left( t_{\rm osc} \right)$.

It remains to study large deviations of periodic orbits with sliding segments \cite{HiMe13}.
For systems with discontinuous drift, the small noise asymptotics of large deviations
may be fundamentally different to that of smooth systems \cite{GrHe01}.
In addition, it remains to investigate the effects of noise on sliding bifurcations
at which a segment of sliding motion is created or destroyed.

\appendix
\section{Calculations for the relay control example in the absence of noise}
\label{sec:DET}
\setcounter{equation}{0}

With $\ee = 0$, (\ref{eq:relayControlSystem2}) is the piecewise-linear ODE system
\begin{equation}
\dot{\bX} = \left\{ \begin{array}{lc}
A \bX + B \;, & X_1 < 0 \\
A \bX - B \;, & X_1 > 0
\end{array} \right. \;,
\label{eq:relayControlSystem4}
\end{equation}
where $A$ (\ref{eq:ABCDvalues3d}) has eigenvalues
$-\lambda$ and $-\omega \zeta \pm {\rm i} |\omega| \sqrt{1 - \zeta^2}$.
For the parameter values (\ref{eq:paramValues}), $\lambda$ is relatively small,
thus solutions to each linear half-system of (\ref{eq:relayControlSystem4})
rapidly approach the eigenspace corresponding to the eigenvalue $-\lambda$.
The eigenvector of $A$ for $-\lambda$ is
\begin{equation}
v_{-\lambda} = \left[ 1,~2 \zeta \omega,~\omega^2 \right]^{\sf T} \;,
\end{equation}
and the equilibria of the left and right half-systems of (\ref{eq:relayControlSystem4}) are, respectively,
\begin{equation}
\bX^{*(L)} = \left[ \begin{array}{c}
\frac{1}{\lambda \omega^2} \\
-1 + \frac{2 \zeta}{\lambda \omega} + \frac{1}{\omega^2} \\
2 + \frac{2 \zeta}{\omega} + \frac{1}{\lambda}
\end{array} \right] \;, \qquad
\bX^{*(R)} = \left[ \begin{array}{c}
-\frac{1}{\lambda \omega^2} \\
1 - \frac{2 \zeta}{\lambda \omega} - \frac{1}{\omega^2} \\
-2 - \frac{2 \zeta}{\omega} - \frac{1}{\lambda}
\end{array} \right] \;.
\end{equation}
$\bX^{*(L)}$ and $\bX^{*(R)}$ are both {\em virtual} equilibria of (\ref{eq:relayControlSystem4}).
The weak stable manifold for each equilibrium is the line that passes through the equilibrium in the direction $v_{-\lambda}$.
These manifolds intersect $X_1 = 0$ at
\begin{equation}
\bX_{\rm int}^{(L)} = \left[
0,~-1+\frac{1}{\omega^2},~2+\frac{2 \zeta}{\omega}
\right]^{\sf T} \;, \qquad
\bX_{\rm int}^{(R)} = \left[
0,~1-\frac{1}{\omega^2},~-2-\frac{2 \zeta}{\omega}
\right]^{\sf T} \;.
\label{eq:XintLR}
\end{equation}
Consequently $\Gamma$ arrives at $X_1 = 0$ at points extremely close to
$\bX_{\rm int}^{(L)}$ and $\bX_{\rm int}^{(R)}$.
For the purposes of applying the coordinate change described in \S\ref{sub:COORDRCS},
it is appropriate to approximate the point with $X_3 > 0$ at which $\Gamma$ returns to $X_1 = 0$ by $\bX_{\rm int}^{(L)}$.

Stable sliding motion occurs on $X_1 = 0$
when $\dot{X_1} > 0$ for the left half-system of (\ref{eq:relayControlSystem4})
and $\dot{X_1} < 0$ for the right half-system of (\ref{eq:relayControlSystem4}).
By (\ref{eq:ABCDvalues3d}), stable sliding motion occurs on the strip
\begin{equation}
\left\{ (0,X_2,X_3)^{\sf T} ~\big|~ -1 < X_2 < 1 \right\} \;.
\label{eq:stableSlidingRegion}
\end{equation}
Sliding motion is specified by Filippov's solution \cite{Fi88,Fi60}, which yields
\begin{equation}
\left[ \begin{array}{c}
\dot{X}_2 \\ \dot{X}_3
\end{array} \right] =
\left[ \begin{array}{cc}
2 & 1 \\ -1 & 0
\end{array} \right]
\left[ \begin{array}{c}
X_2 \\ X_3
\end{array} \right] \;.
\label{eq:slidingDyns}
\end{equation}
Equation (\ref{eq:slidingDyns}) has the explicit solution
\begin{equation}
\left[ \begin{array}{c}
X_{{\dd},2}(t;\bX_0) \\
X_{{\dd},3}(t;\bX_0)
\end{array} \right] =
{\rm e}^t \left[ \begin{array}{c} 1 \\ -1 \end{array} \right] X_{0,2} +
{\rm e}^t \left(
t \left[ \begin{array}{c} 1 \\ -1 \end{array} \right] +
\left[ \begin{array}{c} 0 \\ 1 \end{array} \right] \right)
(X_{0,2}+X_{0,3}) \;,
\label{eq:slidingSoln}
\end{equation}
for any initial point $\bX_0 = (0,X_{0,2},X_{0,3})$ in the stable sliding region.

The upper sliding segment of $\Gamma$ ends at $X_2 = 1$.
With the approximation that the sliding segment starts at $\bX_{\rm int}^{(L)}$, (\ref{eq:XintLR}),
the deterministic sliding time, here call it $T$, is therefore determined by 
\begin{equation}
X_{{\dd},2} \left( T;\bX_{\rm int}^{(L)} \right) = 1 \;,
\end{equation}
and sliding ends at
\begin{equation}
(0,1,Z)^{\sf T} \;, {\rm ~where~}
Z \equiv X_{{\dd},3} \left( T;\bX_{\rm int}^{(L)} \right) \approx 2.561 \;.
\end{equation}

In the transformed coordinates (\ref{eq:Xtox}), the initial point for the sliding phase
and the end point for the excursion phase are, respectively,
\begin{eqnarray}
\bx_{\rm int}^{(L)} = P \bX_{\rm int}^{(L)} + Q &=&
\left[ \begin{array}{c}
0 \\
-2+\frac{1}{\omega^2} \\
2+\frac{2 \zeta}{\omega} - Z - \left( 2-\frac{1}{\omega^2} \right)
\frac{1}{Z+2}
\end{array} \right] \;, \\
\bx_{\rm int}^{(R)} = P \bX_{\rm int}^{(R)} + Q &=&
\left[ \begin{array}{c}
0 \\
-\frac{1}{\omega^2} \\
-2-\frac{2 \zeta}{\omega} - Z - \frac{1}{\omega^2(Z+2)}
\end{array} \right] \;. \label{eq:xRint}
\end{eqnarray}

\section{Calculations of the regular phase for relay control}
\label{sec:CALEX}
\setcounter{equation}{0}

Here we provide details of calculations for the relay control example
that were outlined in \S\ref{sub:REGRELAY}.

The deterministic solution to (\ref{eq:relayControlSystem5R}) is given by
\begin{equation}
\bx_{\dd}(t) = 
{\rm e}^{\mathcal{A} t} \left( \bx_0 + \mathcal{A}^{-1} \mathcal{B}^{(R)} \right) -
\mathcal{A}^{-1} \mathcal{B}^{(R)} \;,
\label{eq:bxdetrelay}
\end{equation}
where $\bx_0 = \bx_{\dd}(0)$ denotes the initial point.
Here we take $\bx_0 = \bx_{\Gamma}^E$ (the deterministic end point of the
previous escaping phase, refer to Fig.~\ref{fig:phaseSchem})
with which first passage to the switching manifold
occurs at $\bx_{\Gamma}^R \approx \bx_{\rm int}^{(R)}$ (\ref{eq:xRint}), see \S\ref{sub:COORDRCS}.

Through elementary use of (\ref{eq:calA})-(\ref{eq:calBLBRD}),
the coefficients in the PDE for $\mathcal{P}$ (\ref{eq:fpe3}) are found to be
\begin{eqnarray}
\phi_1^{(R)}(\bx_{\dd}^R) &=& x_{\Gamma,2}^R \;, \label{eq:phi1R} \\
\phi_2^{(R)}(\bx_{\dd}^R) &=& \frac{-1}{Z+2} x_{\Gamma,2}^R + x_{\Gamma,3}^R + Z + 2 \;, \label{eq:phi2R} \\
\phi_3^{(R)}(\bx_{\dd}^R) &=& \frac{-1}{(Z+2)^2} x_{\Gamma,2}^R + \frac{1}{Z+2} x_{\Gamma,3}^R \;, \label{eq:phi3R}
\end{eqnarray}
\begin{equation}
\frac{\partial \phi_1^{(R)}}{\partial x_2}(\bx_{\dd}^R) = 1 \;, \qquad
\frac{\partial \phi_1^{(R)}}{\partial x_3}(\bx_{\dd}^R) = 0 \;,
\end{equation}
\begin{equation}
\left( D D^{\sf T} \right)_{1,1} = 1 \;, \qquad
\left( D D^{\sf T} \right)_{2,1} = -2 \;, \qquad
\left( D D^{\sf T} \right)_{3,1} = \frac{Z}{Z+2} \;.
\end{equation}
Then by substituting
\begin{equation}
\bx_{\dd} \left( \sqrt{\ee} \tau + t_{\Gamma}^R \right) =
\bx_{\Gamma}^R + \sqrt{\ee} \phi^{(R)}(\bx_{\dd}^R) \tau + O(\ee) \;,
\end{equation}
with (\ref{eq:phi1R})-(\ref{eq:phi3R}) into
the expression for the free-space PDF (\ref{eq:pfs}),
we obtain an expression for $f^{(0)}$ by the absorbing boundary condition (\ref{eq:exBC0}).
Specifically
\begin{equation}
f^{(0)}(u_2,u_3,\tau) =
\frac{1}{(2 \pi)^{\frac{3}{2}} \sqrt{\det(K(t_{\Gamma}^R))}}
\,{\rm exp} \left( -\frac{1}{2} \chi^{\sf T} K(t_{\Gamma}^R)^{-1} \chi \right) \;,
{\rm ~where~}
\chi = \left[ \begin{array}{c}
-\phi_1^{(R)}(\bx_{\dd}^R) \tau \\
u_2 - \phi_2^{(R)}(\bx_{\dd}^R) \tau \\
u_3 - \phi_3^{(R)}(\bx_{\dd}^R) \tau
\end{array} \right] \;,
\label{eq:f0}
\end{equation}
which is used in (\ref{eq:calP02}) to obtain $\mathcal{P}^{(0)}$.
The function $g^{(1)}$ (which appears in the second term of the expression for $\mathcal{P}^{(1)}$ (\ref{eq:calP12}))
is determined from (\ref{eq:calPj}) and is given by
\begin{eqnarray}
g^{(1)}(u_2,u_3,\tau) &=&
- \frac{2}{\left( D D^{\sf T} \right)_{1,1}} \left( \frac{\partial \phi_1^{(R)}}{\partial x_2}(\bx_{\dd}^R) u_2
+ \frac{\partial \phi_1^{(R)}}{\partial x_3}(\bx_{\dd}^R) u_3 \right) f^{(0)}
- \frac{1}{\phi_1^{(R)}(\bx_{\dd}^R)} f^{(0)}_{\tau} \nonumber \\
&&+~\left( \frac{2 \left( D D^{\sf T} \right)_{2,1}}{\left( D D^{\sf T} \right)_{1,1}} - \frac{\mu_1}{\phi_1^{(R)}(\bx_{\dd}^R)}
\right) f^{(0)}_{u_2}
+ \left( \frac{2 \left( D D^{\sf T} \right)_{3,1}}{\left( D D^{\sf T} \right)_{1,1}} - \frac{\mu_2}{\phi_1^{(R)}(\bx_{\dd}^R)}
\right) f^{(0)}_{u_3} \;.
\label{eq:g1}
\end{eqnarray}

\subsubsection*{Calculation of $\mathbb{E}[t^R]$}

From (\ref{eq:meantex1}) we can write
\begin{equation}
\mathbb{E}[t^R] =
\int_0^\infty \int_0^\infty \int_{-\infty}^\infty \int_{-\infty}^\infty
p_f(\bx,t) \,dx_3 \,dx_2 \,dx_1 \,dt +
\int_0^\infty \int_0^\infty \int_{-\infty}^\infty \int_{-\infty}^\infty
\mathcal{P}(z,u_2,u_3,\tau) \,dx_3 \,dx_2 \,dx_1 \,dt \;.
\label{eq:meantex4}
\end{equation}
Using
$\Psi(s) \equiv -\frac{x_{{\dd},1} \left( s+t_{\Gamma}^R \right)}
{\sqrt{2 \kappa_{11} \left( s+t_{\Gamma}^R \right)}}$,
$\xi = \frac{x_1 - x_{{\dd},1}(t)}{\sqrt{2 \kappa_{11}(t)}}$
and $s = t - t_{\Gamma}^R$,
the first integral in (\ref{eq:meantex4}) is
\begin{eqnarray}
&& \int_0^\infty \int_0^\infty \int_{-\infty}^\infty \int_{-\infty}^\infty
p_f(\bx,t) \,dx_3 \,dx_2 \,dx_1 \,dt \nonumber \\
&=& \int_0^\infty \int_0^\infty \int_{-\infty}^\infty \int_{-\infty}^\infty
\frac{1}{(2 \pi \ee)^{\frac{3}{2}} \sqrt{\det(K(t))}}
{\rm e}^{-\frac{1}{2 \ee} (\bx - \bx_{\dd}(t))^{\sf T}
K(t)^{-1} (\bx - \bx_{\dd}(t))}
\,dx_3 \,dx_2 \,dx_1 \,dt \nonumber \\
&=& \int_0^\infty \int_0^\infty
\frac{1}{\sqrt{2 \pi \ee \kappa_{11}(t)}}
{\rm e}^{-\frac{(x_1 - x_{{\dd},1}(t))^2}{2 \ee \kappa_{11}(t)}}
\,dx_1 \,dt \nonumber \\
&=& \frac{1}{\sqrt{\pi \ee}} \int_{-t_{\Gamma}^R}^\infty \int_{\Psi(s)}^\infty
{\rm e}^{-\frac{\xi^2}{\ee}} \,d\xi \,ds \;.
\end{eqnarray}
Then reversing the order of integration and expanding $s = \Psi^{-1}(\xi)$
as a Taylor series centred at $\xi = 0$ produces
\begin{eqnarray}
\frac{1}{\sqrt{\pi \ee}} \int_{-t_{\Gamma}^R}^\infty \int_{\Psi(s)}^\infty
{\rm e}^{-\frac{\xi^2}{\ee}} \,d\xi \,ds
&=& \frac{1}{\sqrt{\pi \ee}} \int_{-\infty}^\infty
\left(
t - \frac{\sqrt{2 \kappa_{11}}}{\dot{x}_{{\dd},1}} \xi
+ \left( \frac{\dot{\kappa}_{11}}{\dot{x}_{{\dd},1}^2}
- \frac{\kappa_{11} \ddot{x}_{{\dd},1}}{\dot{x}_{{\dd},1}^3} \right) \xi^2 + O(\xi^3)
\right) \Bigg|_{t = t_{\Gamma}^R}
{\rm e}^{-\frac{\xi^2}{\ee}} \,d\xi \nonumber \\
&=& t_{\Gamma}^R
+ \frac{1}{2} \left( \frac{\dot{\kappa}_{11}}{\left( \phi_1^{(R)}(\bx_{\dd}^R) \right)^2}
+ \frac{\kappa_{11} \ddot{x}_{{\dd},1}}{\left( \phi_1^{(R)}(\bx_{\dd}^R) \right)^3} \right)
\Bigg|_{t = t_{\Gamma}^R} \ee + O(\ee^2) \;.
\label{eq:intPart1}
\end{eqnarray}
The second integral in (\ref{eq:meantex4}) is
\begin{eqnarray}
&& \int_0^\infty \int_0^\infty \int_{-\infty}^\infty \int_{-\infty}^\infty
\mathcal{P}(z,u_2,u_3,\tau) \,dx_3 \,dx_2 \,dx_1 \,dt \nonumber \\
&=& -\ee \int_{-\frac{t_{\Gamma}^R}{\sqrt{\ee}}}^\infty
\int_0^\infty \int_{-\infty}^\infty \int_{-\infty}^\infty
f^{(0)}(u_2,u_3,\tau) {\rm e}^{2 x_{\Gamma,2}^R z}
\,du_3 \,du_2 \,dz \,d\tau + O \left( \ee^{\frac{3}{2}} \right) \nonumber \\
&=& -\frac{\ee}{2 \left( \phi_1^{(R)}(\bx_{\dd}^R) \right)^2} + O \left( \ee^{\frac{3}{2}} \right) \;,
\label{eq:intPart2}
\end{eqnarray}
and the sum of (\ref{eq:intPart1}) and (\ref{eq:intPart2}) produces (\ref{eq:meantex2}).

\subsubsection*{Calculation of $\mathbb{E}[\bx^R]$}

Here we briefly describe the manner by which we evaluate $\mathbb{E}[\bx^R]$ numerically.

Equation (\ref{eq:meanxjex1}) gives\removableFootnote{
The $\frac{1}{\ee}$ arises from $z = \frac{x_1}{\ee}$.
}
\begin{equation}
\mathbb{E}[x_j^R] = \frac{\ee}{2} (D D^{\sf T})_{11}
\int_0^\infty \int_{-\infty}^\infty \int_{-\infty}^\infty
x_j^R \left( \frac{\partial p_f}{\partial x_1}(0,x_2,x_3,t) 
+ \frac{1}{\ee} \frac{\partial \mathcal{P}}{\partial z}(0,u_2,u_3,\tau) \right)
\,dx_2 \,dx_3 \,dt \;,
\end{equation}
for $j = 2,3$, and changing to the local variables (\ref{eq:regularScaling}) yields
\begin{eqnarray}
\mathbb{E}[x_j^R] &=& \frac{\ee^{\frac{5}{2}}}{2} (D D^{\sf T})_{11}
\int_{\frac{t_{\Gamma}^R}{\sqrt{\ee}}}^\infty \int_{-\infty}^\infty \int_{-\infty}^\infty
\left( \sqrt{\ee} u_j + x_{{\Gamma},j}^R \right)
\bigg( \frac{\partial p_f}{\partial x_1}
(0,\sqrt{\ee} u_2 + x_{{\Gamma},2}^R,\sqrt{\ee} u_3 + x_{{\Gamma},3}^R,
\sqrt{\ee} \tau + t_{\Gamma}^R) \nonumber \\
&&+~\frac{1}{\ee} \frac{\partial \mathcal{P}}{\partial z}
(0,u_2,u_3,\tau) \bigg)
\,du_2 \,du_3 \,d\tau \;.
\end{eqnarray}
Since $p_f$ is Gaussian with covariance matrix, $K(t)$,
it is straightforward to derive
\begin{eqnarray}
\frac{\partial p_f}{\partial x_1}(0,x_2,x_3,t) &=&
-\frac{1}{\ee \det(K)}
\Big( -(\kappa_{22} \kappa_{33} - \kappa_{23}^2) x_{\Gamma,1}^R
+ (\kappa_{13} \kappa_{23} - \kappa_{12} \kappa_{33}) (x_2 - x_{\Gamma,2}^R) \nonumber \\
&&+~(\kappa_{12} \kappa_{23} - \kappa_{13} \kappa_{22}) (x_3 - x_{\Gamma,3}^R) \Big)
p_f(0,x_2,x_3,t) \;.
\end{eqnarray}
We also have from (\ref{eq:exBC0})
\begin{eqnarray}
\frac{\partial \mathcal{P}}{\partial z}(0,u_2,u_3,\tau) &=&
\frac{1}{\ee^{\frac{3}{2}}}
\left( -\frac{2 \phi_1^{(R)}(\bx_{\dd}^R)}{(D D^{\sf T})_{11}} \left( f^{(0)} + \sqrt{\ee} f^{(1)} \right)
+ \sqrt{\ee} g^{(1)} + O(\ee) \right) \nonumber \\
&=& -2 x_{\Gamma,2}^R p_f|_{x_1 = 0} 
+ \frac{1}{\ee} g^{(1)} + O \left( \frac{1}{\sqrt{\ee}} \right) \;.
\end{eqnarray}
Finally we obtain
\begin{eqnarray}
\mathbb{E}[x_j^R] &=& \frac{\ee^{\frac{3}{2}}}{2}
\int_{\frac{t_{\Gamma}^R}{\sqrt{\ee}}}^\infty \int_{-\infty}^\infty \int_{-\infty}^\infty
\left( \sqrt{\ee} u_j + x_{{\Gamma},j}^R \right)
\bigg( -\frac{1}{\det(K)} \Big( -(\kappa_{22} \kappa_{33} - \kappa_{23}^2) x_{\Gamma,1}^R \nonumber \\
&&+~(\kappa_{13} \kappa_{23} - \kappa_{12} \kappa_{33}) (x_2 - x_{\Gamma,2}^R)
+ (\kappa_{12} \kappa_{23} - \kappa_{13} \kappa_{22}) (x_3 - x_{\Gamma,3}^R) \Big) \nonumber \\
&&\times~p_f \left( 0,\sqrt{\ee} u_2 + x_{{\Gamma},2}^R,\sqrt{\ee} u_3 + x_{{\Gamma},3}^R,
\sqrt{\ee} \tau + t_{\Gamma}^R \right) \nonumber \\
&&+~\frac{2 \kappa}{\alpha} p_f \big|_{x_1 = 0}
+ \frac{1}{\ee} \,g^{(1)}(u_2,u_3,\tau) \bigg)
\,du_2 \,du_3 \,d\tau + O \left( \ee^{\frac{3}{2}} \right) \;.
\label{eq:meanxjex2}
\end{eqnarray}
To produce the black lines in panels B and C of Fig.~\ref{fig:checkExcur} we have numerically evaluated
the leading order component of (\ref{eq:meanxjex2}), which is $O(\ee)$.

\section{Calculation of $\sigma$}
\label{sec:SIGMA}
\setcounter{equation}{0}

Here we derive the formula (\ref{eq:sigma2}):
\begin{equation}
\sigma \sigma^{\sf T} = \frac{(b_L-b_R) (b_L-b_R)^{\sf T}}{(a_L+a_R)^2} \;,
\label{eq:sigma4}
\end{equation}
where $\sigma$ appears in (\ref{eq:linearDiffusionApprox}).
This is achieved by employing a linear diffusion approximation to reduce the drift term,
$\big( F_0(z(t),\by_{\dd}(t)) - \Omega(\by_{\dd}(t)) \big) \,dt$,
of (\ref{eq:linearDiffusionApprox0}), to a diffusion term that
approximates this drift term and in the limit $\ee \to 0$ has an equivalent distribution.
This is possible because the evolution of $z(t)$ is fast relative to that of $\by_{\dd}(t)$.

Since we are taking the limit $\ee \to 0$,
we may neglect higher order terms
in the stochastic differential equation for $z(t)$, (\ref{eq:dhatx4}).
Furthermore, the vector noise term in (\ref{eq:dhatx4}) is equivalent to a scalar noise term
$\sqrt{\alpha} \,dW(t)$, where $\alpha = (D D^{\sf T})_{11}$.
It is convenient to further replace this term with simply $dW(t)$,
as the noise amplitude $\sqrt{\alpha}$ appears as only a multiplicative factor in the final result.
We let
\begin{equation}
r = \frac{t}{\ee} \;,
\end{equation}
represent the fast time-scale.
Then (\ref{eq:dhatx4}) becomes
\begin{equation}
dq(r;\by) = \left\{ \begin{array}{lc}
a_L(\by) \;, & q < 0 \\
-a_R(\by) \;, & q > 0
\end{array} \right\} \,dr + \,dW(r) \;,
\label{eq:dxcheck}
\end{equation}
where we have replaced $z$ with the symbol $q$ to indicate that
changes mentioned above have been made.
In (\ref{eq:dxcheck}) $\by$ is treated as a constant,
so (\ref{eq:dxcheck}) represents {\em Brownian motion with two-valued drift} \cite{KaSh91}.

In order to approximate the behaviour of the drift term,
$\big( F_0(z(t),\by_{\dd}(t)) - \Omega(\by_{\dd}(t)) \big) \,dt$, in distribution, we let
\begin{equation}
R(r,\by) = \mathbb{E} \left[
\left( F_0(q(\tilde{r}+r;\by),\by) - \Omega(\by) \right)
\left( F_0(q(\tilde{r};\by),\by) - \Omega(\by) \right)^{\sf T}
\right] \;.
\label{eq:R}
\end{equation}
For $r \ge 0$, $R(r,\by)$ denotes the autocovariance
of the function $F_0$ (\ref{eq:F0}) with (\ref{eq:dxcheck}).
In (\ref{eq:R}), we take $q(\tilde{r};\by)$ to be at steady-state
and thus $R(r,\by)$ is independent of the value of $\tilde{r}$.
By stochastic averaging theory \cite{FrWe12,PaSt08,MoCu11,Kh66b},
in the limit $\ee \to 0$ the drift term may be replaced by
the diffusion term $\sigma(\by_{\dd}(t)) \sqrt{\alpha \ee} \,dW(t)$, where
\begin{equation}
\sigma(\by) \sigma(\by)^{\sf T} =
2 \int_0^\infty R(r,\by) \,dr \;.
\label{eq:sigma}
\end{equation}
Below we derive (\ref{eq:sigma4}) by evaluating (\ref{eq:sigma}).

Let $p(q,r|q_0)$ denote the transitional PDF of (\ref{eq:dxcheck}) with $q(0) = q_0$.
When $a_L, a_R > 0$, (\ref{eq:dxcheck}) has the steady-state PDF
\begin{equation}
p_{\rm ss}(q) = \frac{2 a_L a_R}{a_L+a_R}
\left\{ \begin{array}{lc}
{\rm e}^{2 a_L q} \;, & q \le 0 \\
{\rm e}^{-2 a_R q} \;, & q \ge 0
\end{array} \right. \;.
\label{eq:pss2}
\end{equation}
Then, by (\ref{eq:sigma}) we can write
\begin{equation}
\sigma(\by) \sigma(\by)^{\sf T} =
2 \int_0^\infty \int_{-\infty}^\infty \int_{-\infty}^\infty
\big( F_0(q,\by) - \Omega(\by) \big)
\big( F_0(q_0,\by) - \Omega(\by) \big)^{\sf T}
p(q,r|q_0) p_{\rm ss}(q_0) \,dq \,dq_0 \,dr \;,
\end{equation}
where $F_0$ is given by (\ref{eq:F0}).
Since,
\begin{equation}
\mathbb{E} \left[ F_0(q,\by) - \Omega(\by) \right] \equiv 0 \;,
\end{equation}
it follows that
\begin{equation}
\sigma(\by) \sigma(\by)^{\sf T} =
2 \int_0^\infty \int_{-\infty}^\infty \int_{-\infty}^\infty
\big( F_0(q,\by) - \Omega(\by) \big)
\big( F_0(q_0,\by) - \Omega(\by) \big)^{\sf T}
\big( p(q,r|q_0) - p_{\rm ss}(q) \big)
p_{\rm ss}(q_0) \,dq \,dq_0 \,dr \;.
\end{equation}
By (\ref{eq:Omega}), (\ref{eq:meanFcal}) and (\ref{eq:F0}) we have
\begin{equation}
F_0(q,\by) - \Omega(\by) =
\left\{ \begin{array}{lc}
\frac{a_L (b_L-b_R)}{a_L+a_R} \;, & q < 0 \\
\frac{-a_R (b_L-b_R)}{a_L+a_R} \;, & q > 0
\end{array} \right. \;.
\end{equation}
Therefore we can write
\begin{eqnarray}
\sigma(\by) \sigma(\by)^{\sf T} &=&
2 \frac{(b_L-b_R)(b_L-b_R)^{\sf T}}{(a_L+a_R)^2} \Bigg(
a_L^2 \int_0^\infty \int_{-\infty}^0 \int_{-\infty}^0
\big( p(q,r|q_0) - p_{\rm ss}(q) \big)
p_{\rm ss}(q_0) \,dq \,dq_0 \,dr \nonumber \\
&&-~a_L a_R \int_0^\infty \int_{-\infty}^0 \int_0^\infty
\big( p(q,r|q_0) - p_{\rm ss}(q) \big)
p_{\rm ss}(q_0) \,dq \,dq_0 \,dr \nonumber \\
&&-~a_L a_R \int_0^\infty \int_0^\infty \int_{-\infty}^0
\big( p(q,r|q_0) - p_{\rm ss}(q) \big)
p_{\rm ss}(q_0) \,dq \,dq_0 \,dr \nonumber \\
&&+~a_R^2 \int_0^\infty \int_0^\infty \int_0^\infty
\big( p(q,r|q_0) - p_{\rm ss}(q) \big)
p_{\rm ss}(q_0) \,dq \,dq_0 \,dr \Bigg) \;.
\label{eq:sigma3}
\end{eqnarray}
Next we show that
\begin{equation}
\int_0^\infty p(q,r|q_0) - p_{\rm ss}(q) \,dr = \left\{ \begin{array}{lc}
\left( \frac{a_L^3+a_R^3}{a_L a_R (a_L+a_R)^2} +
\frac{2 a_R}{a_L+a_R} (q+q_0) \right)
{\rm e}^{2 a_L q} \\
+~\frac{1}{a_R} \left( {\rm e}^{-a_R(q-q_0)-a_R|q-q_0|} - {\rm e}^{-2 a_R q} \right)
\;, & q_0 \le 0,\, q \le 0 \\
\left( \frac{a_L^3+a_R^3}{a_L a_R (a_L+a_R)^2} -
\frac{2 a_L}{a_L+a_R} q + \frac{2 a_R}{a_L+a_R} q_0 \right)
{\rm e}^{-2 a_R q} \;, & q_0 \le 0,\, q \ge 0 \\
\left( \frac{a_L^3+a_R^3}{a_L a_R (a_L+a_R)^2} +
\frac{2 a_R}{a_L+a_R} q - \frac{2 a_L}{a_L+a_R} q_0 \right)
{\rm e}^{2 a_L q} \;, & q_0 \ge 0,\, q \le 0 \\
\left( \frac{a_L^3+a_R^3}{a_L a_R (a_L+a_R)^2} -
\frac{2 a_L}{a_L+a_R} (q+q_0) \right)
{\rm e}^{-2 a_R q} \\
+~\frac{1}{a_L} \left( {\rm e}^{a_L(q-q_0)-a_L|q-q_0|} - {\rm e}^{2 a_L q} \right)
\;, & q_0 \ge 0,\, q \ge 0
\end{array} \right. \;,
\label{eq:LapDiff}
\end{equation}
and from (\ref{eq:pss2}) and (\ref{eq:LapDiff}) straight-forward integration reveals that 
the integrals that appear in (\ref{eq:sigma3}) are given simply by 
\begin{equation}
\begin{split}
\int_0^\infty \int_{-\infty}^0 \int_{-\infty}^0
\big( p(q,r|q_0) - p_{\rm ss}(q) \big)
p_{\rm ss}(q_0) \,dq \,dq_0 \,dr &=
\frac{1}{2 (a_L+a_R)^2} \;, \\
\int_0^\infty \int_{-\infty}^0 \int_0^\infty
\big( p(q,r|q_0) - p_{\rm ss}(q) \big)
p_{\rm ss}(q_0) \,dq \,dq_0 \,dr &=
\frac{-1}{2 (a_L+a_R)^2} \;, \\
\int_0^\infty \int_0^\infty \int_{-\infty}^0
\big( p(q,r|q_0) - p_{\rm ss}(q) \big)
p_{\rm ss}(q_0) \,dq \,dq_0 \,dr &=
\frac{-1}{2 (a_L+a_R)^2} \;, \\
\int_0^\infty \int_0^\infty \int_0^\infty
\big( p(q,r|q_0) - p_{\rm ss}(q) \big)
p_{\rm ss}(q_0) \,dq \,dq_0 \,dr &=
\frac{1}{2 (a_L+a_R)^2} \;,
\end{split}
\label{eq:sigma3aux}
\end{equation}
with which we immediately arrive at the desired result (\ref{eq:sigma4}).

To prove (\ref{eq:LapDiff}), we first note that,
as shown in \cite{KaSh84}, $p(q,r|q_0)$ is given by
\begin{equation}
p(q,r|q_0) = \left\{ \begin{array}{lc}
2 {\rm e}^{2 a_L q} \int_0^\infty 
h(r,b,a_R) * h(r,b-q-q_0,a_L) \,db +
G(q,r,a_L|q_0) \;, & q_0 \le 0,\, q \le 0 \\
2 {\rm e}^{-2 a_R q} \int_0^\infty
h(r,b+q,a_R) * h(r,b-q_0,a_L) \,db \;, & q_0 \le 0,\, q \ge 0 \\
2 {\rm e}^{2 a_L q} \int_0^\infty
h(r,b+q_0,a_R) * h(r,b-q,a_L) \,db \;, & q_0 \ge 0,\, q \le 0 \\
2 {\rm e}^{-2 a_R q} \int_0^\infty
h(r,b+q+q_0,a_R) * h(r,b,a_L) \,db +
G(q,r,-a_R|q_0) \;, & q_0 \ge 0,\, q \ge 0
\end{array} \right. \;,
\label{eq:p}
\end{equation}
where
\begin{eqnarray}
h(r,q_0,\omega) &=& \frac{|q_0|}{\sqrt{2 \pi r^3}}
{\rm e}^{-\frac{(q_0 - \omega r)^2}{2 r}} \;, \label{eq:h} \\
G(q,r,\omega|q_0) &=& \frac{1}{\sqrt{2 \pi r}}
{\rm e}^{-\frac{(q-q_0-\omega r)^2}{2 r}} - {\rm e}^{-2 \omega q_0} 
\frac{1}{\sqrt{2 \pi r}}
{\rm e}^{-\frac{(q+q_0-\omega r)^2}{2 r}} \;, \label{eq:Gabsorb}
\end{eqnarray}
and $*$ denotes convolution with respect to $r$.
Here we derive (\ref{eq:LapDiff}) from (\ref{eq:p})-(\ref{eq:Gabsorb}) for $q_0, q \ge 0$.
The case $q_0 \ge 0$, $q \le 0$ is similar and
the remaining two cases follow by symmetry.

For $q_0, q \ge 0$, direct integration yields\removableFootnote{
$\int_0^\infty G(q,r,-a_R|q_0) \,dr =
\frac{1}{a_R} \left( {\rm e}^{-a_R(q-q_0)-a_R|q-q_0|} - {\rm e}^{-2 a_R q} \right)$
}
\begin{align}
&\int_0^\infty p(q,r|q_0) - p_{\rm ss}(q) \,dr =
\frac{1}{a_R} \left( {\rm e}^{-a_R(q-q_0)-a_R|q-q_0|} - {\rm e}^{-2 a_R q} \right) \nonumber \\
&+~\lim_{\nu \to 0^+}
\mathcal{L} \left( 2 {\rm e}^{-2 a_R q} \int_0^\infty
h(r,b+q+q_0,a_R) * h(r,b,a_L) \,db - p_{\rm ss}(q) \right) \;,
\end{align}
where
\begin{equation}
\mathcal{L}[f(r)] = \int_0^\infty {\rm e}^{-\nu r} f(r) \,dr \;.
\end{equation}
denotes a Laplace transform in $r$.
Next, we recall (\ref{eq:pss2}) and note that
\begin{equation}
\mathcal{L}[h(r,q_0,\omega)] = {\rm e}^{\omega q_0 - \sqrt{\omega^2 + 2 \nu} |q_0|} \;,
\end{equation}
to obtain
\begin{align}
&\int_0^\infty p(q,r|q_0) - p_{\rm ss}(q) \,dr
= \frac{1}{a_R} \left( {\rm e}^{-a_R(q-q_0)-a_R|q-q_0|} - {\rm e}^{-2 a_R q} \right) \nonumber \\
&+~2 {\rm e}^{-2 a_R q} \lim_{\nu \to 0^+}
\left( \frac{{\rm e}^{\left( a_R - \sqrt{a_R^2 + 2 \nu} \right)(q+q_0)}}
{-a_R + \sqrt{a_R^2 + 2 \nu} - a_L + \sqrt{a_L^2 + 2 \nu}} -
\frac{a_L a_R}{\nu (a_L+a_R)} \right) \;.
\end{align}
Finally, by substituting
$-a + \sqrt{a^2 + 2 \nu} = \frac{\nu}{a} - \frac{\nu^2}{2 a^3} + O(\nu^3)$,
with $a = a_L,a_R$ in the above equation,
terms involving $\frac{1}{\nu}$ vanish and we arrive at (\ref{eq:LapDiff}) for $q_0, q \ge 0$.

\section{Derivations of formulas for ${\rm Diff} \left( t^S \right)$ and ${\rm Var}(t^S)$}
\label{sec:TSL}
\setcounter{equation}{0}

In this section we derive (\ref{eq:Difftsl2}) and (\ref{eq:combvartsl}):
\begin{align}
{\rm Diff} \left( t^S \right) &= {\rm Diff} \left( t^S \big| \bx_{\Gamma}^M \right) 
+ \rD_{\bx} t_{\dd}^S \left( \bx_{\Gamma}^M \right)^{\sf T} {\rm Diff} \left( \bx^M \right)
+ \sum_{i=1}^N \sum_{j=1}^N \rD_{\bx}^2 t_{\dd}^S \left( \bx_{\Gamma}^M \right)_{i,j}
{\rm Cov} \left( x_{\Gamma}^M \right)_{i,j} + O \left( \ee^{\frac{3}{2}} \right) \;, \label{eq:DifftslA} \\
{\rm Var} \left( t^S \right) &= {\rm Var} \left( t^S \big| \bx_{\Gamma}^M \right)
+ \rD_{\bx} t_{\dd}^S \left( \bx_{\Gamma}^M \right)^{\sf T}
{\rm Cov} \left( \bx^M \right)
\rD_{\bx} t_{\dd}^S \left( \bx_{\Gamma}^M \right) + O \left( \ee^{\frac{3}{2}} \right) \;, \label{eq:combvartslA}
\end{align}
which express the leading order terms for ${\rm Diff} \left( t^S \right)$ and ${\rm Var}(t^S)$
in terms of conditioned quantities and may be evaluated using the results of \S\ref{sec:SLIDE}.
Analogous formulas in \S\ref{sec:COMB} relating to other components
of the stochastic dynamics may be derived in the same fashion\removableFootnote{
For example,
\begin{eqnarray}
{\rm Cov} \left( \bx^S \right)
&=& \int \left( \bx^S - \mathbb{E} \left[ \bx^S \right] \right)
\left( \bx^S - \mathbb{E} \left[ \bx^S \right] \right)^{\sf T} p \left( \bx^S \right) \,d\bx^S \nonumber \\
&=& \int \int \left( \bx^S - \mathbb{E} \left[ \bx^S \big| \bx^M \right]
+ \mathbb{E} \left[ \bx^S \big| \bx^M \right] - \mathbb{E} \left[ \bx^S \right] \right)
\left( \cdots \right)^{\sf T}
p \left( \bx^S \big| \bx^M \right) \,d\bx^S p \left( \bx^M \right) \,d\bx^M \nonumber \\
&=& \int \int \left( \bx^S - \mathbb{E} \left[ \bx^S \big| \bx^M \right]
+ \rD_{\bx} \bx_{\dd}^S \left( \bx_{\Gamma}^M \right) \left( \bx^M - \bx_{\Gamma}^M \right)
+ O(\ee) \right) \left( \cdots \right)^{\sf T} p \left( \bx^S \big| \bx^M \right)
\,d\bx^S p \left( \bx^M \right) \,d\bx^M \nonumber \\
&=& \int {\rm Cov} \left[ \bx^S \big| \bx^M \right] p \left( \bx^M \right) \,d\bx^M \nonumber \\
&&+~\rD_{\bx} \bx_{\dd}^S \left( \bx_{\Gamma}^M \right)
\int \left( \bx^M - \bx_{\Gamma}^M \right) \left( \bx^M - \bx_{\Gamma}^M \right)^{\sf T}
p \left( \bx^M \right) \,d\bx^M \rD_{\bx} \bx_{\dd}^S \left( \bx_{\Gamma}^M \right)^{\sf T}
+ O \left( \ee^{\frac{3}{2}} \right) \nonumber \\
&=& {\rm Cov} \left( \bx^S \big| \bx_{\Gamma}^M \right)
+ \rD_{\bx} \bx_{\dd}^S \left( \bx_{\Gamma}^M \right)
{\rm Cov} \left( \bx^M \right)
\rD_{\bx} \bx_{\dd}^S \left( \bx_{\Gamma}^M \right)^{\sf T} \;.
\end{eqnarray}
}.

First, by definition,
\begin{equation}
{\rm Diff} \left( t^S \right) \equiv \mathbb{E} \left[ t^S \right] - t_{\Gamma}^S
= \int t^S p \left( t^S \right) \,dt^S - t_{\Gamma}^S \;,
\end{equation}
where throughout this exposition $p(\cdot)$ denotes the PDF of the indicated variable.
Conditioning over the starting point $\bx^M$ gives
\begin{equation}
{\rm Diff} \left( t^S \right) =
\int t^S \int p \left( t^S \big| \bx^M \right) p \left( \bx^M \right) \,d\bx^M \,dt^S - t_{\Gamma}^S \;.
\end{equation}
By then reversing the order of integration and using
${\rm Diff} \left( t^S \big| \bx^M \right) \equiv
\mathbb{E} \left[ t^S \big| \bx^M \right] - t_{\dd}^S \left( \bx^M \right)$, we obtain
\begin{equation}
{\rm Diff} \left( t^S \right) =
\int \left( t_{\dd}^S \left( \bx^M \right) + {\rm Diff} \left( t^S \big| \bx^M \right) \right)
p \left( \bx^M \right) \,d\bx^M - t_{\Gamma}^S \;.
\label{eq:DifftslA2}
\end{equation}
By replacing $t_{\dd}^S \left( \bx^M \right)$ in (\ref{eq:DifftslA2})
with its Taylor series centred at the deterministic value $\bx^M = \bx_{\Gamma}^M$:
\begin{equation}
t_{\dd}^S \left( \bx^M \right) =
t_{\dd}^S \left( \bx_{\Gamma}^M \right) +
\rD_{\bx} t_{\dd}^S \left( \bx_{\Gamma}^M \right)^{\sf T} \left( \bx^M - \bx_{\Gamma}^M \right) +
\left( \bx^M - \bx_{\Gamma}^M \right)^{\sf T}
\rD_{\bx}^2 t_{\dd}^S \left( \bx_{\Gamma}^M \right)
\left( \bx^M - \bx_{\Gamma}^M \right) + O \left( \ee^{\frac{3}{2}} \right) \;,
\label{eq:tdetSTaylor}
\end{equation}
and evaluating the integral in (\ref{eq:DifftslA2}) we arrive at (\ref{eq:DifftslA}).
The error term in (\ref{eq:tdetSTaylor}) is $O \left( \ee^{\frac{3}{2}} \right)$
because $\bx^M - \bx_{\Gamma}^M = O \left( \sqrt{\ee} \right)$.

Second, to derive (\ref{eq:combvartslA}) we begin by writing
\begin{equation}
{\rm Var} \left( t^S \right)
= \int \left( t^S - \mathbb{E} \left[ t^S \right] \right)^2 p \left( t^S \right) \,dt^S \;.
\end{equation}
Conditioning over $\bx^M$ gives
\begin{equation}
{\rm Var} \left( t^S \right)
= \int \left( t^S - \mathbb{E} \left[ t^S \right] \right)^2 \int p\left( t^S \big| \bx^M \right)
p \left( \bx^M \right) \,d\bx^M \,dt^S \;.
\end{equation}
Reversing the order of integration and adding and subtracting $\mathbb{E} \left[ t^S \big| \bx^M \right]$ produces
\begin{equation}
{\rm Var} \left( t^S \right)
= \int \int  \left( t^S - \mathbb{E} \left[ t^S \big| \bx^M \right]
+ \mathbb{E} \left[ t^S \big| \bx^M \right] - \mathbb{E} \left[ t^S \right] \right)^2
p \left( t^S \big| \bx^M \right) \,dt^S p \left( \bx^M \right) \,d\bx^M \;.
\label{eq:combvartslA2}
\end{equation}
Since the mean values differ from their deterministic values by $O(\ee)$, we have
\begin{equation}
\mathbb{E} \left[ t^S \big| \bx^M \right] - \mathbb{E} \left[ t^S \right] =
t_{\Gamma}^S \left( \bx^M \right) - t_{\Gamma}^S + O(\ee) \;.
\label{eq:meanSubstraction}
\end{equation}
By substituting (\ref{eq:tdetSTaylor}) and (\ref{eq:meanSubstraction}) into (\ref{eq:combvartslA2}),
and noting $t_{\dd}^S \left( \bx_{\Gamma}^M \right) = t_{\Gamma}^S$,
we obtain
\begin{equation}
{\rm Var} \left( t^S \right)
= \int \int  \left( t^S - \mathbb{E} \left[ t^S \big| \bx^M \right]
+ \rD_{\bx} t_{\dd}^S \left( \bx_{\Gamma}^M \right)^{\sf T} \left( \bx^M - \bx_{\Gamma}^M \right) + O(\ee) \right)^2
p \left( t^S \big| \bx^M \right) \,dt^S p \left( \bx^M \right) \,d\bx^M \;.
\label{eq:combvartslA3}
\end{equation}
Finally, by expanding the square in (\ref{eq:combvartslA3}) and evaluating the
double integral we arrive at (\ref{eq:combvartslA}).

\end{document}